%% file: main.tex
\documentclass[table,a4paper,reqno]{amsart}
\usepackage{graphicx} 
\usepackage{adjustbox}
\usepackage{comment}
\usepackage{amsmath}
\usepackage{amsthm}
\usepackage{amssymb}
\usepackage{mathrsfs}
\usepackage{adjustbox}
\usepackage{diagbox}
\usepackage{mathtools}
\usepackage{tikz, pgfplots}
\usetikzlibrary{positioning}
\usepackage{rotating}
\usepackage{lscape}
\usepackage{multicol}
\usepackage{multirow}
\usepackage[thinlines]{easytable}
\usepackage{cite}
\usepackage{multicol}
\usepackage{esvect}
\usepackage{mathdots}
\usepackage{xcolor}
\usepackage{young}
\usepackage{ytableau}

\numberwithin{equation}{section}
\newtheorem{theorem}{Theorem}[section]
\newtheorem*{remark}{Remark}
\newtheorem{lemma}[theorem]{Lemma}
\newtheorem{proposition}[theorem]{Proposition}
\newtheorem{corollary}[theorem]{Corollary}
\theoremstyle{definition}
\newtheorem{definition}[theorem]{Definition}

\begin{document}
\title[The Haar State values of monomials on $\mathcal{O}(U_q(3))$]{The Haar State values of monomials and a method to pick orthonormal bases on $\mathcal{O}(U_q(3))$}
\author{Ting Lu}
\address{Ting Lu, Institute for Advanced Study in Mathematics, Harbin Institute of Technology, Harbin, Heilongjiang Province, China}
\email{ting\_lu@hit.edu.cn}

\begin{abstract}
    In this paper, we investigate the evaluation problem of the Haar state on the quantum group $\mathcal{O}(U_q(n))$~\cite{faddeev1988quantization,klimyk2012quantum} ($n\ge 3$) which is a $q$-deformation of the Haar measure on the Lie group $U(n)$. The relation between the Haar state values of monomials on $\mathcal{O}(U_q(n))$ is studied. On $\mathcal{O}(U_q(3))$, the Haar state values of monomials are explicitly computed and these values are expressed as a finite summation of rational polynomials in $q$. As an application, we compute the Gram matrices of the irreducible co-representations of $\mathcal{O}(U_q(3))$ which is essential to the method of constructing orthonormal bases on $\mathcal{O}(U_q(3))$ proposed by Noumi, Yamada, and Mimachi~\cite{noumi1993finite}. New connections between the Haar state values of monomials and basic hypergeometric multi-summations are found during our computation.
\end{abstract}
\maketitle

\setcounter{tocdepth}{2}

\tableofcontents

\input{Chapters/New_intro}

\input{Chapters/Chapter_2}

\input{Chapters/Chapter_3}

\input{Chapters/Chapter_4}

\input{Chapters/Chapter_5}

\section*{Acknowledgment}
This paper contains the work of the author during his PhD study at Texas A\&M University. The author sincerely appreciates all the helps and advice from his PhD advisor Prof. Jeffrey Kuan from Texas A\&M University and PhD committee member Prof. Michael Brannan from the University of Waterloo. The author also feels grateful for the helpful discussion with Prof. Simeng Wang at Harbin Institute of Technology during the writing of this paper.
\bibliographystyle{plain}
\bibliography{reference}

\appendix
\input{Chapters/Appendix_A}
\newpage
\input{Chapters/Appendix_B}
\newpage
\input{Chapters/Appendix_C}

\end{document}

%% file: Chapters/New_intro.tex
\section{Introduction}
 The existence and uniqueness of the Haar measure, a translation invariant probability measure, is one of the fundamental properties of compact topological groups. Its applications ranges from pure math to various scientific fields including representation theory, random matrix theory, theoretical physics, mathematical physics, statistics and many other fields. In this paper, we present the first systematic investigation on the evaluation problem of the Haar state $h$ on the quantum group $\mathcal{O}(U_q(n))$~\cite{faddeev1988quantization,klimyk2012quantum} ($n\ge 3$) which can be considered as a $q$-deformation of the Haar measure on $U(n)$. A brief review of the algebraic property of $\mathcal{O}(U_q(n))$ and its Haar state $h$ can be found in Section 2. All results in this paper is directly applicable to Woronowicz's compact quantum group $SU_q(n)$~\cite{woronowicz1988tannaka,woronowicz1987compact}.
 
The interest of evaluating Haar integrals stems from studies in theoretical physics. Let $U_{i,j}$'s ($1\le i,j\le n$) be coordinate functions on a matrix Lie group $G$ and $\overline{U_{i,j}}$ the conjugate of $U_{i,j}$. Haar integrals in the following form 
\begin{equation}\label{eq:1.1}
    \int_{U\in G}U_{i_1,j_1}\cdots U_{i_m,j_m}\overline{U_{l_1,k_1}\cdots U_{l_m,k_m}}\text{d}U
\end{equation}
appear in lattice gauge theory and theoretical physicists become the first group of researchers to investigate the evaluation of (\ref{eq:1.1}) on compact matrix Lie groups. Related notable works include but not limit to Weingarten~\cite{weingarten1978asymptotic}, Creutz~\cite{creutz1978invariant}, and Brezin and Gross~\cite{brezin1980external}. The evaluation of (\ref{eq:1.1}) is also related to the method of moments in the study of random matrix theory. In the work of Diaconis and Shahshahani~\cite{diaconis1994eigenvalues}, the evaluation of a Haar integral similar to (\ref{eq:1.1}) is the essential step in understanding the distribution of eigenvalues of random matrices on compact matrix Lie groups. Following their work there is a surge of interest in searching for a systematic method to evaluate Haar integrals on compact groups. This quest came to an end when the Weingarten calculus is proposed in the work of Collins~\cite{collins2003moments} and Collins and {\'S}niady~\cite{collins2006integration}. 

The key point of the Weingarten calculus is that the evaluation of Haar integral in (\ref{eq:1.1}) is converted to computing certain matrices associated with the invariant tensor space under the action of $G$. Readers may find more details on the Weingarten calculus in the work of Collins, Matsumoto, and Novak~\cite{collins2022weingarten}. Since its emergence, the Weingarten calculus has been applied to numerous studies across multiple fields and proved to be a fruitful theory as summarized by Collins~\cite{collins2022moment}. Especially, the theory has been successfully extended to free compact quantum groups~\cite{wang1995free,banica2007integration,banica2007permutationintegration,banica2008integration}. 

However, on $\mathcal{O}(U_q(n))$, the evaluation provided by the Weingarten calculus is not as efficient as the classical case since the Haar state is non-tracial ($h(xy)\ne h(yx)$). In particular, the consequence of $h$ being non-tracial cannot be neglected where $q$ approaches $0$. Our strategy to this issue is to break the integrand in (\ref{eq:1.1}) into smaller pieces and consider the evaluation problem of
\begin{equation}\label{eq:1.5}
    h(x_{i_1,j_1}\cdots x_{i_m,j_m}\cdot \text{det}_q^{-a_m}).
\end{equation}
We will call the integrand in (\ref{eq:1.5}) a monomial on $\mathcal{O}(U_q(n))$. 

On $SU_q(2)$~\cite{sl1987twisted} which is a special case of $\mathcal{O}(U_q(2))$, the evaluation of (\ref{eq:1.5}) is studied by Woronowicz~\cite{sl1987twisted} and Masuda \textit{et al.}~\cite{masuda1991representations} independently. Both of their works rely on the fact that every monomial on $SU_q(2)$ with a non-zero Haar state value can be written as a linear combination of a set of special monomials and those special monomials are mutually commutative. On $\mathcal{O}(U_q(n))$ ($n\ge 3$), such a set of special monomials can be found but they are not mutually commutative. Hence, their approaches on $SU_q(2)$ cannot be applied on higher rank $\mathcal{O}(U_q(n))$ directly. 

In Section 3, we discuss the evaluation of (\ref{eq:1.5}) on general $\mathcal{O}(U_q(n))$. A special set of monomials that serves like a basis of the support of the Haar state is chosen in Section 3.1. The method of finding linear relations between the Haar state values of our chosen monomials is given in Section 3.2. The techniques of Woronowicz and Masuda \textit{et al.} are combined and we use these new techniques to construct the Source matrix in Section 3.4.  Although a complete solution is still in searching, the Source matrix allows us to find $h\left((\prod_{i=1}^nx_{i,n+1-i})^m\cdot \text{det}_q^{-m}\right)$. The importance of $h\left((\prod_{i=1}^nx_{i,n+1-i})^m\cdot \text{det}_q^{-m}\right)$ in the evaluation of (\ref{eq:1.5}) will be explained in Section 4. 

In Section 4, we give a complete evaluation of (\ref{eq:1.5}) on $\mathcal{O}(U_q(3))$ (Theorem 4.1). Our approach rely on the invariant property of the Haar state under the action of $\mathcal{U}_q(gl_3)$ on $\mathcal{O}(U_q(3))$. This approach is inspired by the work of Noumi, Yamada, and Mimachi~\cite{noumi1993finite} in which the evaluation of (\ref{eq:1.1}) on the homogeneous quantum space $\mathcal{O}(U_q(n))\backslash \mathcal{O}(U_q(n-1))$ is given. In Corollary 4.2, the result in Theorem 4.1 is extended to general $\mathcal{O}(U_q(n))$.

Just like the classical case, $\mathcal{O}(U_q(n))$ is a pre-Hilbert space with the inner product $\langle a,b\rangle_L=h(a^*b)$ (or $\langle a,b\rangle_R=h(ab^*)$) as proved by Noumi, Yamada, and Mimachi~\cite{noumi1993finite}. They further proposed a method to construct orthonormal bases on $\mathcal{O}(U_q(n))$ given that we can find orthonormal bases of irreducible corepresentations of $\mathcal{O}(U_q(n))$. In Section 5, we apply our result to compute the Gram matrices of irreducible corepresentations of $\mathcal{O}(U_q(3))$(Theorem 5.8 and 5.9). With these Gram matrices, one can construct orthonormal bases on the irreducible corepresentations by the Gram-Schmidt process and then an orthonormal basis of $\mathcal{O}(U_q(3))$ is obtained by the method of Noumi, Yamada, and Mimachi~\cite{noumi1993finite}.

The evaluations of Haar integrals on compact (quantum) groups are closely related to combinatorial problems. This connection has been studied in the work of Brouwer and Beenakker~\cite{brouwer1996diagrammatic}, Collins, Gurau, and Lionni~\cite{collins2023tensor}, and Magee and Puder~\cite{magee2019matrix}. Moreover, Kumar and Landsberg~\cite{kumar2015connections} show that the evaluation of a special case of (\ref{eq:1.5}) on $SU(n)$ is related to the Alon-Tarsi conjecture. In our opinion, it would be interesting to explore a combinatorial explanation for our result (\ref{eq:4.73}) and (\ref{eq:5.35+}).

The study of $\mathcal{O}(SU_q(n))$ is also connected with basic hypergeometric series~\cite{masuda1991representations,koelink1989clebsch,vaksman1990algebra,bergeron2019suq}. Basic hypergeometric series appearing in these studies are single summations. In our computations, multi-summations of basic hypergeometric terms appear as the inner products between basis vectors of irreducible corepresentions. By comparing our computations with the results of Noumi, Yamada, and Mimachi~\cite{noumi1993finite}, we achieve an identity of basic hypergeometric multi-summation in Proposition 5.3. A more complicate basic hypergeometric multi-summation identity can be derived from Proposition 5.6. We hope our work can provide a quantum group perspective to the study of basic hypergeometric multi-summations.

%% file: Chapters/Chapter_2.tex
\section{Preliminaries on $\mathcal{O}(GL_q(n))$}
In this paper, we work in the framework of Faddeev, Reshetikhin, and Takhtajan~\cite{faddeev1988quantization,dijkhuizen1994cqg}. $\mathcal{O}(U_q(n))$ is the "compact" real form of $\mathcal{O}(GL_q(n))$. The existence and uniqueness of the Haar state on $\mathcal{O}(GL_q(n))$ relies on the fact that $\mathcal{O}(GL_q(n))$ is a co-semisimple Hopf algebra~\cite{sweedler1969hopf}. In particular, the algebra of coordinate functions on $U(n)$ is co-semisimple and the unique Haar state is the integral with respect to the normalized Haar measure. Readers can find a detailed study of the algebraic properties and representation theory on $\mathcal{O}(GL_q(n))$ from Noumi, Yamada, and Mimachi~\cite{noumi1993finite}. Here we give a brief review of the Hopf algebra $\mathcal{O}(GL_q(n))$, the relation between $\mathcal{U}_q(gl_n)$ and $\mathcal{O}(GL_q(n))$, the Peter-Weyl decomposition of $\mathcal{O}(GL_q(n))$, and the Haar state on $\mathcal{O}(GL_q(n))$. Then, we introduce the method of picking orthonormal basis on $\mathcal{O}(U_q(n))$ proposed by Noumi, Yamada, and Mimachi~\cite{noumi1993finite}.

\input{Chapters/Chapter_2/algebraic_properties}

\input{Chapters/Chapter_2/action_on_the_Hopf}

\input{Chapters/Chapter_2/finite_dim_rep}

\input{Chapters/Chapter_2/haar_state}

%% file: Chapters/Chapter_2/algebraic_properties.tex
\subsection{Algebraic properties of $\mathcal{O}(GL_q(n))$}
Before we introduce $\mathcal{O}(GL_q(n))$, we define the bi-algebra $\mathcal{O}(Mat_q(n))$. Denote the generators of $\mathcal{O}(Mat_q(n))$ as $\{x_{i,j}\}_{i,j=1}^n$. We will assume that these generators are put in a "generator matrix" in the following way.
\begin{align}
    \begin{bmatrix}
    x_{1,1}&\cdots&x_{1,n}\\
    \vdots&\cdots&\vdots\\  
    x_{n,1}&\cdots&x_{n,n}
    \end{bmatrix}\tag{2.0}\label{mat:generator_n}
\end{align}
The bialgebra $\mathcal{O}(Mat_q(n))$ is generated by the commutation relations ($i<l$, $j<k$):
\begin{align}\stepcounter{equation}
    \tag{2.1.a}x_{i,j}x_{l,j}=q\cdot x_{l,j}x_{i,j}\ \ \ \ \ \ \ \ \ x_{i,j}x_{i,k}=q\cdot x_{i,k}x_{i,j},&\label{eq:2.1.a} \\
    x_{ik}x_{lj}=x_{lj}x_{ik},&\tag{2.1.b}\label{eq:2.1.b}\\
    x_{i,j}x_{l,k}=x_{l,k}x_{i,j}+(q-q^{-1})x_{ik}x_{lj}.&\tag{2.1.c}\label{eq:2.1.c}
\end{align}
The co-multiplication $\Delta:\mathcal{O}(Mat_q(n))\rightarrow \mathcal{O}(Mat_q(n))\otimes\mathcal{O}(Mat_q(n))$ and co-unit $\epsilon: \mathcal{O}(Mat_q(n))\rightarrow\mathbb{C}$ is defined by
\begin{align}
    \Delta(x_{i,j})=\sum_{k=1}^nx_{i,k}\otimes x_{k,j}\quad\text{ and }\quad\varepsilon(x_{i,j})=\delta_{i,j}.\label{eq:2.2}
\end{align}
Notice that $\Delta$ and $\epsilon$ are $\mathbb{C}$-algebra homomorphisms. If we set $q=1$, then we get the bi-algebra of coordinate functions on $n\times n$ matrices.

Let $I=\{i_1,i_2,\cdots,i_r\}$ ($i_1<i_2<\cdots<i_r$) and $J=\{j_1,j_2,\cdots,j_r\}$ ($j_1<j_2<\cdots<j_r$) be two subsets of $\{1,2,\cdots,n\}$. The quantum $r-$minor determinant $\xi_J^I$ can be defined in two equivalent forms
\begin{align}\stepcounter{equation}
    \xi_J^I&=\sum_{\tau\in\mathfrak{S}_r}(-q)^{l(\tau)}x_{i_1,j_{\tau(1)}}\cdots x_{i_r,j_{\tau(r)}}\tag{2.3.a}\label{eq:2.3.a}\\
    \xi_J^I&=\sum_{\tau\in\mathfrak{S}_r}(-q)^{l(\tau)}x_{i_{\tau(1)},j_1}\cdots x_{i_{\tau(r)},j_r}\tag{2.3.b}\label{eq:2.3.b}
\end{align}
where $\mathfrak{S}_r$ is the permutation group on $r$ letters and $l(\tau)$ is the number of inversions of $\tau$. Noumi, Yamada, and Mimachi~\cite{noumi1993finite} showed that $\xi_J^I$'s are actually the matrix coefficients of the alternating tensor representation of $\mathcal{O}(Mat_q(n))$ of degree $r$. The quantum determinant $D_q=\xi_{\{1,\cdots,n\}}^{\{1,\cdots,n\}}$ is defined by
\begin{align}\stepcounter{equation}
    D_q&=\sum_{\sigma\in\mathfrak{S}_n}(-q)^{l(\sigma)}x_{\sigma(1),1}\cdots x_{\sigma(n),n}\tag{2.4.a}\label{2.4.a}\\
    &=\sum_{\tau\in\mathfrak{S}_n}(-q)^{l(\tau)}x_{1,\tau(1)}\cdots x_{n,\tau(n)}.\tag{2.4.b}\label{eq:2.4.b}
\end{align}
We remark that since $\xi_J^I$'s are matrix coefficients they satisfy
\begin{align}\stepcounter{equation}
    \Delta(\xi_J^I)=\sum_{|K|=r}\xi_K^I\otimes \xi_J^K.\tag{2.5.a}\label{eq:2.5.a}
\end{align}
Especially we have
\begin{align}
    \Delta(D_q)=D_q\otimes D_q. \tag{2.5.b}\label{eq:2.5.b}
\end{align}
Let $\widehat{k}=\{1,\cdots,k-1,k+1,\cdots,n\}$. By Noumi, Yamada, and Mimachi~\cite{noumi1993finite}, we have
\begin{align}\stepcounter{equation}
    \sum_{k=1}^n(-q)^{i-k}\xi_{\widehat{i}}^{\widehat{k}}x_{k,j}=\delta_{i,j}D_q,\tag{2.6.a}\label{eq:2.6.a}\\
    \sum_{k=1}^n(-q)^{k-j}x_{i,k}\xi^{\widehat{j}}_{\widehat{k}}=\delta_{i,j}D_q,\tag{2.6.b}\label{eq:2.6.b}\\
    \sum_{k=1}^n(-q)^{k-j}x_{k,i}\xi_{\widehat{j}}^{\widehat{k}}=\delta_{i,j}D_q,\tag{2.6.c}\label{eq:2.6.c}\\
    \sum_{k=1}^n(-q)^{i-k}\xi^{\widehat{i}}_{\widehat{k}}x_{j,k}=\delta_{i,j}D_q.\tag{2.6.d}\label{eq:2.6.d}
\end{align}
The Hopf algebra $\mathcal{O}(GL_q(n))$ is defined by adjoining an extra generator $\text{det}_q^{-1}$, the inverse of $D_q$, to $\mathcal{O}(Mat_q(n))$ satisfying
\begin{align}\stepcounter{equation}
    &D_q\cdot \text{det}_q^{-1}=\text{det}_q^{-1}\cdot D_q=1,\tag{2.7.a}\label{eq:2.7.a}\\
    & \text{det}_q^{-1}\cdot x_{i,j}=x_{i,j}\cdot \text{det}_q^{-1}.\tag{2.7.b}\label{eq:2.7.b}
\end{align}
If we identify $D_q$ with $1$, we get the Hopf algebra $\mathcal{O}(SL_q(n))$. Noumi, Yamada, and Mimachi~\cite{noumi1993finite} showed that the center of $\mathcal{O}(Mat_q(n))$ is generated by $D_q$ and $\text{det}_q^{-1}$ when $q$ is not a root of unity. The co-multiplication and co-unit of $\mathcal{O}(Mat_q(n))$ can be extended to $\mathcal{O}(GL_q(n))$ by
\begin{align}
    \Delta(\text{det}_q^{-1})=\text{det}_q^{-1}\otimes\text{det}_q^{-1}\quad \text{ and }\quad\varepsilon(\text{det}_q^{-1})=1.\label{eq:2.8}
\end{align}
The antipod $S:\mathcal{O}(GL_q(n))\rightarrow \mathcal{O}(GL_q(n))$ is defined by
\begin{align}
    S(x_{i,j})=(-q)^{i-j}\xi_{\widehat{i}}^{\widehat{j}}\cdot \text{det}_q^{-1}.\label{eq:2.9}
\end{align}
When $q\in\mathbb{R}$, we can define a conjugate linear anti-homomorphism $*:\mathcal{O}(GL_q(n))\rightarrow \mathcal{O}(GL_q(n))$ by
\begin{align}
    x_{i,j}^*=S(x_{j,i})=(-q)^{j-i}\xi_{\widehat{j}}^{\widehat{i}}\cdot \text{det}_q^{-1}.\label{eq:2.10}
\end{align}
With this $*$ map, $\mathcal{O}(GL_q(n))$ becomes a Hopf $*$-algebra and $(\mathcal{O}(GL_q(n)),*)$ is denoted as $\mathcal{O}(U_q(n))$, the "compact" real form of $\mathcal{O}(GL_q(n))$. If we consider Woronowicz's compact quantum group $SU_q(n)$ as a $*$-Hopf algebra, then it is $(\mathcal{O}(SL_q(n)),*)$.

Let $I^c$ denote the complement of $I$. Define
\begin{equation*}
    \text{sgn}_q(I;J)=\left\{\begin{array}{lcr}
        0 & \text{if} & I\cap J\neq\varnothing, \\
        (-q)^{l(I;J)} & \text{if} & I\cap J=\varnothing,
    \end{array}\right.
\end{equation*}
where $l(I;J)$ is the number of pair $(i,j), i\in I, j\in J$ such that $i>j$. Noumi, Yamada, and Mimachi showed that
\begin{align}\stepcounter{equation}
    \left(\xi_J^I\right)^*=S\left(\xi_I^J\right)=\frac{\text{sgn}_q(J;J^c)}{\text{sgn}_q(I;I^c)}\xi_{J^c}^{I^c}\cdot \text{det}_q^{-1}\tag{2.11.a}\label{eq:2.11.a}
\end{align}
and in particular
\begin{align}
    \left(D_q\right)^*=\text{det}_{q}^{-1}\text{ and }\left(\text{det}_{q}^{-1}\right)^*=D_q.\tag{2.11.b}\label{eq:2.11.b}
\end{align}

%% file: Chapters/Chapter_2/action_on_the_Hopf.tex
\subsection{The action of $\mathcal{U}_q(gl_n)$ on $\mathcal{O}(GL_q(n))$}
In this subsection, we review the relation between the quantum universal enveloping algebra $\mathcal{U}_q(gl_n)$~\cite{drinfeld1986quantum,jimbo1985aq} and the quantum group $\mathcal{O}(GL_q(n))$.

Let $F_n$ be the free $\mathbb{Z}$-module of rank $n$ with basis $\{\epsilon_1,\cdots,\epsilon_n\}$ with the symmetric bilinear form $\langle\cdot,\cdot\rangle$ such that $\langle\epsilon_i,\epsilon_j\rangle=\delta_{i,j}$. Elements of $F_n$ are called \textit{integral weights}.

The generators of $\mathcal{U}_q(gl_n)$ are denoted as $e_k,f_k(1\le k<n)$ and $q^\lambda(\lambda\in\frac{1}{2}F_n)$. For a detailed description of the algebraic structure of $\mathcal{U}_q(gl_n)$, we refer to \cite{noumi1993finite}. The Hopf algebra structure on $\mathcal{U}_q(gl_n)$ is defined by
\begin{align}\stepcounter{equation}
    &\Delta(q^\lambda)=q^\lambda\otimes q^\lambda,\quad \varepsilon(q^\lambda)=1,\quad S(q^\lambda)=q^{-\lambda},\tag{2.12.a}\label{eq:2.13.a}\\
    &\Delta(e_k)=e_k\otimes q^{-(\epsilon_k-\epsilon_{k+1})/2}+q^{(\epsilon_k-\epsilon_{k+1})/2}\otimes e_k,\tag{2.12.b}\label{eq:2.13.b}\\
    &\varepsilon(e_k)=0,\quad S(e_k)=-q^{-1}e_k,\notag\\
    &\Delta(f_k)=f_k\otimes q^{-(\epsilon_k-\epsilon_{k+1})/2}+q^{(\epsilon_k-\epsilon_{k+1})/2}\otimes f_k,\tag{2.12.c}\label{eq:2.13.c}\\
    &\varepsilon(f_k)=0,\quad S(f_k)=-qf_k.\notag
\end{align}
Let $\mathcal{A}$ and $\mathcal{B}$ be two Hopf algebras. By a \textit{pairing of Hopf algebras} we mean a bilinear form $(\cdot,\cdot):\mathcal{A}\times \mathcal{B}\rightarrow\mathbb{C}$ satisfying
\begin{align}\stepcounter{equation}
    (a,\varphi\psi)=(\Delta_\mathcal{A}(a),\varphi\otimes\psi)\quad\text{ and }\quad(a,1)=\varepsilon_\mathcal{A}(a)&,\tag{2.13.a}\label{eq:2.14.a}\\
    (ab,\varphi)=(a\otimes b,\Delta_\mathcal{B}(\varphi))\quad\text{ and }\quad(1,\varphi)=\varepsilon_\mathcal{B}(\varphi)&,\tag{2.13.b}\label{eq:2.14.b}\\
    (S_\mathcal{A}(a),\varphi)=(a,S_\mathcal{B}(\varphi))&\tag{2.13.c}\label{eq:2.14.c}
\end{align}\
for all $a,b\in \mathcal{A}$ and $\varphi,\psi\in \mathcal{B}$. Here, the pairings in (\ref{eq:2.14.a}) and (\ref{eq:2.14.b}) from $\mathcal{A}\otimes\mathcal{A}\times\mathcal{B}\otimes\mathcal{B}$ to $\mathbb{C}$ is induced from the bilinear form $(\cdot,\cdot)$ and we keep the same notation for simplicity.

There is a unique pairing of Hopf algebras between $\mathcal{U}_q(gl_n)$ and $\mathcal{O}(GL_q(n))$ satisfying
\begin{align}\stepcounter{equation}
    &(q^\lambda,x_{i,j})=\delta_{i,j}q^{\langle\lambda,\epsilon_i\rangle},\tag{2.14.a}\label{eq:2.15.a}\\
    &(e_k,x_{i,j})=\delta_{i,k}\delta_{j,k+1},\quad (f_k,x_{i,j})=\delta_{i,k+1}\delta_{j,k},\tag{2.14.b}\label{eq:2.15.b}\\
    &(q^\lambda,\text{det}_q^m)=q^{m\langle\lambda,\epsilon_1+\cdots+\epsilon_n\rangle}\quad\quad (m\in\mathbb{Z}),\tag{2.14.c}\label{eq:2.15.c}\\
    &(e_k,\text{det}_q^m)=(f_k,\text{det}_q^m)=0\quad\quad (m\in\mathbb{Z})\tag{2.14.d}\label{eq:2.15.d}
\end{align}
where $\text{det}_q^m=D_q^m$ when $m\ge 0$. 

With this unique pairing of $\mathcal{U}_q(gl_n)$ and $\mathcal{O}(GL_q(n))$, we may regard $\mathcal{O}(GL_q(n))$ as a bimodule over $\mathcal{U}_q(gl_n)$. The left and right action of $\mathcal{U}_q(gl_n)$ on $\mathcal{O}(GL_q(n))$ can be defined in the following way. If $a\in \mathcal{U}_q(gl_n)$ and $x\in \mathcal{O}(GL_q(n))$ with $\Delta(x)=\sum x_i^1\otimes x_i^2$, then
\begin{align}\stepcounter{equation}
    a\cdot x=(id\otimes (a,\cdot))\circ\Delta(x)=\sum (a,x_i^2)x_i^1,\tag{2.15.a}\label{eq:2.16.a}\\
    x\cdot a=((a,\cdot)\otimes id)\circ\Delta(x)=\sum (a,x_i^1)x_i^2.\tag{2.15.b}\label{eq:2.16.b}
\end{align}
The actions of generators of $\mathcal{U}_q(gl_n)$ on generators of $\mathcal{O}(GL_q(n))$ are given by
\begin{align}\stepcounter{equation}
    &q^\lambda\cdot x_{i,j}=q^{\langle\lambda,\epsilon_j\rangle}x_{i,j} & &x_{i,j}\cdot q^\lambda=q^{\langle\lambda,\epsilon_i\rangle}x_{i,j}\tag{2.16.a}\label{eq:2.17.a}\\ 
    &e_k\cdot x_{i,j}=\delta_{j,k+1}x_{i,j-1} & &x_{i,j}\cdot e_k=\delta_{i,k}x_{i+1,j}, \tag{2.16.b}\label{eq:2.17.b}\\  
    &f_k\cdot x_{i,j}=\delta_{j,k}x_{i,j+1} & &x_{i,j}\cdot f_k=\delta_{i,k+1}x_{i-1,j}.\tag{2.16.c}\label{eq:2.17.c}
\end{align}
If $a\in\mathcal{U}_q(gl_n)$ such that $\Delta(a)=\sum a_i^1\otimes a_i^2$, we have
\begin{align}
    a\cdot(\psi\varphi)=\sum_i(a_i^1\cdot\psi)(a_i^2\cdot\varphi)\quad\text{ and }\quad (\psi\varphi)\cdot a=\sum_i(\psi\cdot a_i^1)(\varphi\cdot a_i^2)\label{eq:2.18}
\end{align}
due to (\ref{eq:2.14.a}). With (\ref{eq:2.18}), we find that
\begin{align}\stepcounter{equation}
    &e_k\cdot\xi_J^I=\left\{ \begin{array}{cl}
        \xi_{J_{k+1\rightarrow k}}^I & \text{if } k\notin J \text{ and } k+1\in J \\
        0 & \mbox{otherwise} 
        \end{array}\right.\tag{2.18.a}\label{eq:2.19.a}\\
        &f_k\cdot\xi_J^I=\left\{ \begin{array}{cl}
        \xi_{J_{k\rightarrow k+1}}^I & \text{if } k+1\notin J \text{ and } k\in J \\
        0 & \mbox{otherwise} 
        \end{array}\right.\tag{2.18.b}\label{eq:2.19.b}
\end{align}
where $J_{k+1\rightarrow k}$ ($J_{k\rightarrow k+1}$) is the set obtained by removing $k+1$ ($k$) from $J$ and adding $k$ ($k+1$) to $J$.

%% file: Chapters/Chapter_2/finite_dim_rep.tex
\subsection{Irreducible co-representations and the Peter-Weyl decomposition}

A quantum subgroup of a quantum group $G$ is a quantum group $K$ endowed with a surjective homomorphism $\pi_G:G\rightarrow K$ of Hopf algebras.

Define the diagonal subgroup $\mathcal{O}(D_n)$ of $\mathcal{O}(GL_q(n))$ as the commutative ring 
\begin{align}
    \mathcal{O}(D_n)=\mathbb{C}\left[t_1,t_1^{-1}\cdots,t_n,t_n^{-1}\right]\label{eq:2.20}
\end{align}
of Laurent polynomials in $n$ variables $(t_1,\cdots,t_n)$ with co-multiplication $\Delta(t_i)=t_i\otimes t_i$ and co-unit $\varepsilon(t_i)=1$. The surjective homomorphism $\pi_{D_n}:\mathcal{O}(GL_q(n))\rightarrow \mathcal{O}(D_n)$ is given by $\pi_{D_n}(x_{i,j})=\delta_{i,j}t_i$ and $\pi_{D_n}(\text{det}_q^{-1})=\prod_{i=1}^nt_i^{-1}$. 

Define the Borel subgroups $\mathcal{O}(B_n^+)$ ($\mathcal{O}(B_n^-)$) of $\mathcal{O}(GL_q(n))$ by setting $x_{i,j}=0$ for $i>j$ (resp. $i<j$) and $x_{i,j}=0$ for $i<j$ (resp. $i>j$). This is the analogue of the "upper (lower) triangular matrices". Denote $b_{i,j}$ as the modulo class of $x_{i,j}$. Then $\mathcal{O}(B_n^+)$ and $\mathcal{O}(B_n^-)$ is in the form of
\begin{align}\stepcounter{equation}
    \mathcal{O}(B_n^+)&=\mathbb{C}\left[b_{i,j}(i\le j),b_{1,1}^{-1},\cdots ,b_{n,n}^{-1} \right]\tag{2.20.a}\label{eq:2.21.a}\\
    \mathcal{O}(B_n^-)&=\mathbb{C}\left[b_{i,j}(i\ge j),b_{1,1}^{-1},\cdots ,b_{n,n}^{-1} \right]\tag{2.20.b}\label{eq:2.21.b}
\end{align}
with commutation relations induced from (2.1). Notice that $b_{ii}$'s ($1\le i\le n$) on the diagonal commute with each other. The surjective homomorphism $\pi_{B_n^\pm}:\mathcal{O}(GL_q(n))\rightarrow \mathcal{O}(B_n^\pm)$ is given by $\pi_{B_n^\pm}(x_{i,j})=b_{i,j}$ and $\pi_{B_n^\pm}(\text{det}_q^{-1})=\prod_{i=1}^nb_{i,i}^{-1}$.

Given a quantum group $G$, an element $\chi\in G$ is called a linear character of $G$ if
\begin{equation*}
    \Delta(\chi)=\chi\otimes\chi\text{ and }\epsilon(\chi)=1.
\end{equation*}
Given a left $G$ co-module $X$ with the $\mathbf{C}$-algebra homomorphism $L_G:X\rightarrow G\otimes X$, we say $\psi\in X$ is a left relative $G$-invariant with character $\chi$ if $L_G(\psi)=\chi\otimes\psi$. The subspace of all left relative $G$-invariants in $X$ is denoted as
\begin{equation*}
    \mathcal{O}(G\backslash X;\chi)=\{\psi\in X;L_G(\psi)=\chi\otimes\psi\}.
\end{equation*}
If $X$ is a right $G$ co-module with co-action $R_G:X\rightarrow X\otimes G$, the subspace of all right relative $G$-invariants in $X$ is denoted as
\begin{equation*}
    \mathcal{O}(X/G;\chi)=\{\psi\in X;R_G(\psi)=\psi\otimes\chi\}.
\end{equation*}

Let $\lambda=\lambda_1\epsilon_1+\cdots+\lambda_n\epsilon_n$ be an integral weight with all $\lambda_i\ge 0$. We define characters $t^\lambda\in \mathcal{O}(D_n)$ and $b^\lambda\in \mathcal{O}(B_n^\pm)$ by
\begin{align}\stepcounter{equation}
    t^\lambda&=t_{1}^{\lambda_1}\cdots t_{n}^{\lambda_n},\tag{2.21.a}\label{eq:2.22.a}\\
    b^\lambda&=b_{11}^{\lambda_1}\cdots b_{nn}^{\lambda_n}.\tag{2.21.b}\label{eq:2.22.b}
\end{align}

Let $L_{D_n}=(\pi_{D_n}\otimes id)\circ\Delta$ and $R_{D_n}=(id\otimes\pi_{D_n})\circ\Delta$. $\mathcal{O}(GL_q(n))$ becomes a left (resp. right) $\mathcal{O}(D_n)$ co-module with the co-action defined by $L_{D_n}$ (resp. $R_{D_n}$). To describe the left and right relative $\mathcal{O}(D_n)$-invariants in $\mathcal{O}(GL_q(n))$, we define the row sum vector $\alpha(A)$ and column sum vector $\beta(A)$ for a $n\times n$ matrix $A=\{a_{i,j}\}_{i,j=1}^n$ with integer entries as
\begin{align}
    \alpha(A)=\sum_{i=1}^n\left(\sum_{j=1}^na_{i,j}\epsilon_i\right)\in F_n\text{ and }\beta(A)=\sum_{j=1}^n\left(\sum_{i=1}^na_{i,j}\epsilon_j\right)\in F_n.\label{eq:2.23}
\end{align}
If we use $a_{i,j}$ to denote the number of appearance of $x_{i,j}$ in a monomial $x\in \mathcal{O}(GL_q(n))$ and write $A=\{a_{i,j}\}_{i,j=1}^n$, then the left and right character of $x$ is
\begin{align}\stepcounter{equation}
    L_{D_n}(x)&=t^{\alpha(A)}\otimes x,\tag{2.23.a}\label{eq:2.24.a}\\
    R_{D_n}(x)&=x\otimes t^{\beta(A)}.\tag{2.23.b}\label{eq:2.24.b}
\end{align}
If $I=\{i_1,\cdots,i_s\}\subset \{1,\cdots,n\}$, we denote $\vv{I}=\epsilon_{i_1}+\cdots+\epsilon_{i_s}$. The left and right character of quantum minor determinant $\xi_J^I$ in (2.3) is
\begin{align}\stepcounter{equation}
    L_{D_n}(\xi_J^I)&=t^{\vv{I}}\otimes \xi_J^I,\tag{2.24.a}\label{eq:2.25.a}\\
    R_{D_n}(\xi_J^I)&=\xi_J^I\otimes t^{\vv{J}}.\tag{2.24.b}\label{eq:2.25.b}
\end{align}

Let $L_{B_m^-}=(\pi_{B_m^-}\otimes id)\circ\Delta$ and $R_{B_m^+}=(id\otimes\pi_{B_m^+})\circ\Delta$. $\mathcal{O}(GL_q(n))$ becomes a left $\mathcal{O}(B_n^-)$ (resp. right $\mathcal{O}( B_n^+)$) co-module with the co-action defined by $L_{B_m^-}$ (resp. $R_{B_m^+}$). To describe the relative $\mathcal{O}(B_n^-)$-invariants in $\mathcal{O}(GL_q(n))$, we write $\xi_J=\xi_{j_1\cdots j_r}$ (resp. $\xi^J=\xi^{j_1\cdots j_r}$) as the abbreviation of $\xi_{j_1,\cdots,j_r}^{1,\cdots,r}$ (resp. $\xi^{j_1,\cdots,j_r}_{1,\cdots,r}$) with $J=\{j_1<\cdots<j_r\}$ and $\Lambda_r=\epsilon_1+\cdots+\epsilon_r$. Let $J_l$ ($1\le l\le s$) be subsets of $\{1,\cdots,n\}$ and denote $m_l=|J_l|$. If $\psi=\xi_{J_1}\cdots\xi_{J_s}$ and $\varphi=\xi^{J_1}\cdots\xi^{J_s}$, we have
\begin{align}\stepcounter{equation}
    L_{B_m^-}(\psi)&=b^\lambda\otimes\psi\tag{2.25.a}\label{eq:2.26.a}\\
    R_{B_m^+}(\varphi)&=\varphi\otimes b^\lambda\tag{2.25.b}\label{eq:2.26.b}
\end{align}
where $\lambda=\Lambda_{m_1}+\cdots+\Lambda_{m_s}$. Moreover, we have
\begin{align}\label{eq:2.27}
     L_{B_m^-}(\text{det}_q^{-1})=b^{-\Lambda_n}\otimes \text{det}_q^{-1}\text{ and }R_{B_m^+}(\text{det}_q^{-1})=\text{det}_q^{-1}\otimes b^{-\Lambda_n}.
\end{align}
Let $l\in\mathbb{Z}$. Based on (\ref{eq:2.27}) the relative invariant subspaces $\mathcal{O}(B_n^-\backslash GL_q(n);b^\lambda)$ and $\mathcal{O}( GL_q(n)/B_n^+;b^\lambda)$ satisfies
\begin{align}\stepcounter{equation}
    \mathcal{O}(B_n^-\backslash GL_q(n);b^\lambda)=\text{det}_q^{-l}\mathcal{O}\left(B_n^-\backslash GL_q(n);b^{\lambda+l(\epsilon_1+\cdots+\epsilon_n)}\right)\tag{2.27.a}\label{eq:2.28.a}\\
    \mathcal{O}( GL_q(n)/B_n^+;b^\lambda)=\text{det}_q^{-l}\mathcal{O}\left( GL_q(n)/B_n^+;b^{\lambda+l(\epsilon_1+\cdots+\epsilon_n)}\right)\tag{2.27.b}\label{eq:2.28.b}
\end{align}
where $\text{det}_q^{-l}=(D_q)^{-l}$ when $l\le 0$.

Noumi, Yamada, and Mimachi~\cite{noumi1993finite} showed that $\mathcal{O}(B_n^-\backslash GL_q(n);b^\lambda)$ and\\ $\mathcal{O}( GL_q(n)/B_n^+;b^\lambda)$ are irreducible co-representations of $ \mathcal{O}(GL_q(n))$.Moreover, they gave a method to pick $\mathbb{C}$-bases for these co-representations. Here we cite their results. Let $\lambda=\lambda_1\epsilon_1+\cdots+\lambda_n\epsilon_n\in F_n$ with $\lambda_1\ge\lambda_2\ge\cdots\ge\lambda_n$. We will call such a $\lambda$ as a \textit{dominant integral weight}. By a Young diagram of shape $\lambda$ we mean a Young diagram with $\lambda_i$ cubes in the $i$-th row. Denote the set of all semi-standard tableaux of shape $\lambda$ with labels in $\{1,\cdots,n\}$ as $\text{SSYT}_n(\lambda)$. For $\mathbf{T}=(T_{r,s})\in \text{SSYT}_n(\lambda)$, let $J_s$ ($1\le s\le \lambda_1$) be the content of the $s$-th column of $\mathbf{T}$. Denote $\xi_{\mathbf{T}}$ and $\xi^{\mathbf{T}}$ as the product of quantum minor determinants
\begin{align}\stepcounter{equation}
    \xi_{\mathbf{T}}&=\xi_{J_1}\cdots\xi_{J_{\lambda_1}},\tag{2.28.a}\label{eq:2.29.a}\\
    \xi^{\mathbf{T}}&=\xi^{J_1}\cdots\xi^{J_{\lambda_1}}.\tag{2.28.b}\label{eq:2.29.b}
\end{align}
With the notation in (\ref{eq:2.28.a}) and (\ref{eq:2.28.b}), we have:
\begin{theorem}[NYM]\label{thm:2.1}
    Let $\lambda=(\lambda,\cdots,\lambda_n)$ be a dominant integral weight and $l\in\mathbb{Z}$ with $\lambda_n\ge -l$. Monomials in the form of $\text{det}_q^{-l}\xi_{\mathbf{T}}$ (resp. $\text{det}_q^{-l}\xi^{\mathbf{T}}$) indexed by $\mathbf{T}\in \text{SSYT}_n(\lambda+l(\epsilon_1+\cdots+\epsilon_n))$ form a $\mathbb{C}$-basis for the right (resp. left) $\mathcal{O}(GL_q(n))$ co-module $\mathcal{O}(B_n^-\backslash GL_q(n);b^\lambda)$ (resp. $\mathcal{O}( GL_q(n)/B_n^+;b^\lambda)$).
\end{theorem}


In the following, we give the definitions needed for the Peter-Weyl decomposition of $\mathcal{O}(GL_q(n))$. We will always assume that $\lambda$ is a dominant integral weight and denote
\begin{equation}\label{eq:irrep}
        V^R(\lambda)=\mathcal{O}\left(B_n^-\backslash GL_q(n);b^\lambda\right)\text{ and }V^L(\lambda)=\mathcal{O}\left( GL_q(n)/B_n^+;b^\lambda\right).
\end{equation}
The dual of $V^R(\lambda)$ is denoted as $V^{*}_R(\lambda)=\text{Hom}_\mathbb{C}(V^R(\lambda),\mathbb{C})$. $V^{*}_R(\lambda)$ has a unique left $\mathcal{O}(GL_q(n))$ co-module structure $L^*_{B_m^-}:V^{*}_R(\lambda)\rightarrow \mathcal{O}(GL_q(n))\otimes V^{*}_R(\lambda)$ induced by $L_{B_m^-}$
\begin{align}
    \left(L^*_{B_m^-}(u),v\right)=\left(u,L_{B_m^-}(v)\right)\text{ for all }u\in V^{*}_R(\lambda)\text{ and }v\in V^R(\lambda),\label{eq:2.31}
\end{align}
where $(\cdot,\cdot)$ is the contraction $V^{*}_R(\lambda)\otimes V^R(\lambda)\rightarrow \mathbb{C}$. Let $\{v_i\}_{i\in I}$ be a $\mathbb{C}$-basis of $V^R(\lambda)$ and $\omega_{i,j}$ ($i,j\in I$) be the corresponding matrix coefficients such that
\begin{align}
    L_{B_m^-}(v_j)=\sum_{i\in I}v_i\otimes\omega_{i,j}.\label{eq:2.32}
\end{align}
If $\{u_i\}_{i\in I}$ is the basis of $V^{*}_R(\lambda)$ dual to $\{v_i\}_{i\in I}$, i.e., $(u_i,v_j)=\delta_{i,j}$, then we have
\begin{align}
    L^*_{B_m^-}(u_i)=\sum_{j\in I} \omega_{i,j}\otimes u_j.\label{eq:2.33}
\end{align}
Define homomorphism $\Phi_\lambda:V^{*}_R(\lambda)\otimes V^R(\lambda)\rightarrow \mathcal{O}(GL_q(n))$ as
\begin{align}
    \Phi_\lambda(u_i\otimes v_j)=\omega_{i,j}.\label{eq:2.34}
\end{align}
Denote the image of $\Phi_\lambda$ as $W(\lambda)$. We remark that $W(\lambda)$ is a two-sided irreducible $\mathcal{O}(GL_q(n))$ co-module since we have
\begin{align}
    \Delta(\omega_{i,j})=\sum_{k\in I}\omega_{i,k}\otimes\omega_{k,j}\text{ and }\varepsilon(\omega_{i,j})=\delta_{i,j}.\label{eq:2.35}
\end{align}
Noumi, Yamada, and Mimachi showed that
\begin{theorem}[NYM]\label{thm:2.2} The following statements are true.
    \begin{enumerate}
        \item [1)] $V^{*}_R(\lambda)$ is isomorphic to $V^L(\lambda)$.
        \item[2)] The homomorphism $\Phi_\lambda$ is an isomorphism and  $W(\lambda)$ is isomorphic to $V^L(\lambda)\otimes V^R(\lambda)$.
    \end{enumerate}
\end{theorem}
\noindent Then, they give the Peter-Weyl decomposition of $\mathcal{O}(GL_q(n))$.
\begin{theorem}[NYM]\label{thm:2.3}
    $\mathcal{O}(GL_q(n))$ is decomposed in the the direct sum of irreducible two-sided $\mathcal{O}(GL_q(n))$ co-modules
    \begin{align}
        \mathcal{O}(GL_q(n))=\bigoplus_\lambda W(\lambda),\label{eq:2.36}
    \end{align}
    where the summation runs over all dominant integral weight $\lambda$.
\end{theorem}
Readers can find the explicit expressions of $\omega_{i,j}$'s on $\mathcal{O}(SL_q(2))$ in the book of Klimyk and Schm{\"u}dgen~\cite{klimyk2012quantum}. On $\mathcal{O}(SU_q(3))$, the explicit expressions of a part of $\omega_{i,j}$'s are calculated in the work of Bergeron, Koelink, and Vinet~\cite{bergeron2019suq}. Their method can be used to calculate all $\omega_{i,j}$'s on $\mathcal{O}(SU_q(3))$. Currently, the explicit expressions of matrix coefficients on higher rank $\mathcal{O}(U_q(n))$ are unknown.

%% file: Chapters/Chapter_2/haar_state.tex
\subsection{The Haar state and the method of  Noumi, Yamada, and Mimachi}
The Haar state $h:\mathcal{O}(GL_q(n))\rightarrow \mathbb{C}$ is the unique linear functional satisfying the translation invariant property
\begin{align}
    (id\otimes h)\circ\Delta(x)=h(x)\cdot 1=(h\otimes id)\circ\Delta(x)\text{ for all }x\in \mathcal{O}(GL_q(n))
    \label{eq:2.37}
\end{align}
where $1$ is the multiplicative unit of $\mathcal{O}(GL_q(n))$. Noumi, Yamada, and Mimachi showed that $h$ is actually the projection of $\mathcal{O}(GL_q(n))$ onto the trivial $\mathcal{O}(GL_q(n))$ co-module $W((0,\cdots,0))=\mathbb{C}$ whose matrix coefficient is $1\in \mathcal{O}(GL_q(n))$. Hence, if we take the limit $q\rightarrow 1$ of the Haar state $h$ on $\mathcal{O}(GL_q(n))$, we get the Haar state on $\mathcal{O}(U_1(n))$ which is the integration with respect to the normalized Haar measure on the unitary group $U(n)$. Moreover, shown by Reshetikhin and Yakimov~\cite{reshetikhin2001quantum}, we can write the Haar state on the quantized algebra of coordinate functions on the compact real form of a simply connected Lie group as an integral. In their expression, the Haar state is an integral over the maximal torus of the compact real form and the integrand is the trace of an operator on an infinite dimensional space.

Recall the action of $a\in\mathcal{U}_q(gl_n)$ on $\varphi\in\mathcal{O}(GL_q(n))$ (\ref{eq:2.16.a}) and (\ref{eq:2.16.b}). By the invariance property of $h$, we have
\begin{align}\stepcounter{equation}
    h(a\cdot\varphi)&=(a,(h\otimes id)\circ\Delta(\varphi))=(a,h(\varphi)\cdot 1)=h(\varphi)\varepsilon_{\mathcal{U}_q}(a)\tag{2.37.a}\label{eq:2.38}\\
    h(\varphi\cdot a)&=(a,(id\otimes h)\circ\Delta(\varphi))=(a,h(\varphi)\cdot 1)=h(\varphi)\varepsilon_{\mathcal{U}_q}(a)\tag{2.37.b}\label{eq:2.38+}
\end{align}

On the compact real form $\mathcal{O}(U_q(n))$ ($q\in\mathbb{R}$), Noumi, Yamada, and Mimachi defined two positive definite Hermitian forms
\begin{align}
    \langle\varphi,\psi\rangle_L=h(\varphi^*\psi)\text{ and }\langle\varphi,\psi\rangle_R=h(\varphi\psi^*)\label{eq:2.39}
\end{align}
and showed matrix coefficients of different dominant integral weights are orthogonal with respect to these two Hermitian forms. They also gave an explicit description of the modular property of the between the two Hermitian forms. Define the \textit{modular automorphism} $\rho:\mathcal{O}(U_q(n))\rightarrow \mathcal{O}(U_q(n))$ as
\begin{align}
    \rho(x_{i,j})=q^{2n+2-2i-2j}\cdot x_{i,j}\quad\text{and}\quad \rho(det_q^{-1})=det_q^{-1}.\label{eq:2.40}
\end{align}
\begin{theorem}[NYM]\label{thm:2.4}
    The modular automorphism $\rho$ satisfies
    \begin{align}
        h(\varphi\psi)=h(\rho(\psi)\varphi)\text{  for all  }\varphi,\psi\in \mathcal{O}(U_q(n)).\label{eq:2.41}
    \end{align}
    Hence, we have
    \begin{align}
        \langle\varphi,\psi\rangle_L=\langle\rho(\psi),\varphi\rangle_R\text{  for all  }\varphi,\psi\in \mathcal{O}(U_q(n)).\label{eq:2.42}
    \end{align}
\end{theorem}
Finally, we discuss the method to compute the norm of matrix coefficients of the Peter-Weyl decomposition proposed by Noumi, Yamada, and Mimachi. Let $\lambda\in F_n$ be a dominant integral weight and fix $\mathbf{T}\in\text{SSYT}(\lambda)$. Denote the columns of $\mathbf{T}$ by $J_s$ ($1\le s\le \lambda_1$). Then the right character of the basis vector $\xi_{\mathbf{T}}$ is
\begin{align}
    R_{D_n}(\xi_{\mathbf{T}})=\xi_{\mathbf{T}}\otimes t^{\vv{J_1}+\cdots+\vv{J_{\lambda_1}}}.\label{eq:2.43}
\end{align}
Hence, $V^R(\lambda)$ has the following decomposition
\begin{align}
    V^R(\lambda)=\bigoplus_{\mu\in F_n}V^R_\mu(\lambda);\text{ where } V^R_\mu(\lambda)=\left\{v\in V^R(\lambda); R_{D_n}(v)=v\otimes t^\mu\right\}.\label{eq:2.44}
\end{align}
This decomposition is orthogonal with respect to $\langle\cdot,\cdot\rangle_L$ and $\langle\cdot,\cdot\rangle_R$. In other words, if $v_\eta\in V^R_\eta(\lambda)$ and $v_\gamma\in V^R_\gamma(\lambda)$ with $\eta\ne\gamma$, then $\langle v_\eta,v_\gamma\rangle_L=0$ and $\langle v_\eta,v_\gamma\rangle_R=0$.
Let $\{u_i\}_{i\in I(\lambda)}$ be a $\mathbb{C}$-basis for $V^R(\lambda)$ with $R_{D_n}(u_i)=u_i\otimes t^{\lambda(i)}$. Denote $\{\omega_{i,j}\}_{i,j\in I(\lambda)}$ as the matrix coefficients with respect to $\{u_i\}_{i\in I(\lambda)}$. Noumi, Yamada, and Mimachi gave the formula of the square lenght of $\omega_{i,j}$ with the assumption that $\langle u_i,u_j\rangle_L=\delta_{i,j}$.
\begin{theorem}[NYM]
    Assume that $\{u_i\}_{i\in I(\lambda)}$ is an orthonormal basis of $V^R(\lambda)$ with respect to $\langle\cdot,\cdot\rangle_L$ and $R_{D_n}(u_i)=u_i\otimes t^{\lambda(i)}$. The basis $\{\omega_{i,j}\}_{i,j\in I(\lambda)}$ for the two-sided $\mathcal{O}(U_q(n))$ co-modules $W(\lambda)$ is orthogonal with respect to $\langle\cdot,\cdot\rangle_L$ and $\langle\cdot,\cdot\rangle_R$. The square length of $\omega_{i,j}$ is given by the formula
    \begin{align}\stepcounter{equation}
        \langle\omega_{i,j},\omega_{i,j}\rangle_L&=\frac{q^{2\langle\rho,\lambda(i)\rangle}}{d_\lambda},\tag{2.44.a}\label{eq:2.45.a}\\
        \langle\omega_{i,j},\omega_{i,j}\rangle_R&=\frac{q^{-2\langle\rho,\lambda(j)\rangle}}{d_\lambda},\tag{2.44.b}\label{eq:2.45.b}
    \end{align}
    where $2\rho=\sum_{k=1}^n(n+1-2k)\epsilon_k$. $d_\lambda$ is the $q$-analogue of the dimension of $V^R(\lambda)$ given by
    \begin{align}
        d_\lambda=\sum_{i\in I(\lambda)}q^{2(\rho,\lambda(i))}.\label{eq:2.46}
    \end{align}
\end{theorem}
We remark that a similar result holds for orthonormal basis of $V^L(\lambda)$. For more detail, see Theorem 3.7 in Noumi, Yamada, and Mimachi's paper~\cite{noumi1993finite}.

%% file: Chapters/Chapter_3.tex
\section{The Haar state values of monomials on $\mathcal{O}(GL_q(n))$}
Recall that on $U(n)$ we can write the conjugate of the coordinate function $U_{i,j}(g)$ as
\begin{equation}
    \begin{split}
        &\overline{U_{i,j}(g)}=U^{-1}_{i,j}(g)=\frac{\text{det}_{\hat{j},\hat{i}}(g)}{\text{det}(g)}
    \end{split}
\end{equation}
where $\text{det}_{\hat{j},\hat{i}}(g)$ denotes the determinant of the $(j,i)$ cofactor of matrix $g\in U(n)$. Since $\text{det}_{\hat{j},\hat{i}}(g)$ is polynomial consisting of $U_{i,j}(g)$'s, the evaluation of (\ref{eq:1.1}) can be broke down to the evaluation of 
\begin{equation}\label{eq:1.4}
    \int_{U\in G}\frac{U_{i_1,j_1}\cdots U_{i_m,j_m}}{\text{det}^{a_m}(U)}\text{d}U.
\end{equation}
The analogue of (\ref{eq:1.4}) on $\mathcal{O}(GL_q(n))$ is 
\begin{equation}\tag{1.5}
    h(x_{i_1,j_1}\cdots x_{i_m,j_m}\cdot \text{det}_q^{-a_m})
\end{equation}
and we will discuss the evaluation problem of (\ref{eq:1.5}) in this section. We start our discussion with the characterization of those monomials with non-zero Haar state values on $\mathcal{O}(GL_q(n))$. Then we select a special set of monomials called psudo-basis of order $m$ such that every monomial with a non-zero Haar state value can be decomposed as a linear combination of these psudo-basis vectors. We will give a method to construct linear relations between the Haar state values of these psudo-basis vectors. As an application, we solve the Haar state values of monomials in our psudo-basis of order $1$. Although for large rank $n$ and order $m$ a method to enumerate psudo-basis vectors is unknown, we propose the construction of Source matrix of order $m$ on $\mathcal{O}(U_q(n))$ which is used to solve $h\left((\prod_{i=1}^nx_{i,n+1-i})^m\cdot \text{det}_q^{-m}\right)$. Finally, we present two (anti)automorphisms that preserves the Haar state value.

\input{Chapters/Chapter_3/Character_of_non-zero}

\input{Chapters/Chapter_3/Linear_relation_Haar_states}

\input{Chapters/Chapter_3/Linear_system_of_order_1}

\input{Chapters/Chapter_3/Construction_of_the_Source_matrix}

\input{Chapters/Chapter_3/Value_preserving_maps}

%% file: Chapters/Chapter_3/Character_of_non-zero.tex
\subsection{Monomials with non-zero Haar state values}

Let $\lambda,\mu\in F_n$ be integral weight. Define the left and right relative $\mathcal{O}(D_n)$-invariant subspace $\mathcal{A}[\lambda,\mu]$ of $\mathcal{O}(GL_q(n))$ as
\begin{align}\label{eq:3.1}
    \mathcal{A}[\lambda,\mu]=\{x\in \mathcal{O}(GL_q(n)); L_{D_n}(x)=t^\lambda\otimes x, R_{D_n}(x)=x\otimes t^\mu\}.
\end{align}
Notice that $L_{D_n}(\text{det}_q^{-1})=t^{-\vv{1}}\otimes \text{det}_q^{-1}, R_{D_n}(\text{det}_q^{-1})=\text{det}_q^{-1}\otimes t^{-\vv{1}}$ where $-\vv{1}=-\epsilon_1-\cdots-\epsilon_n$. We can decompose $\mathcal{O}(GL_q(n))$ as 
\begin{align}\label{eq:3.2}
    \mathcal{O}(GL_q(n))=\bigoplus_{\lambda,\mu\in F_n}\mathcal{A}[\lambda,\mu].
\end{align}
and $\mathcal{O}(GL_q(n))$ has a grading structure induced by $\mathcal{A}[\lambda_1,\mu_1]\cdot\mathcal{A}[\lambda_2,\mu_2]\subset\mathcal{A}[\lambda_1+\lambda_2,\mu_1+\mu_2]$. 

Define $L_h:\mathcal{O}(GL_q(n))\rightarrow\mathcal{O}(D_n)$ and $R_h:\mathcal{O}(GL_q(n))\rightarrow\mathcal{O}(D_n)$ as
\begin{align}\addtocounter{equation}{1}
    L_h&=(\pi_{D_n}\otimes h)\circ\Delta,\tag{3.5.a}\label{eq:3.3.a}\\
    R_h&=(h\otimes \pi_{D_n})\circ\Delta.\tag{3.5.b }\label{eq:3.3.b}
\end{align}
For $x\in \mathcal{A}[\lambda,\mu]$, we have
\begin{align}\addtocounter{equation}{1}
    L_h(x)&=(id\otimes h)\circ L_{D_n}(x)=h(x)\cdot t^{\lambda},\tag{3.6.a}\label{eq:3.4.a}\\
    R_h(x)&=(h \otimes id)\circ R_{D_n}(x)=h(x)\cdot t^{\mu}.\tag{3.6.b}\label{eq:3.4.b}
\end{align}
On the other hand, by (\ref{eq:2.37}) we have
\begin{align}\addtocounter{equation}{1}
    L_h(x)&=(L_{D_n}\otimes id)(id\otimes h)\circ\Delta(x)=h(x)\cdot L_{D_n}(1)=h(x)\cdot 1,\tag{3.7.a}\label{eq:3.5.a}\\
    R_h(x)&=(id\otimes R_{D_n})(h\otimes id)\circ\Delta(x)=h(x)\cdot R_{D_n}(1)=h(x)\cdot 1.\tag{3.7.b}\label{eq:3.5.b}
\end{align}
Comparing (3.6) and (3.7), we conclude that
\begin{proposition}\label{prop:3.3}
    If $x\in\mathcal{O}(GL_q(n))$ and $h(x)\ne 0$, then $x\in \mathcal{A}\left[\vv{0},\vv{0}\right]$ where $\vv{0}$ is the zero vector in $F_n$.
\end{proposition}

Next, we give the characterization of those $x\in\mathcal{A}\left[\vv{0},\vv{0}\right]$ after making several definitions. By a \textit{monomial} on $\mathcal{O}(GL_q(n))$, we mean an element in the form of $x\cdot\text{det}_q^{-m}$ ($m\ge 0$) where $x$ is a finite product of generators $x_{i,j}\in \mathcal{O}(GL_q(n))$.
\begin{definition}\label{def:3.1}
    Given a monomial $x\cdot\text{det}_q^{-m}\in \mathcal{O}(GL_q(n))$, the counting matrix of $x\cdot\text{det}_q^{-m}$, denoted as $\theta(x)$, is an $n\times n$ matrix whose $(i,j)$-th entry equals to the number of appearance of generator $x_{i,j}$ in $x$. The counting matrix of monomials in the form of $\text{det}_q^{-m}$ is the zero matrix. The row sum vector and column sum vector of $\theta(x)$ is denoted by $\alpha(x)$ and $\beta(x)$, respectively.
\end{definition}
\begin{definition}\label{def:3.4}
    Let $A$ be an $n\times n$ matrix with integer entries and $k\in\mathbb{N}$. We say that $A$ is an $n\times n$ $m$-doubly stochastic matrix if $\alpha(A)=\beta(A)=m\cdot\vv{1}$.
\end{definition}
\begin{theorem}\label{thm:3.4}
    Let $x\cdot\text{det}_q^{-m}\in \mathcal{O}(GL_q(n))$ be a monomial such that $h\left(x\cdot\text{det}_q^{-m}\right)\ne 0$. Then, the counting matrix of $x\cdot\text{det}_q^{-m}$, $\theta(x)$, is an $m$-doubly stochastic matrix.
\end{theorem}
\begin{proof}
    By (\ref{eq:3.4.a}) and (\ref{eq:3.5.a}), we find that 
    $$h\left(x\cdot\text{det}_q^{-m}\right)\cdot t^{\alpha(x)-m\cdot\vv{1}}=h\left(x\cdot\text{det}_q^{-m}\right)\cdot 1.$$
    Since $h\left(x\cdot\text{det}_q^{-m}\right)\ne 0$, we must have $\alpha(x)-m\cdot\vv{1}=\vv{0}$, i.e., $\alpha(x)=m\cdot\vv{1}$. Similarly, by (\ref{eq:3.4.b}) and (\ref{eq:3.5.b}), we have $\beta(x)=m\cdot\vv{1}$. Thus, $\theta(x)$ is an $m$-doubly stochastic matrix.
\end{proof}
\noindent Denote 
$$\mathcal{A}\left[\vv{0},\vv{0}\right]_m=\{x\cdot\text{det}_q^{-m}\in\mathcal{O}(GL_q(n))| \theta(x)\text{ is an }m\text{-doubly stochastic matrix}\}.$$
We say that monomials in $\mathcal{A}\left[\vv{0},\vv{0}\right]_m$ are \textit{monomials of order} $m$. It follows that
$$\mathcal{A}\left[\vv{0},\vv{0}\right]=\bigoplus_{m\in\mathbb{N}}\mathcal{A}\left[\vv{0},\vv{0}\right]_m.$$
Since $h$ is a linear functional, we want to find a $\mathbb{C}$-basis for $\mathcal{A}\left[\vv{0},\vv{0}\right]_m$ consisting of monomials. To serve this purpose, we define a linear order on the set of $m\times n$ matrices induced by the lexicographical order on vectors. For each $m\times n$ matrix $C=\{c_{i,j}\}_{i=1}^m|_{j=1}^n$ we associate the vector $\mathcal{V}(C)\in\mathbb{R}^{mn}$ by
\begin{align*}
    \mathcal{V}(C)=(c_{1,1},c_{1,2},\cdot,c_{1,n},c_{2,1}\cdots,c_{m,n})
\end{align*}
We will say $C\prec D$ if $\mathcal{V}(C)\prec\mathcal{V}(D)$ in the lexicographical order of vectors. We make the following observation on the process of reordering the generators of $\mathcal{O}(GL_q(n))$ in a monomial $x\cdot\text{det}_q^{-m}$. 

\hfill

\noindent{\it Let $x=P\cdot x_{ik}x_{jl}\cdot Q\cdot\text{det}_q^{-m}$ be a monomial of order $m$ where $i<j,k<l$ and $P,Q$ are two monomials containing no generator $\text{det}_q^{-1}$. If we switch the order of $x_{ik}$ and $x_{jl}$ so that}
\begin{equation*}
    x\cdot\text{det}_q^{-m}=y\cdot\text{det}_q^{-m}+(q-q^{-1})z\cdot\text{det}_q^{-m},
\end{equation*}
{\it where $y=P\cdot x_{jl}x_{ik}\cdot Q$ and $z=P\cdot x_{il}x_{jk}\cdot Q$, then $z\cdot\text{det}_q^{-m}$ is a monomial of order $m$ and $\theta(z)\prec\theta(x)=\theta(y)$.}

\hfill

\noindent In other words, $\theta(z)$ has the same row sum and column sum vector as $\theta(x)$ and $\theta(y)$ but the order of $\theta(z)$ drops. This observation provide a criterion to pick monomial basis for each $\mathcal{A}\left[\vv{0},\vv{0}\right]_m$. Let $B_n(m)$ be the set of all $n\times n$ $m$-doubly stochastic matrices.
\begin{proposition}\label{prop:3.5}
    Suppose that for each $M\in B_n(m)$ we fix one and only one monomial $x_M\cdot\text{det}_q^{-m}\in \mathcal{A}\left[\vv{0},\vv{0}\right]_m$. Then, every monomial $\varphi\in \mathcal{A}\left[\vv{0},\vv{0}\right]_m$ can be written as a linear combination of these $x_M\cdot\text{det}_q^{-m}$'s.
\end{proposition}
\begin{proof}
    Based on our observation, we have
    \begin{equation}\label{eq:3.6}
        \varphi=c\cdot x_{\theta(\varphi)}\cdot\text{det}_q^{-m}+\sum_{i\in I}c_i\cdot y_i\cdot\text{det}_q^{-m}
    \end{equation}
    where $I$ is a finite set. The coefficients are Laurent polynomials in variable $q$. $y_i$'s are obtained by applying commutation relation (\ref{eq:2.1.c}) with $B_n(m)\ni\theta(y_i)\prec\theta(\varphi)$. We may repeat the process (\ref{eq:3.6}) to each $y_i$ if necessary. Again, any new monomials generated by (\ref{eq:2.1.c}) during the process is a monomial of order $k$ with the order of its counting matrix less than $\theta(y_i)$. The process will terminate in finite steps since $B_n(m)$ is a finite set.
\end{proof}
\begin{definition}\label{def:3.6}
    The set $P_m=\left\{x_M\cdot\text{det}_q^{-m}\in \mathcal{A}\left[\vv{0},\vv{0}\right]_m; M\in B_n(m)\right\}$ chosen in Proposition \ref{prop:3.5} is called a set of psudo-basis of order $m$
\end{definition}
 We remark that each $P_m$ is only linear independent within $\mathcal{A}\left[\vv{0},\vv{0}\right]_m$. The union $P$ of all chosen $P_m$
\begin{equation}
    P=\bigcup_{m\in\mathbb{N}} P_m\label{eq:3.7}
\end{equation}
is not a linear independent set in $\mathcal{A}\left[\vv{0},\vv{0}\right]$ due to the linearity induced by $D_q\cdot\text{det}_q^{-1}=1$.


%% file: Chapters/Chapter_3/Linear_relation_Haar_states.tex
\subsection{Linear relation between the Haar state values of psudo-basis vectors}
Although the set of all psudo-basis vectors does not form a monomial basis of $\mathcal{A}\left[\vv{0},\vv{0}\right]$, they arise naturally in the computation of the translation operation of $h$. To demonstrate the relation between the translation operation of a monomial of order $m$ and psudo-basis vectors of order $m$, we make the following remark on the co-multiplication of a monomial. If $x=\prod_{k\in K}x_{i_k,j_k}\cdot\text{det}_q^{-m}$ is a monomial, then
\begin{align}\label{eq:3.10}
    \Delta(x\cdot\text{det}_q^{-m})&=\prod_{k\in K}\left(\sum_{l_k=1}^nx_{i_k,l_k}\cdot\text{det}_q^{-m}\otimes x_{l_k,j_k}\cdot\text{det}_q^{-m}\right)\\
    &=\sum_{i\in I'} z_i\cdot\text{det}_q^{-m}\otimes y_i\cdot\text{det}_q^{-m}.\notag
\end{align}
where the summation rums over all tensor products  $z_i\cdot\text{det}_q^{-m}\otimes y_i\cdot\text{det}_q^{-m}$ satisfying the \textbf{Order Restriction} of co-multiplication:
\begin{itemize}
    \item[i)] the $k$-th generator of the left component $z_i$ is in the $i_k$-th row
    \item[ii)] the $k$-th generator of the right component $y_i$ is in the $j_k$-th column
    \item[iii)] The column index of the $l$-th generator in $z_i$ equals to the row index of the $l$-th generator in $y_i$.
\end{itemize}
As a direct consequence of the order restriction of co-multiplication, we have the following lemma.
\begin{lemma}\label{lemma:3.6}
    Let $z_i\cdot\text{det}_q^{-m}\otimes y_i\cdot\text{det}_q^{-m}$ be a tensor product appearing in (\ref{eq:3.10}). Then, we have
    \begin{equation}\label{eq:3.11}
        \alpha(x)=\alpha(z_i), \beta(z_i)=\alpha(y_i), \beta(y_i)=\beta(x).
    \end{equation}
\end{lemma}

\begin{proposition}\label{prop:3.7}
    Let $x\cdot\text{det}_q^{-m}$ be a monomial of order $m$. We have
    \begin{equation}\label{eq:3.12}
        (id\otimes h)\Delta(x\cdot\text{det}_q^{-m})=\sum_{i\in I}h(y_i\cdot\text{det}_q^{-m})z_i\cdot\text{det}_q^{-m}
    \end{equation}
    where $I$ is a finite set consisting of those $i$'s such that $h(y_i\cdot\text{det}_q^{-m})\ne 0$. For every $i\in I$, $y_i\cdot\text{det}_q^{-m}$ and $z_i\cdot\text{det}_q^{-m}$ are monomials of order $m$. 
\end{proposition}
\begin{remark}
    Similar result is true for the left translation $(h\otimes id)\Delta(x\cdot\text{det}_q^{-m})$.
\end{remark}
\begin{proof}
    Assume that $\Delta(x)\cdot\text{det}_q^{-m}$ is in the same form as in (\ref{eq:3.10}). Then, we have
    \begin{equation}\label{eq:3.13}
        (id\otimes h)\Delta(x\cdot\text{det}_q^{-m})=\sum_{i\in I'}h(y_i\cdot\text{det}_q^{-m})z_i\cdot\text{det}_q^{-m}.
    \end{equation}
    By Lemma \ref{lemma:3.6} we must have $\beta(y_i)=\beta(x)=m\cdot \Vec{1}$ and $\alpha(z_i)=\alpha(x)=m\cdot \Vec{1}$. If $h(y_i)\ne 0$, then by Theorem \ref{thm:3.4} $\alpha(y_i)=\beta(y_i)=m\cdot \Vec{1}$. By Lemma \ref{lemma:3.6} again, $\beta(z_i)=\alpha(y_i)=m\cdot \Vec{1}$. Hence, $h(y_i\cdot\text{det}_q^{-m})\ne 0$ implies that both $y_i\cdot\text{det}_q^{-m}$ and $z_i\cdot\text{det}_q^{-m}$ are monomials of order $m$. Let $I$ be the subset of $I'$ consisting of those indices $i$ such that $h(y_i\cdot\text{det}_q^{-m})\ne 0$, then we get (\ref{eq:3.12}) in the described form.
\end{proof}
Fix a set of psudo-basis of order $m$, $P_m$, and an $x_M\cdot\text{det}_q^{-m}\in P_m$. Combining Proposition \ref{prop:3.5} with Proposition \ref{prop:3.7}, we get 
\begin{align}\label{eq:3.14}
    (id\otimes h)\Delta(x_M\cdot\text{det}_q^{-m})=\sum_{L}\left(\sum_{K}c^M_{K,L}h(x_K\cdot\text{det}_q^{-m})\right)\cdot x_L\cdot\text{det}_q^{-m}
\end{align}
where $x_K\cdot\text{det}_q^{-m}$ and $x_L\cdot\text{det}_q^{-m}$ run over all psudo-basis of order $m$ in $P_m$. On the other hand, notice that $D_q^m$ is a summation of monomials of order $m$. Hence, we have
\begin{align}\label{eq:3.15}
    1=D_q^m\cdot \text{det}_q^{-m}=\sum_L b_L\cdot x_L\cdot \text{det}_q^{-m}
\end{align}
where $x_L\cdot\text{det}_q^{-m}$ run over all psudo-basis of order $m$ in $P_m$. Substituting (\ref{eq:3.14}) and (\ref{eq:3.15}) into the translation invariance property $(id\otimes h)\Delta(x_M\cdot\text{det}_q^{-m})=h(x_M\cdot\text{det}_q^{-m})\cdot 1$, we get
\begin{equation}\label{eq:3.16}
    \begin{split}   &\sum_{L}\left(\sum_{K}c^M_{K,L}h(x_K\cdot\text{det}_q^{-m})\right)\cdot x_L\cdot\text{det}_q^{-m}\\
    &=\sum_L b_Lh(x_M\cdot\text{det}_q^{-m})\cdot x_L\cdot \text{det}_q^{-m}.
    \end{split}
\end{equation}
By comparing the coefficients of each $x_L\cdot \text{det}_q^{-m}$, we get
\begin{equation}\label{eq:3.17}
\begin{split}
    &\sum_{K}c^M_{K,L}h(x_K\cdot\text{det}_q^{-m})=b_Lh(x_M\cdot\text{det}_q^{-m})\\
    \Longleftrightarrow&\sum_{K}c^M_{K,L}h(x_K\cdot\text{det}_q^{-m})-b_Lh(x_M\cdot\text{det}_q^{-m})=0.
\end{split}
\end{equation}
We will say (\ref{eq:3.17}) is a linear relation derived from \textit{equation basis} $x_M\cdot\text{det}_q^{-m}$ and \textit{comparing basis} $x_L\cdot \text{det}_q^{-m}$. Notice that the equation basis and the comparing basis in (\ref{eq:3.17}) could come from two different sets of psudo-basis of order $m$. The total number of linear relations in the form of (\ref{eq:3.17}) is $|B_n(m)|^2$.

Our aim is to construct a linear system consisting of the Haar states of psudo-basis of order $m$. This will be a linear system with $|B_n(m)|$ unknowns. All linear relations in the form of (\ref{eq:3.17}) have zero right-hand-side. Hence, we need a linear relation with non-zero right-hand-side. Evaluating the Haar state on both sides of (\ref{eq:3.15}), we get
\begin{align}\label{eq:3.18}
    1=h(1)=\sum_L b_Lh(x_L\cdot \text{det}_q^{-m}).
\end{align}
Thus, we can pick linear relation (\ref{eq:3.18}) together with $|B_n(m)|-1$ many linear relations in the form of (\ref{eq:3.17}) to build an invertible linear system. We will call such an invertible linear system as a \textit{linear system of order $m$}. Here we remark that to form an invertible linear system, the $|B_n(m)|-1$ many linear relations in the form of (\ref{eq:3.17}) cannot be chosen randomly. 

%% file: Chapters/Chapter_3/Linear_system_of_order_1.tex
\subsection{Linear system of order $1$ on $\mathcal{O}(GL_q(n))$}
In this subsection we will pick a set of psudo-basis of order $1$ and construct a linear system of order $1$ to find the Haar state values of our chosen psudo-basis vectors.

$B_n(1)$ is just the set of all $n\times n$ permutation matrices and there is an one-to-one correspondence between $B_n(1)$ and the permutation group on $n$ letters $\mathfrak{S}_n$. Hence, we will use permutations $\sigma\in\mathfrak{S}_n$ to denote those matrices in $B_n(1)$. The set of psudo-basis of order $1$ we pick are in the form of
\begin{align}\label{eq:3.19}
    x_\sigma=\prod_{i=1}^nx_{i,\sigma(i)}\cdot\text{det}_q^{-1}.
\end{align}
With this choice we have
\begin{align}\stepcounter{equation}
    1=D_q\cdot\text{det}_q^{-1}=\sum_{\tau\in\mathfrak{S}_n}(-q)^{l(\tau)}\cdot x_{\tau}.\tag{3.20.a}\label{eq:3.20.a}\\
    1=h(D_q\cdot\text{det}_q^{-1})=\sum_{\tau\in\mathfrak{S}_n}(-q)^{l(\tau)}h(x_{\tau})\tag{3.20.b}\label{eq:3.20.b}
\end{align}
Taking the right translation operation of our psudo-basis, we get
\begin{equation}\label{eq:3.21}
\begin{split}
    (id\otimes h)\Delta(x_\sigma)&=\sum_{\tau\in\mathfrak{S}_n}h\left(\prod_{i=1}^nx_{\tau(i),\sigma(i)}\cdot\text{det}_q^{-1}\right)\cdot\prod_{i=1}^nx_{i,\tau(i)}\cdot\text{det}_q^{-1}\\
    &=\sum_{\tau\in\mathfrak{S}_n}h\left(\prod_{i=1}^nx_{\tau(i),\sigma(i)}\cdot\text{det}_q^{-1}\right)\cdot x_{\tau}
\end{split}
\end{equation}
due to Proposition \ref{prop:3.7} and the order restriction of the comultiplication. 

Denote $\sigma_0\in\mathfrak{S}_n$ as the permutation such that $\sigma_0(i)=n+1-i$. for those $\sigma\ne\sigma_0$, we denote $\sigma_1\in\mathfrak{S}_n$ as the permutation such that $\sigma_1(i)=n+1-\sigma(i)$. To pick the $|B_n(1)|-1$ many linear relations in the form of (\ref{eq:3.17}), consider the linear relation  derived from equation basis $x_\sigma$ and comparing basis $x_{\sigma_1}$ for all $\mathfrak{S}_n\ni\sigma\ne\sigma_0$. We have
\begin{align}\label{eq:3.22}
    h\left(\prod_{i=1}^nx_{\sigma_1(i),\sigma(i)}\cdot\text{det}_q^{-1}\right)=(-q)^{l(\sigma_1)}h(x_\sigma).
\end{align}
Notice that $x_{\sigma_1(i),\sigma(i)}=x_{n+1-\sigma(i),\sigma(i)}$ ($1\le i\le n$) is located on the reverse diagonal of the generator matrix (\ref{mat:generator_n}). Hence, they commute with each other and we have
\begin{align}\label{eq:3.23}
    \prod_{i=1}^nx_{\sigma_1(i),\sigma(i)}\cdot\text{det}_q^{-1}=\prod_{i=1}^nx_{i,n+1-i}\cdot\text{det}_q^{-1}=x_{\sigma_0}.
\end{align}
Since $l(\sigma_1)=\frac{n(n+1)}{2}-l(\sigma)$ and $l(\sigma_0)=\frac{n(n+1)}{2}$, (\ref{eq:3.22}) becomes
\begin{align}\label{eq:3.24}
    h(x_{\sigma_0})=(-q)^{l(\sigma_0)-l(\sigma)}h(x_{\sigma}).
\end{align}
Substituting (\ref{eq:3.24}) into (\ref{eq:3.20.b}), we get
\begin{align}\stepcounter{equation}
    &\left(\sum_{\sigma\in\mathfrak{S}_n}(-q)^{2l(\sigma)-l(\sigma_0))}\right)h(x_{\sigma_0})=1\tag{3.25.a}\label{eq:3.25.a}\\
    \Longleftrightarrow&h(x_{\sigma_0})=\frac{(-q)^{l(\sigma_0)}}{\sum_{\sigma\in\mathfrak{S}_n}q^{2l(\sigma)}}.\tag{3.25.b}\label{eq:3.25.b}
\end{align}
Substituting (\ref{eq:3.25.b}) into (\ref{eq:3.24}), we get
\begin{align}\label{eq:3.26}
    h(x_{\sigma})=\frac{(-q)^{l(\sigma)}}{\sum_{\sigma\in\mathfrak{S}_n}q^{2l(\sigma)}}
\end{align}
for all $\sigma\in\mathfrak{S}_n$. let $I_n(k)$ be the number of permutations on $n$ letters with $k$ inversions. The denominator in (\ref{eq:3.26}) can be written as
\begin{equation}
    \sum_{\sigma\in\mathfrak{S}_n}q^{2l(\sigma)}=\sum_{k=0}^{\frac{n(n+1)}{2}}I_n(k)q^{2k}.
\end{equation}
By Andrews~\cite{andrews1998theory}, the generating function of $I_n(k)$ is
\begin{equation}
    \sum_{k=0}^{\frac{n(n+1)}{2}}I_n(k)x^{k}=\prod_{j=1}^n\frac{1-x^j}{1-x}. 
\end{equation}
Denote $[n]_q=\frac{1-q^n}{1-q}$. The denominator of (\ref{eq:3.26}) can be rewritten as
\begin{equation}
    \sum_{\sigma\in\mathfrak{S}_n}q^{2l(\sigma)}=\prod_{j=1}^n\frac{1-q^{2j}}{1-q^2}=[n]_{q^2}!,
\end{equation}
and we get
\begin{equation}\label{eq:3.30}
    h(x_{\sigma})=\frac{(-q)^{l(\sigma)}}{[n]_{q^2}!}
\end{equation}
for all $\sigma\in\mathfrak{S}_n$.

%% file: Chapters/Chapter_3/Construction_of_the_Source_matrix.tex
\subsection{Construction of the Source matrix on $\mathcal{O}(U_q(n))$}

 For large $n$ and $m$ the enumeration of matrices in $B_n(m)$ becomes a primary issue of constructing higher order linear systems. The number of matrices in $B_n(m)$ is known by Stein and Stein~\cite{stein1970enumeration} and Odama and Musiker~\cite{odama2001enumeration}. However, their results are obtained by the method of generating function and cannot be used to enumerate matrices in $B_n(m)$ directly. As a consequence, there is no explicit method to pick a set of psudo-basis of order $m$ on $\mathcal{O}(GL_q(n))$ when $n$ and $m$ are large. 

However, solving $h(x^m_{\sigma_0})$ on $\mathcal{O}(U_q(n))$ is still necessary since $x^m_{\sigma_0}$ is the unique monomial of order $m$ (up to reordering) corresponding to the counting matrix
\begin{align}\label{eq:3.31}
    R_m=\begin{bmatrix}
        0&0&\cdots&0&m\\
        0&0&\cdots&m&0\\
    \vdots&\vdots&\iddots&\vdots&\vdots\\  
    0&m&\cdots&0&0\\
    m&0&\cdots&0&0
    \end{bmatrix}\in B_n(m).
\end{align}
In this sense, the Haar state of psudo-basis vector $h(x_{R_m})$($=h(x^m_{\sigma_0})$) is invariant under our choice of different $P_m$ in Proposition \ref{prop:3.5}. 

In the following we describe the construction of the \textit{Source matrix of order} $m$ by solving which we get $h(x^m_{\sigma_0})$ on $\mathcal{O}(U_q(n))$. We restrict ourselves on $\mathcal{O}(U_q(n))$ since the modular automorphism (\ref{eq:2.40}) is a key ingredient of our construction.

We start with introducing the unknowns of the Source matrix of order $m$. Let $\sigma\in\mathfrak{S}_n$ and $M_\sigma\in B_n(1)$ the permutation matrix corresponding to $\sigma$. Denote $M_m^\sigma=M_\sigma+R_{m-1}$ where $R_{m-1}\in B_n(m-1)$ is defined in (\ref{eq:3.31}). The psudo-basis of order $m$ corresponds to $M_m^\sigma$, denoted as $x_{m}^\sigma$, is defined as 
\begin{align}\label{eq:3.32}
    x_{m}^\sigma=x_\sigma\cdot x_{\sigma_0}^{m-1}=\left(\prod_{i=1}^nx_{i,\sigma(i)}\right)\left(\prod_{i=1}^nx^{m-1}_{i,n+1-i}\right)\cdot \text{det}_q^{-m}.
\end{align}
The Haar states of these $x_{m}^\sigma$'s will be the unknowns of the Source matrix of order $m$. Notice that $h(x_{m}^{\sigma_0})=h(x_{\sigma_0}^{m})$ is what we are seeking for.


The linear relations between $h(x_{m}^\sigma)$'s are derived using (\ref{eq:3.17}) and the property of the modular automorphism (\ref{eq:2.40}). We will fix the equation basis as $x_{\sigma_0}^{m}$. The comparing bases in our construction are defined as follow. Let $L_m$ be the psudo-basis of order $m$ consisting of those monomials of order $m$ whose generators are arranged in the (row) lexicographical order, i.e., generators are ordered first by row index then by column index. The comparing basis we pick are those psudo-basis vectors in $L_m$, denoted as $z_m^\sigma$, corresponding to $M_\sigma+(m-1)I_n\in B_n(m)$ where $I_n$ is the $n\times n$ identity matrix. Each $z_m^\sigma$ can be written as
\begin{align}\stepcounter{equation}
    z_{m}^\sigma=\prod_{r=1}^nz_{r(\sigma)}=z_{1(\sigma)}z_{2(\sigma)}\cdots z_{n(\sigma)}\cdot \text{det}_q^{-m}\tag{3.33.a}\label{eq:3.33.a}
\end{align}
where
\begin{align}\tag{3.33.b}\label{eq:3.33.b}
    z_{r(\sigma)}=\left\{ \begin{array}{rcl}
x_{r,r}^{m-1}\cdot x_{r,\sigma(r)} & \mbox{if}
& r\le \sigma(r) \\ x_{r,\sigma(r)}\cdot x_{r,r}^{m-1} & \mbox{if} & r> \sigma(r)
\end{array}\right.
\end{align}
is the partial monomial representing the $r$-th row of $M_\sigma+(m-1)I_n\in B_n(m)$. 

Fix the order of generators in $x_{\sigma_0}^m$ as  $x_{\sigma_0}^m=\prod_{i=1}^nx^{m}_{i,n+1-i}\cdot \text{det}_q^{-m}$. Applying the right translation operation to $x_{\sigma_0}^m$, we get
\begin{equation}\label{eq:3.33}
\begin{split}
    (id\otimes h)\circ\Delta(x_{\sigma_0}^m)&=\sum_{i\in I}\zeta^i_1\zeta^i_2\cdots \zeta^i_n\cdot\text{det}_q^{-m}\otimes h\left(\eta^i_1\eta^i_2\cdots \eta^i_n\cdot\text{det}_q^{-m}\right)\\
    &=\sum_{i\in I}h\left(\eta^i_1\eta^i_2\cdots \eta^i_n\cdot\text{det}_q^{-m}\right)\zeta^i_1\zeta^i_2\cdots \zeta^i_n\cdot\text{det}_q^{-m}
\end{split}
\end{equation}
where each pair of $(\zeta^i_r, \eta^i_r)$ comes from the comultiplication of $x^m_{r,n+1-r}$ and $I$ is the index set such that each $\zeta^i_1\zeta^i_2\cdots \zeta^i_n\cdot\text{det}_q^{-m}\otimes \eta^i_1\eta^i_2\cdots \eta^i_n\cdot\text{det}_q^{-m}$ is a tensor product of two monomials of order $m$. For a fixed $i$, each $\zeta^i_r$ $(1\le r\le n)$ is in the form of
\begin{align}\label{eq:3.34}
    \zeta^i_r=x_{r,k(i,r,1)}x_{r,k(i,r,2)}\cdots x_{r,k(i,r,m)}
\end{align}
and each $\eta^i_r$ $(1\le r\le n)$ is in the form of
\begin{align}
    \eta^i_r=x_{k(i,r,1),n+1-r}x_{k(i,r,2),n+1-r}\cdots x_{k(i,r,m),n+1-r}
\end{align}
where $k(i,r,j)$ represents the column index of the $j$-th generator in $\zeta^i_r$. Notice that generators in each $\zeta^i_1\zeta^i_2\cdots \zeta^i_n\cdot\text{det}_q^{-m}$ are already ordered by their row index. Hence, to rewrite $\zeta^i_1\zeta^i_2\cdots \zeta^i_n\cdot\text{det}_q^{-m}$ as a linear combination of psudo-basis vectors in $L_m$, it suffices to change the order of generators in each $\zeta^i_r$ using the row switch rule in (\ref{eq:2.1.a}). We remark that this reorder process of $\zeta^i_r$ does not generate extra monomials. Therefore, to compute the linear relation (\ref{eq:3.17}) derived from the equation basis $x_{\sigma_0}^m$ and a comparing basis $x_N\in L_m$, it suffices to consider all $\zeta^i_1\zeta^i_2\cdots \zeta^i_n\cdot\text{det}_q^{-m}$ in (\ref{eq:3.33}) whose counting matrix equals to the counting matrix of $x_N$. 

\textbf{Claim:} The linear relations derived from equation basis $x_{\sigma_0}^m$ and comparing basis $z_m^\sigma$ can be written as linear relations between $h(x_m^\tau)$'s in (\ref{eq:3.32}).

\textbf{Check:} Without loss of generality, we fix $\sigma\in\mathfrak{S}_n$. Recall that the counting matrix of $z_m^\sigma$ is $M_\sigma+(m-1)I_n$. To compute such a linear relation, the tensor products in (\ref{eq:3.33}) under our consideration are those $\zeta^i_1\zeta^i_2\cdots \zeta^i_n\cdot\text{det}_q^{-m}\otimes \eta^i_1\eta^i_2\cdots \eta^i_n\cdot\text{det}_q^{-m}$ whose counting matrix of $\zeta^i_1\zeta^i_2\cdots \zeta^i_n\cdot\text{det}_q^{-m}$ equals to $M_\sigma+(m-1)I_n$. Equation (\ref{eq:3.34}) shows that each $\zeta^i_r$ represents the $r$-th row of the counting matrix of $\zeta^i_1\zeta^i_2\cdots \zeta^i_n\cdot\text{det}_q^{-m}$. Hence, each $\zeta^i_r$ is in the form of 
\begin{align}\label{eq:3.37}
    \zeta^i_r=x_{r,r}^{f_r(i)}\cdot x_{r,\sigma(r)}\cdot x_{r,r}^{m-1-f_r(i)}
\end{align}
where $f_r(i):I\mapsto \{0,1,\cdots,m-1\}$. The corresponding $\eta^i_r$ is in the form of
\begin{align}\label{eq:3.38}
    \eta^i_r=x_{r,n+1-r}^{f_r(i)}\cdot x_{\sigma(r),n+1-r}\cdot x_{r,n+1-r}^{m-1-f_r(i)}.
\end{align}
Hence, to verify our claim, it suffices to show that $h(\eta^i_1\eta^i_2\cdots \eta^i_n\cdot\text{det}_q^{-m})$ with each $\eta^i_r$ in (\ref{eq:3.38}) can be written as a linear combination of $h(x_m^\sigma)$'s in (\ref{eq:3.32}). 
\begin{definition}
    Define the \textit{polarity of a generator} $x_{i,j}$, denoted as $p(x_{i,j})$, by
    \begin{align}
        p(x_{i,j})=\left\{ \begin{array}{ccl}
        \text{positive} & \mbox{if}
        & n+1-i-j>0\\ \text{neutral} & \mbox{if}
        & n+1-i-j=0\\ \text{negative} & \mbox{if}
        & n+1-i-j<0
        \end{array}\right..
    \end{align}
\end{definition}
$p(x_{i,j})$ is positive/neutral/negative if and only if $x_{i,j}$ is above/on/below the reverse diagonal of the generator matrix (\ref{mat:generator_n}), respectively. In $\eta^i_1\eta^i_2\cdots \eta^i_n\cdot\text{det}_q^{-m}$, a positive (negative) $x_{\sigma(r), n+1-r}$ commutes with all non-positive (non-negative) generators to its right (left). Therefore, we are able to reorder $\eta^i_1\eta^i_2\cdots \eta^i_n\cdot\text{det}_q^{-m}$ such that all positive (negative) $x_{\sigma(r), n+1-r}$'s are at the very right (left) end with their relative order unchanged just using switching rules (\ref{eq:2.1.a}) and (\ref{eq:2.1.b}). We remark that no extra monomial is generated during this process.


\textit{Example 3.1:} Let $\sigma=(13542)$. On $\mathcal{O}(U_q(5))$, one of the tensor products in (\ref{eq:3.33}) is 
\begin{equation*}
\begin{split}
    &\underbrace{(x_{1,1})^{m-1}x_{1,3}}_{\zeta_1^\sigma}\underbrace{x_{2,1}(x_{2,2})^{m-1}}_{\zeta_2^\sigma}\underbrace{(x_{3,3})^{m-1}x_{3,5}}_{\zeta_3^\sigma}\underbrace{x_{4,2}(x_{4,4})^{m-1}}_{\zeta_4^\sigma}\underbrace{x_{5,4}(x_{5,5})^{m-1}}_{\zeta_5^\sigma}\cdot\text{det}_q^{-m}\otimes\\
        &\underbrace{(x_{1,5})^{m-1}x_{3,5}}_{\eta_1^\sigma}\underbrace{x_{1,4}(x_{2,4})^{m-1}}_{\eta_2^\sigma}\underbrace{(x_{3,3})^{m-1}x_{5,3}}_{\eta_3^\sigma}\underbrace{x_{2,2}(x_{4,2})^{m-1}}_{\eta_4^\sigma}\underbrace{x_{4,1}(x_{5,1})^{m-1}}_{\eta_5^\sigma}\cdot\text{det}_q^{-m}.
\end{split}
\end{equation*}
In (\ref{eq:3.37}), the $f_r(\sigma)$'s correspond to $\zeta_1^\sigma\zeta_2^\sigma\zeta_3^\sigma\zeta_4^\sigma\zeta_5^\sigma\cdot\text{det}_q^{-m}$ are $f_1(\sigma)=m-1; f_2(\sigma)=0; f_3(\sigma)=m-1; f_4(\sigma)=0; f_5(\sigma)=0$. Notice that generators in $\zeta_1^\sigma\zeta_2^\sigma\zeta_3^\sigma\zeta_4^\sigma\zeta_5^\sigma\cdot\text{det}_q^{-m}$ are in the row lexicographical order. The generators of $\eta_1^\sigma\eta_2^\sigma\eta_3^\sigma\eta_4^\sigma\eta_5^\sigma\cdot\text{det}_q^{-m}$ are listed in the following table.
\begin{center}
\begin{tabular}{|c|c|c|c|c|}
    \hline
        &  &  &\cellcolor{cyan!25} $x_{1,4}$ &\cellcolor{green!20} $x_{1,5}$ \\
    \hline     
        &\cellcolor{cyan!25} $x_{2,2}$ &  &\cellcolor{green!20} $x_{2,4}$ & \\
    \hline
        &  &\cellcolor{green!20} $x_{3,3}$ &  &\cellcolor{red!15}  $x_{3,5}$\\
    \hline
        \cellcolor{cyan!25}$x_{4,1}$&\cellcolor{green!20} $x_{4,2}$& & & \\
    \hline
        \cellcolor{green!20}$x_{5,1}$& &\cellcolor{red!15}$x_{5,3}$& &\\
    \hline
        \multicolumn{1}{c}{$\eta_5^\sigma$} & 
        \multicolumn{1}{c}{$\eta_4^\sigma$} & 
        \multicolumn{1}{c}{$\eta_3^\sigma$} &
        \multicolumn{1}{c}{$\eta_2^\sigma$} &
        \multicolumn{1}{c}{$\eta_1^\sigma$}
\end{tabular}
\end{center}
In the table, positive/neutral/negative generators are colored in blue/green/red, respectively. Notice that the $\eta_i$'s in the table are listed in decreasing order. The left side of the table corresponds to the right side of $\eta_1\eta_2\eta_3\eta_4\eta_5\cdot\text{det}_q^{-m}$. By observing the table, we have get
\begin{samepage}
\begin{equation}\label{eq:3.40}
    \begin{split}
        &\eta_1\eta_2\eta_3\eta_4\eta_5\cdot\text{det}_q^{-m}\\
        =&q^{5m-5}\cdot x_{3,5}x_{5,3}(x_{1,5}x_{2,4}x_{3,3}x_{4,2}x_{5,1})^{m-1}\cdot\text{det}_q^{-m}x_{1,4}x_{2,2}x_{4,1}.
    \end{split}
\end{equation}\qed
\end{samepage}

Observe from Example 3.1 that the relative order between positive (and negative) generators in $\eta^i_1\eta^i_2\cdots \eta^i_n\cdot\text{det}_q^{-m}$ is the reversed column lexicographical order. Assume that we have reordered generators in $\eta^i_1\eta^i_2\cdots \eta^i_n\cdot\text{det}_q^{-m}$ as
\begin{align}\label{eq:3.41}
    \eta^i_1\eta^i_2\cdots \eta^i_n\cdot\text{det}_q^{-m}=q^c\cdot n_\sigma\left(\prod_{i=1}^nx_{i,n+1-i}^{m-1}\right)\cdot\text{det}_q^{-m} e_\sigma p_\sigma
\end{align}
where $p_\sigma$ ($n_\sigma$) is the product of all positive (negative) generators in the monomial in the reversed column lexicographical order and $e_\sigma$ is the product of extra neutral generators not presented in $\left(\prod_{i=1}^nx_{i,n+1-i}^{m-1}\right)$. Evaluating the Haar state on both side of (\ref{eq:3.41}) and applying the modular automorphism $\rho$ in (\ref{eq:2.40}), we get
\begin{equation}\label{eq:3.42}
    \begin{split}
        h\left(\eta^i_1\eta^i_2\cdots \eta^i_n\cdot\text{det}_q^{-m}\right)&=q^c h\left(\rho(e_\sigma p_\sigma)n_\sigma\left(\prod_{i=1}^nx_{i,n+1-i}^{m-1}\right)\cdot\text{det}_q^{-m}\right)\\
        &=q^{c+d} h\left(e_\sigma p_\sigma n_\sigma\left(\prod_{i=1}^nx_{i,n+1-i}^{m-1}\right)\cdot\text{det}_q^{-m}\right).
    \end{split}
\end{equation}
Notice that $e_\sigma p_\sigma n_\sigma$ is a reordering of $\prod_{r=1}^nx_{\sigma(r),n+1-r}$ and its counting matrix belongs to $B_n(1)$. Hence, by the same argument used in the proof of Proposition (\ref{prop:3.5}), we know that $e_\sigma p_\sigma n_\sigma$ can be written as 
\begin{equation}\label{eq:3.43}
    \begin{split}
        e_\sigma p_\sigma n_\sigma=\sum_{\tau\in\mathfrak{S}_n}c_{\tau}\cdot\prod_{i=1}^nx_{i,\tau(i)},
    \end{split}
\end{equation}
where $c_{\tau}$'s are Laurent polynomials in variable $q$ generated by the commutation rule (2.1.a-c). Substituting (\ref{eq:3.43}) into (\ref{eq:3.42}), we get
\begin{equation}\label{eq:3.44}
    \begin{split}
        &h\left(\eta^i_1\eta^i_2\cdots \eta^i_n\cdot\text{det}_q^{-m}\right)\\
        =&q^{c+d} \sum_{\tau\in\mathfrak{S}_n}c_{\tau} h\left(\left(\prod_{i=1}^nx_{i,\tau(i)}\right)\left(\prod_{i=1}^nx_{i,n+1-i}^{m-1}\right)\cdot\text{det}_q^{-m}\right)\\
        =&q^{c+d} \sum_{\tau\in\mathfrak{S}_n}c_{\tau}h(x^\tau_m).
    \end{split}
\end{equation}

\textit{Example 3.2:} In (\ref{eq:3.40}), $n_\sigma=x_{3,5}x_{5,2}$, $e_\sigma=1$, and $p_\sigma=x_{1,4}x_{2,2}x_{4,1}$. Evaluating the Haar state on both sides of (\ref{eq:3.40}) and applying the modular automorphism, we get
\begin{samepage}
\begin{equation}\label{eq:3.45}
    \begin{split}
        &h\left(\eta_1\eta_2\eta_3\eta_4\eta_5\cdot\text{det}_q^{-m}\right)\\
        =&q^{5m-5}h\left(x_{3,5}x_{5,3}(x_{1,5}x_{2,4}x_{3,3}x_{4,2}x_{5,1})^{m-1}\cdot\text{det}_q^{-m}x_{1,4}x_{2,2}x_{4,1}\right)\\
        =&q^{5m-5}h\left(\rho(x_{1,4}x_{2,2}x_{4,1})x_{3,5}x_{5,3}(x_{1,5}x_{2,4}x_{3,3}x_{4,2}x_{5,1})^{m-1}\cdot\text{det}_q^{-m}\right)\\
        =&q^{5m-5+2\cdot 4}h\left((x_{1,4}x_{2,2}x_{4,1})x_{3,5}x_{5,3}(x_{1,5}x_{2,4}x_{3,3}x_{4,2}x_{5,1})^{m-1}\cdot\text{det}_q^{-m}\right)\\
        =&q^{5m+3}h\left(x_{1,4}x_{2,2}x_{3,5}x_{4,1}x_{5,3}(x_{1,5}x_{2,4}x_{3,3}x_{4,2}x_{5,1})^{m-1}\cdot\text{det}_q^{-m}\right)\\
    \end{split}
\end{equation}\qed
\end{samepage}

Recall that we can write $\zeta^i_1\zeta^i_2\cdots \zeta^i_n\cdot\text{det}_q^{-m}=q^{\lambda(i)}z_m^\sigma$ since generators in $\zeta^i_1\zeta^i_2\cdots \zeta^i_n\cdot\text{det}_q^{-m}$ are ordered by their row index. Hence by (\ref{eq:3.44}), we can write
\begin{equation}\label{eq:3.46}
    \begin{split}
        &h\left(\eta^i_1\eta^i_2\cdots \eta^i_n\cdot\text{det}_q^{-m}\right)\zeta^i_1\zeta^i_2\cdots \zeta^i_n\cdot\text{det}_q^{-m}\\
        =&q^{c+d+\lambda(i)} \left(\sum_{\tau\in\mathfrak{S}_n}c_{\tau}h(x^\tau_m)\right)\cdot z_m^\sigma.
    \end{split}
\end{equation}
Recall the construction of (\ref{eq:3.17}). (\ref{eq:3.46}) implies that the linear relation derived from equation basis $x_{\sigma_0}^m$ and comparing basis $z_m^\sigma$ can be written as a linear combination of $h(x_m^\tau)$'s in (\ref{eq:3.32}).\qed

We remark that the usage of the modular automorphism $\rho$ in (\ref{eq:3.42}) is essential since in general $\eta^i_1\eta^i_2\cdots \eta^i_n\cdot\text{det}_q^{-m}$ in (\ref{eq:3.38}) cannot be written as a linear combination of $x_m^\sigma$ in (\ref{eq:3.32}). Finding the linear decomposition of a monomial of order $m$ into a set of psudo-basis of order $m$ is different from finding the linear decomposition of \textit{the Haar state of} a monomial of order $m$ into \textit{the Haar state of} a set of psudo-basis of order $m$.

\textit{Example 3.3:} Let $m=2$. (\ref{eq:3.45}) gives
\begin{equation}
    \begin{split}
        &h\left((x_{1,5}x_{3,5})(x_{1,4}x_{2,4})(x_{3,3}x_{5,3})(x_{2,2}x_{4,2})(x_{4,1}x_{5,1})\cdot\text{det}_q^{-2}\right)\\
    =&q^{13}h\left((x_{1,4}x_{2,2}x_{3,5}x_{4,1}x_{5,3})(x_{1,5}x_{2,4}x_{3,3}x_{4,2}x_{5,1})\cdot\text{det}_q^{-2}\right).
    \end{split}
\end{equation}
On the other hand, we have
\begin{equation}\label{eq:3.48}
    \begin{split}
        &(x_{1,5}x_{3,5})(x_{1,4}x_{2,4})(x_{3,3}x_{5,3})(x_{2,2}x_{4,2})(x_{4,1}x_{5,1})\cdot\text{det}_q^{-2}\\
        =&q^{-2}(x_{1,5}x_{3,5})(x_{1,4}x_{5,3}x_{2,4}x_{3,3})(x_{2,2}x_{4,1}x_{4,2}x_{5,1})\cdot\text{det}_q^{-2}\\
        =&q^{-3}(x_{3,5}x_{1,4}x_{5,3}x_{2,2}x_{4,1})(x_{1,5}x_{2,4}x_{3,3}x_{4,2}x_{5,1})\cdot\text{det}_q^{-2}\\
        &-(q-q^{-1})q^{-2}(x_{1,5}x_{3,5})(x_{1,4}x_{5,3}x_{2,4}x_{3,2})(x_{2,3}x_{4,1}x_{4,2}x_{5,1})\cdot\text{det}_q^{-2}.
    \end{split}
\end{equation}
In (\ref{eq:3.48}), $(x_{1,5}x_{3,5})(x_{1,4}x_{5,3}x_{2,4}x_{3,2})(x_{2,3}x_{4,1}x_{4,2}x_{5,1})\cdot\text{det}_q^{-2}$ in the last line is obtained by switching $x_{3,3}$ and $x_{2,2}$ in monomial of the second line. Notice that the counting matrix of $(x_{1,5}x_{3,5})(x_{1,4}x_{5,3}x_{2,4}x_{3,2})(x_{2,3}x_{4,1}x_{4,2}x_{5,1})\cdot\text{det}_q^{-2}$ is not in the form of $M_2^\sigma$ for any $\sigma\in\mathfrak{S}_n$. Hence, we conclude that $(x_{1,5}x_{3,5})(x_{1,4}x_{2,4})(x_{3,3}x_{5,3})$\\$(x_{2,2}x_{4,2})(x_{4,1}x_{5,1})\cdot\text{det}_q^{-2}$ is not a linear combination of $x_2^\sigma$ in (\ref{eq:3.32}).\qed

The Source matrix of order $m$ consists of the $n!-1$ linear relations derived from equation basis $x^m_{\sigma_0}$ and comparing basis $z_m^\sigma$ with $\sigma\ne id$ and the linear relation with non-zero right-hand-side
\begin{align}\label{eq:3.49}
    h(x_{\sigma_0}^{m-1})=&h(D_q\cdot \text{det}_q^{-1}x_{\sigma_0}^{m-1})=\sum_{\sigma\in\mathfrak{S}_n}(-q)^{l(\sigma)}h(x_m^\sigma).
\end{align}
Hence, Source matrices of different orders form a recursive linear system. The solution of the Source matrix of order $m-1$ serves as the right-hand-side of the Source matrix of order $m$. The Source matrix of order $1$ is already solved in Subsection 3.3.

The primary computations involved in constructing the Source matrix of order $m$ are
\begin{enumerate}
    \item[1] the coefficients of $z_m^\sigma$ in $D_q^m$ as in (\ref{eq:3.15}), and
    \item[2]  for each $\sigma\in\mathfrak{S}_n$, the decomposition of $e_\sigma p_\sigma n_\sigma$ into a linear combination of monomials in the form of $\prod_{i=1}^nx_{i,\tau(i)}$ as in (\ref{eq:3.43}). 
\end{enumerate}
Notice that for each $\mathcal{O}(U_q(n))$, the coefficients in (\ref{eq:3.43}) for each $\sigma\in\mathfrak{S}_n$ is independent of $m$ and only need to compute once. 

The Source matrix of order $m$ is invertible since each linear relation derived from equation basis $x^m_{\sigma_0}$ and comparing basis $z_m^\sigma$ has different leading term (in the order of counting matrix). More precisely, for a fixed $\sigma\in\mathfrak{S}_n$, let $\varsigma=\sigma_0\circ\sigma^{-1}$ so that $\varsigma$ maps $\sigma(r)$ to $n+1-r$. Then, the leading term in (\ref{eq:3.43}) is $\prod_{i=1}^nx_{i,\varsigma(i)}$. Hence, the leading term of the linear relation derived from equation basis $x^m_{\sigma_0}$ and comparing basis $z_m^\sigma$ is $h(x_m^\varsigma)$. As $\sigma\ne id$ goes over $\mathfrak{S}_n$, the leading term $h(x_m^\varsigma)$ goes over our unknowns in the Source matrix of order $m$ except $h(x_m^{\sigma_0})$. We can find the linear relations between $h(x_m^\varsigma)$'s and $h(x^m_{\sigma_0})$ recursively according to the order of their counting matrices using these linear relations derived from equation basis $x^m_{\sigma_0}$ and different comparing basis $z_m^\sigma$'s. Substituting all the linear relations between $h(x_m^\varsigma)$ and $h(x^m_{\sigma_0})$ into (\ref{eq:3.49}), we get the linear relation between $h(x^m_{\sigma_0})$ and $h(x^{m-1}_{\sigma_0})$ and the Source matrix of order $m$ is solved. This process is an analogue of the method we used in solving the linear system of order $1$ is Subsection 3.3.

%% file: Chapters/Chapter_3/Value_preserving_maps.tex
\subsection{Two value-preserving (anti)automorphisms on $\mathcal{O}(GL_q(n))$}
In this subsection we introduce two (anti)automorphisms that preserve the Haar state values of elements in $\mathcal{O}(GL_q(n))$.

Define the \textit{diagonal flip automorphism} $\gamma:\mathcal{O}(Mat_q(n))\rightarrow \mathcal{O}(Mat_q(n))$ by
\begin{align}\label{eq:3.50}
    \gamma(x_{i,j})=x_{j,i}.
\end{align}
\begin{lemma}\label{lemma:3.8}
    The following statements are true.
    \begin{enumerate}
        \item[1] $\gamma(D_q)=D_q$, and
        \item[2] $\Delta\circ\gamma=\iota\circ(\gamma\otimes\gamma)\circ\Delta$ where $\iota$ is the map that flips the left and right component of a tensor product. 
    \end{enumerate}
\end{lemma}
\begin{proof}
\begin{enumerate}
    \item[1.] Let $\sigma\in\mathfrak{S}_n$. We have
    \begin{equation}
        \begin{split}
            \gamma(D_q)&=\sum_{\sigma\in\mathfrak{S}_n}(-q)^{l(\sigma)}\gamma\left(\prod_{i=1}^nx_{i,\sigma(i)}\right)\\
            &=\sum_{\sigma\in\mathfrak{S}_n}(-q)^{l(\sigma)}\prod_{i=1}^nx_{\sigma(i),i}=D_q.
        \end{split}
    \end{equation}
    \item[2.] It suffices to verify the relation on generators of $\mathcal{O}(Mat_q(n))$. We have
    \begin{equation}
        \begin{split}
            \Delta\circ\gamma(x_{i,j})&=\Delta(x_{j,i})=\sum_{k=1}^nx_{j,k}\otimes x_{k,i}=\tau\left(\sum_{k=1}^nx_{k,i}\otimes x_{j,k}\right)\\
        =&\tau\circ(\gamma\otimes\gamma)\left(\sum_{k=1}^nx_{i,k}\otimes x_{k,j}\right)=\tau\circ(\gamma\otimes\gamma)\circ\Delta(x_{i,j}).
        \end{split}
    \end{equation}
\end{enumerate}
\end{proof}
By Lemma (\ref{lemma:3.8}.1), we can extend $\gamma$ onto $\mathcal{O}(GL_q(n))$ by 
\begin{align}
    \gamma(\text{det}_q^{-1})=\text{det}_q^{-1}.
\end{align}

Define the \textit{double flip anti-automorphism} $\omega:\mathcal{O}(Mat_q(n))\rightarrow \mathcal{O}(Mat_q(n))$ by
\begin{align}\label{eq:3.52}
    \omega(x_{i,j})=x_{n+1-i,n+1-j}.
\end{align}
\begin{lemma}\label{lemma:3.9}
    The following statements are true.
    \begin{enumerate}
        \item[1] $\omega(D_q)=D_q$, and
        \item[2] $\Delta\circ\omega=(\omega\otimes\omega)\circ\Delta$. 
    \end{enumerate}
\end{lemma}
\begin{proof}
    \begin{enumerate}
        \item[1] Let $\sigma\in\mathfrak{S}_n$. We have
        \begin{equation}
        \begin{split}
            \omega\left(\prod_{i=1}^nx_{i,\sigma(i)}\right)&=\prod_{i=n}^1x_{n+1-i,n+1-\sigma(i)}=\prod_{j=1}^nx_{j,n+1-\sigma(n+1-j)}=\prod_{j=1}^nx_{j,\tau(j)}.
        \end{split}
        \end{equation}
        where $\tau=\sigma_0\sigma\sigma_0$. It is easy to check that $l(\tau)=l(\sigma)$. We have
        \begin{equation}
            \begin{split}
                \omega(D_q)&=\sum_{\sigma\in\mathfrak{S}_n}(-q)^{l(\sigma)}\omega\left(\prod_{i=1}^nx_{i,\sigma(i)}\right)\\
                &=\sum_{\sigma\in\mathfrak{S}_n}(-q)^{l(\sigma)}\prod_{i=1}^nx_{i,\sigma_0\sigma\sigma_0(i)}\\
                &=\sum_{\tau\in\mathfrak{S}_n}(-q)^{l(\tau)}\prod_{i=1}^nx_{i,\tau(i)}=D_q\\
            \end{split}
        \end{equation}
        \item[2] It suffices to verify the relation on generators of $\mathcal{O}(Mat_q(n))$. We have
        \begin{equation}
        \begin{split}
            \Delta\circ\omega(x_{i,j})&=\sum_{k=1}^n x_{n+1-i,k}\otimes x_{k,n+1-j}\\
            &=\sum_{k=1}^n x_{n+1-i,n+1-k}\otimes x_{n+1-k,n+1-j}\\
            &=\sum_{k=1}^n \omega(x_{i,k})\otimes \omega(x_{k,j})=(\omega\otimes\omega)\circ\Delta(x_{i,j}).
        \end{split}
    \end{equation}
    \end{enumerate}
\end{proof}
By Lemma (\ref{lemma:3.9}.1), we can extend $\omega$ onto $\mathcal{O}(GL_q(n))$ by 
\begin{align}
    \omega(\text{det}_q^{-1})=\text{det}_q^{-1}.
\end{align}
\begin{theorem}\label{thm:3.10}
    On $\mathcal{O}(GL_q(n))$, the following statements are true.
    \begin{enumerate}
        \item[1] $h\circ\gamma$ is identical with the Haar state $h$.
        \item[2] $h\circ\omega$ is identical with the Haar state $h$.
    \end{enumerate}
\end{theorem}
\begin{proof}
   By the uniqueness of the Haar state, it suffices to show that $h\circ\gamma$ and $h\circ\omega$ are states on $\mathcal{O}(GL_q(n))$ and possess the translation invariant property. In the proof, we will write $\Delta(x)=\sum x_{(1)}\otimes x_{(2)}$. Notice that $\gamma^2=id$ and $\gamma(1)=1$, and $\omega^2=id$ and $\omega(1)=1$.
    \begin{enumerate}
        \item[1] By the right translation invariant property of $h$, we have
        \begin{equation}\label{eq:3.61}
            \begin{split}
            h(\gamma(x))\cdot 1=&(id\otimes h)\circ\Delta\circ(\gamma(x))\\
            =&(id\otimes h)\circ\iota\circ(\gamma\otimes\gamma)\circ\Delta(x)\\
            =&(h\otimes id)\circ(\gamma\otimes\gamma)\circ\Delta(x)\\
            =&\sum h(\gamma(x_{(1)}))\cdot \gamma(x_{(2)}).
            \end{split}
        \end{equation}
        Applying $\gamma$ on both sides of (\ref{eq:3.61}), we get
        \begin{align}\label{eq:3.62}
            h(\gamma(x))\cdot \gamma(1)=\sum h(\gamma(x_{(1)}))\cdot \gamma^2(x_{(2)})=\sum h(\gamma(x_{(1)}))\cdot x_{(2)}.
        \end{align}
        Since $\gamma(1)=1$, (\ref{eq:3.62}) implies that
        \begin{align}
            (h\circ\gamma(x))\cdot 1=((h\circ\gamma)\otimes id)\circ\Delta(x)
        \end{align}
        which means that $h\circ\gamma$ possess the left translation invariant property. Similarly, the left translation invariant property of $h$ implies the right translation invariant property of $h\circ \gamma$. We also have $h\circ\gamma(1)=h(1)=1$. Hence, the uniqueness of $h$ implies that $h\circ\gamma=h$.
        \item[2] By the right translation invariant property of $h$, we have
        \begin{equation}\label{eq:3.64}
            \begin{split}
            h(\omega(x))\cdot 1=&(id\otimes h)\circ\Delta\circ(\omega(x))\\
            =&(id\otimes h)\circ(\omega\otimes\omega)\circ\Delta(x)\\
            =&\sum h(\omega(x_{(2)}))\cdot \omega(x_{(1)}).
            \end{split}
        \end{equation}
        Applying $\omega$ on both sides of (\ref{eq:3.64}), we conclude that
        \begin{equation}
            \begin{split}
                h\circ\omega(x)\cdot 1&=h(\omega(x))\cdot \omega(1)=\sum h(\omega(x_{(2)}))\cdot \omega^2(x_{(1)})\\
                &=\sum h(\omega(x_{(2)}))\cdot x_{(1)}=(id\otimes h\circ\omega)\circ\Delta(x).
            \end{split}
        \end{equation}
        Hence, the right translation invariant property of $h$ implies that of $h\circ\omega$. The left translation invariant property of $h\circ\omega$ is proved similarly. We also have $h\circ\omega(1)=h(1)=1$. Hence, the uniqueness of $h$ implies that $h\circ\omega=h$.
    \end{enumerate}
\end{proof}

%% file: Chapters/Chapter_4.tex
\section{The Haar state values of monomials on $\mathcal{O}(U_q(3))$}
In this section we will pick a set of psudo-basis of order $m$ for each $m\in\mathbb{N}$ and give the Haar state values of our chosen monomials as finite summations of basic-hypergeometric terms in variable $q$. All results obtained in this section is applicable to the compact quantum group $SU_q(3)$ if we identify the generator $\text{det}_q^{-1}$ as $1$ and set $D_q=1$. The computation splits into two independent parts:
\begin{enumerate}
    \item[1)] finding $h(c^me^mg^m)$ by constructing and solving the Source matrix of order $m$,
    \item[2)] applying the invariant property of the Haar state under the action of $\mathcal{U}_q(gl_3)$ on $\mathcal{O}(U_q(3))$ to find the Haar state value of the general form (\ref{eq:4.4}).
\end{enumerate}
For simplicity, we denote the generators of $\mathcal{O}(Mat_q(3))$ as
\begin{equation}\label{mat:generator_3}
    \begin{bmatrix}
    a&b&c\\
    d&e&f\\  
    g&h&k
    \end{bmatrix}.
\end{equation}
Notice that all matrices $M\in B_3(m)$ can be written as
\begin{equation}\label{mat:4.1}
    \begin{bmatrix}
    s&m-s-r&r\\
    m-s-l&s+r+l+t-m&m-r-t\\  
    l&m-l-t&t
    \end{bmatrix}
\end{equation}
for integers $s,r,l,t$ such that every entry in (\ref{mat:4.1}) is non-negative, i.e.,
\begin{align}\stepcounter{equation}
    s,r,l,t&\ge 0, \tag{4.3.a}\label{eq:4.3.a}\\
    m-s-r,m-s-l,m-r-t,m-l-t&\ge 0,\tag{4.3.b}\label{eq:4.3.b}\\
    s+r+l+t-m&\ge 0.\tag{4.3.c}\label{eq:4.3.c}
\end{align}
Hence, the monomials we picked as a set of psudo-basis of order $m$ are
\begin{align}\label{eq:4.4}
    a^sb^{m-s-r}c^rd^{m-s-l}e^{s+r+l+t-m}f^{m-r-t}g^lh^{m-l-t}k^t\cdot\text{det}_q^{-m}.
\end{align}
For simplicity, the Haar states of monomials in the form of (\ref{eq:4.4}) are denoted as
\begin{align}\label{eq:4.5}
    h(m;s,r,l,t).
\end{align}

\input{Chapters/Chapter_4/New_Source_matrix}

\input{Chapters/Chapter_4/Haar_state_of_others}

\input{Chapters/Chapter_4/Haar_state_general}

\input{Chapters/Chapter_4/Generalization}

%% file: Chapters/Chapter_4/New_Source_matrix.tex
\subsection{The Source matrix of order $m$ on $\mathcal{O}(U_q(3))$}
In this subsection we construct the Source matrix of order $m$ on $\mathcal{O}(U_q(3))$. We will apply a method different from the description in Section 3.4. In our opinion, this construction is more efficient on $\mathcal{O}(U_q(3))$. Currently, it is unclear how to extend this construction into higher rank case. 

The unknowns of Source matrix of order $m$ are the Haar state values of 
\begin{multicols}{2}
    \begin{enumerate}
        \item[1.] $(aek)(ceg)^{m-1}\cdot\text{det}_q^{-m}$,
        \item[2.] $(afh)(ceg)^{m-1}\cdot\text{det}_q^{-m}$,
        \item[3.] $(bdk)(ceg)^{m-1}\cdot\text{det}_q^{-m}$,
        \item[4.]  $(bfg)(ceg)^{m-1}\cdot\text{det}_q^{-m}$,
        \item[5.] $(cdh)(ceg)^{m-1}\cdot\text{det}_q^{-m}$,
        \item[6.] $(ceg)^{m}\cdot\text{det}_q^{-m}$.
    \end{enumerate}
\end{multicols}
\begin{center}
    Table 4.1
\end{center}
To construct linear relations between the Haar state values of these monomials, we will apply the invariant property of the Haar state under the action of $\mathcal{U}_q(gl_3)$ on $\mathcal{O}(U_q(3))$ in (\ref{eq:2.38}) and (\ref{eq:2.38+}). In our computation, we will omit the coefficient generated from $(q^{\lambda},\text{det}_q^{-m})$ (See (\ref{eq:2.15.c})) for simplicity.\\
1. Consider the right action of $e_1\in\mathcal{U}_q(gl_3)$ on monomial $bc^{m-1}e^{m}g^m\cdot\text{det}_q^{-m}$. We have
\begin{equation}\label{eq:4.6}
    \begin{split}
        &bc^{m}e^{m-1}g^m\cdot\text{det}_q^{-m}\cdot e_1\\
        =&e(q^{-1/2}c)^{m}(q^{1/2}e)^{m-1}g^m\cdot\text{det}_q^{-m}\\
        &+\sum_{z=1}^{m}(q^{1/2}b)(q^{1/2}c)^{z-1}f(q^{-1/2}c)^{m-z}(q^{1/2}e)^{m-1}g^m\cdot\text{det}_q^{-m}\\
        =&q^{-1/2}(ceg)^m\cdot\text{det}_q^{-m}+q^{-3/2}\frac{q^2-q^{2m+2}}{1-q^2}(bfg)(ceg)^{m-1}\cdot\text{det}_q^{-m}.
    \end{split}
\end{equation}
Evaluating the Haar state values on both sides of (\ref{eq:4.6}) and applying the invariant property (\ref{eq:2.38+}), we have
\begin{equation}\label{eq:4.7}
    h((bfg)(ceg)^{m-1}\cdot\text{det}_q^{-m})=\frac{-q^{-1}(1-q^2)}{1-q^{2m}}h((ceg)^m\cdot\text{det}_q^{-m}).
\end{equation}
2. Recall that the double flip anti-automorphism $\omega$ (\ref{eq:3.52}). Since this map preserve the Haar state value, we have
\begin{equation}
\begin{split}
    &h((bfg)(ceg)^{m-1}\cdot\text{det}_q^{-m})=h\circ\omega((bfg)(ceg)^{m-1}\cdot\text{det}_q^{-m})\\
    =&h(\text{det}_q^{-m}\cdot(ceg)^{m-1}(cdh))=h((cdh)(ceg)^{m-1}\cdot\text{det}_q^{-m}).
\end{split}
\end{equation}
Hence, we get
\begin{equation}\label{eq:4.9}
    h((cdh)(ceg)^{m-1}\cdot\text{det}_q^{-m})=\frac{-q^{-1}(1-q^2)}{1-q^{2m}}h((ceg)^m\cdot\text{det}_q^{-m}).
\end{equation}
3. Consider the right action of $f_2\in\mathcal{U}_q(gl_3)$ on monomial $bkc^{m-1}e^{m-1}g^{m}\cdot\text{det}_q^{-m}$. We have 
\begin{equation}\label{eq:4.10}
    \begin{split}
        &bkc^{m-1}e^{m-1}g^{m}\cdot\text{det}_q^{-m}\cdot f_2\\
        =&bfc^{m-1}(q^{-1/2}e)^{m-1}(q^{1/2}g)^{m}\cdot\text{det}_q^{-m}\\
        &+\sum_{z=1}^mb(q^{-1/2}k)c^{m-1}(q^{1/2}e)^{m-1}(q^{-1/2}g)^{z-1}d(q^{1/2}g)^{m-z}\cdot\text{det}_q^{-m}\\
        =&q^{1/2}bfc^{m-1}e^{m-1}g^{m}\cdot\text{det}_q^{-m}\\
        &+q^{2m-1/2}\frac{q^{-2}-q^{-2(m+1)}}{1-q^{-2}}bkc^{m-1}e^{m-1}g^{m-1}d\cdot\text{det}_q^{-m}.
    \end{split}
\end{equation}
Evaluating the Haar state values on both sides of (\ref{eq:4.10}) and applying the invariant property (\ref{eq:2.38+}), we have
\begin{equation}\label{eq:4.11}
    \begin{split}
        h(bkc^{m-1}e^{m-1}g^{m-1}d\cdot\text{det}_q^{-m})=-\frac{q^{3}-q}{q^{2m}-1}h((bfg)(ceg)^{m-1}\cdot\text{det}_q^{-m}).
    \end{split}
\end{equation}
On the other hand, we have
\begin{equation}\label{eq:4.12}
\begin{split}
    &h(bkc^{m-1}e^{m-1}g^{m-1}d\cdot\text{det}_q^{-m})=h(\rho(d)bkc^{m-1}e^{m-1}g^{m-1}\cdot\text{det}_q^{-m})\\
    &=q^2h(bdk(ceg)^{m-1}\cdot\text{det}_q^{-m}).
\end{split}
\end{equation}
Combining (\ref{eq:4.11}) and (\ref{eq:4.12}), we find that
\begin{equation}\label{eq:4.13}
    \begin{split}
        h(bdk(ceg)^{m-1}\cdot\text{det}_q^{-m})=\frac{q^{-2}(1-q^2)^2}{(1-q^{2m})^2}h((ceg)^m\cdot\text{det}_q^{-m}).
    \end{split}
\end{equation}
4. Using the double flip anti-automorphism $\omega$ and the modular automorphism $\rho$, we have
\begin{equation}\label{eq:4.14}
    \begin{split}
        &h(bdk(ceg)^{m-1}\cdot\text{det}_q^{-m})=h\circ\omega(bdk(ceg)^{m-1}\cdot\text{det}_q^{-m})\\
        &=h(\text{det}_q^{-m}\cdot(ceg)^{m-1}afh)=h(\rho(afh)(ceg)^{m-1}\cdot\text{det}_q^{-m})\\
        &=h(afh(ceg)^{m-1}\cdot\text{det}_q^{-m}).
    \end{split}
\end{equation}
Combining (\ref{eq:4.13}) and (\ref{eq:4.14}), we find that 
\begin{equation}\label{eq:4.15}
    \begin{split}
        h(afh(ceg)^{m-1}\cdot\text{det}_q^{-m})=\frac{q^{-2}(1-q^2)^2}{(1-q^{2m})^2}h((ceg)^m\cdot\text{det}_q^{-m}).
    \end{split}
\end{equation}
5. Consider the left action of  $e_1\mathcal{U}_q(gl_3)$ on monomial $bkc^{m-1}e^{m}g^{m-1}\cdot\text{det}_q^{-m}$. We have
\begin{equation}\label{eq:4.16}
    \begin{split}
        &e_1\cdot bkc^{m-1}e^{m}g^{m-1}\cdot\text{det}_q^{-m}\\
        =&akc^{m-1}(q^{1/2}e)^{m}(q^{-1/2}g)^{m-1}\cdot\text{det}_q^{-m}\\
        &+\sum_{z=1}^m(q^{-1/2}b)kc^{m-1}(q^{-1/2}e)^{z-1}d(q^{1/2}e)^{m-z}(q^{-1/2}g)^{m-1}\cdot\text{det}_q^{-m}\\
        =&q^{1/2}akc^{m-1}e^{m}g^{m-1}\cdot\text{det}_q^{-m}\\
        &+q^{2m-1/2}\frac{q^{-2}-q^{-2(m+1)}}{1-q^{-2}}bkc^{m-1}e^{m-1}g^{m-1}d\cdot\text{det}_q^{-m}\\
    \end{split}
\end{equation}
Evaluating the Haar state values on both sides of (\ref{eq:4.16}) and applying the invariant property (\ref{eq:2.38}), we have
\begin{equation}\label{eq:4.17}
    \begin{split}
        h(akc^{m-1}e^{m}g^{m-1}\cdot\text{det}_q^{-m})=\frac{1-q^{2m}}{q^{3}-q}h(bkc^{m-1}e^{m-1}g^{m-1}d\cdot\text{det}_q^{-m})
    \end{split}
\end{equation}
Applying the modular automorphism $\rho$, we get
\begin{equation}\label{eq:4.18}
\begin{split}
    &h(akc^{m-1}e^{m}g^{m-1}\cdot\text{det}_q^{-m})\\
    =&h(\rho(e)akc^{m-1}e^{m-1}g^{m-1}\cdot\text{det}_q^{-m})=h(eakc^{m-1}e^{m-1}g^{m-1}\cdot\text{det}_q^{-m})\\
    =&h(aekc^{m-1}e^{m-1}g^{m-1}\cdot\text{det}_q^{-m})\\
    &-(q-q^{-1})h(bdkc^{m-1}e^{m-1}g^{m-1}\cdot\text{det}_q^{-m}).
\end{split}
\end{equation}
Combining (\ref{eq:4.13}), (\ref{eq:4.17}) and (\ref{eq:4.18}), we get
\begin{equation}\label{eq:4.19}
\begin{split}
    &h(aekc^{m-1}e^{m-1}g^{m-1}\cdot\text{det}_q^{-m})\\
    &=-\left(\frac{(1-q^2)^2}{1-q^{2m}}+q^2\right)\frac{q^{-3}(1-q^2)}{(1-q^{2m})}h((ceg)^m\cdot\text{det}_q^{-m}).
\end{split}
\end{equation}
The linear relation with non-zero right-hand-side is
\begin{equation}\label{eq:4.20}
    h(D_q\cdot(ceg)^{m-1}\cdot\text{det}_q^{-m})=h((ceg)^{m-1}\cdot\text{det}_q^{-m+1}).
\end{equation}
Substituting (\ref{eq:4.6}), (\ref{eq:4.9}), (\ref{eq:4.13}), (\ref{eq:4.15}) and (\ref{eq:4.19}) into (\ref{eq:4.20}), we find that
\begin{equation}\label{eq:4.21}
    h\left((ceg)^{m}\cdot\text{det}_q^{-m}\right)=\frac{-q^3(q^{2m}-1)^2\cdot h\left((ceg)^{m-1}\cdot\text{det}_q^{-m+1}\right)}{(1-q^{2m+2})(1-q^{2m+4})}.
\end{equation}
Recall that by (\ref{eq:3.30}) we have
\begin{align}\label{eq:4.36}
    h(ceg\cdot\text{det}_q^{-1})=\frac{(-q)^3(1-q^2)^2}{(q^4-1)(q^6-1)}.
\end{align}
Hence, by (\ref{eq:4.21}) and (\ref{eq:4.36}), we find that
\begin{equation}\label{eq:4.37}
    h\left(c^me^mg^m\cdot\text{det}_q^{-m}\right)=\frac{(-q)^{3m}(q^2-1)^2(q^4-1)}{(q^{2m+2}-1)^2(q^{2m+4}-1)}.
\end{equation}

%% file: Chapters/Chapter_4/Haar_state_of_others.tex
\subsection{The Haar states of some special psudo-basis vectors}
In this subsection, we compute the Haar state values of some special monomials in our psudo-basis that will serve as the initial and boundary conditions of the recursive relation used to find the general form of $h(m;s,r,l,t)$.

\input{Chapters/Chapter_4/Haar_states_special/Haar_1}

\input{Chapters/Chapter_4/Haar_states_special/Haar_2}

\input{Chapters/Chapter_4/Haar_states_special/Haar_3}
\input{Chapters/Chapter_4/Haar_states_special/Haar_4}

\input{Chapters/Chapter_4/Haar_states_special/Haar_5}

%% file: Chapters/Chapter_4/Haar_states_special/Haar_1.tex
\subsubsection{The value of $h(m;0,m,l,0)$}
The monomial corresponds to $h(m;0,m,l,0)$ is
\begin{align}
    c^md^{m-l}e^lg^lh^{m-l}\cdot\text{det}_q^{-m}.
\end{align}
Consider the left action of $f_1\in\mathcal{U}_q(gl_3)$ on the monomial $c^md^{m-l}e^lg^{l+1}h^{m-l-1}$. We have
\begin{equation}\label{eq:4.40}
    \begin{split}
        &f_1\cdot c^md^{m-l}e^lg^{l+1}h^{m-l-1}\cdot\text{det}_q^{-m}\\
        =&q^{m-l+1/2}\frac{q^{-2}-q^{-2(m-l+1)}}{1-q^{-2}}c^md^{m-l-1}e^{l+1}g^{l+1}h^{m-l-1}\cdot\text{det}_q^{-m}\\
        &+q^{m-l+3/2}\frac{q^{-2}-q^{-2(l+2)}}{1-q^{-2}}c^md^{m-l}e^lg^{l}h^{m-l}\cdot\text{det}_q^{-m}.
    \end{split}
\end{equation}
Evaluating the Haar state on both sides of (\ref{eq:4.40}) and applying invariant property (\ref{eq:2.38}), we get
\begin{equation}\label{eq:4.41.b}
    \begin{split}
        \Longleftrightarrow &h\left(c^md^{m-l}e^lg^{l}h^{m-l}\cdot\text{det}_q^{-m}\right)\\
        &=-q^{1-2m+4l}\frac{q^{2(m-l)}-1}{q^{2(l+1)}-1}h\left(c^md^{m-l-1}e^{l+1}g^{l+1}h^{m-l-1}\cdot\text{det}_q^{-m}\right).
    \end{split}
\end{equation}
Applying (\ref{eq:4.41.b}) repeatedly, we get
\begin{equation}\label{eq:4.42}
\begin{split}
    h\left(c^md^{m-l}e^lg^{l}h^{m-l}\cdot\text{det}_q^{-m}\right)=\frac{(-q)^{(m-l)(2l-1)}}{{m \choose l}_{q^2}}\cdot h\left(c^me^{m}g^{m}\right).
\end{split}
\end{equation}

%% file: Chapters/Chapter_4/Haar_states_special/Haar_2.tex
\subsubsection{The value of $h(m;0,r,l,0)$}
The monomial corresponds to $h(m;0,r,l,0)$ is
\begin{align}
     b^{m-r}c^{r}d^{m-l}e^{r-m+l}f^{m-r}g^{l}h^{m-l}\cdot\text{det}_q^{-m}.
\end{align}
Consider the left action of $e_2\in\mathcal{U}_q(gl_3)$ on $b^{m-r-1}c^{r+1}d^{m-l}e^{r-m+l}f^{m-r}g^{l}h^{m-l}\cdot\text{det}_q^{-m}$. We have
\begin{equation}\label{eq:4.44}
    \begin{split}
        &e_2\cdot b^{m-r-1}c^{r+1}d^{m-l}e^{r-m+l}f^{m-r}g^{l}h^{m-l}\cdot\text{det}_q^{-m}\\
        =&q^{m-3r-5/2}\frac{q^2-q^{2r+4}}{1-q^2}b^{m-r}c^{r}d^{m-l}e^{r-m+l}f^{m-r}g^{l}h^{m-l}\cdot\text{det}_q^{-m}\\
        &+q^{-3/2+l+r-2m}\frac{q^2-q^{2(m-r+1)}}{1-q^2}\\
        &\cdot b^{m-r-1}c^{r+1}d^{m-l}e^{r+1-m+l}f^{m-r-1}g^{l}h^{m-l}\cdot\text{det}_q^{-m}\\
    \end{split}
\end{equation}
Evaluating the Haar state on both sides of (\ref{eq:4.44}) and applying invariant property (\ref{eq:2.38}), we get
\begin{equation}\label{eq:4.45}
    \begin{split}
        &h\left(b^{m-r}c^{r}d^{m-l}e^{r-m+l}f^{m-r}g^{l}h^{m-l}\cdot\text{det}_q^{-m}\right)\\
        =&-q^{1+l+4r-3m}\frac{1-q^{2(m-r)}}{1-q^{2r+2}}\\
        &\cdot h\left(b^{m-r-1}c^{r+1}d^{m-l}e^{r+1-m+l}f^{m-r-1}g^{l}h^{m-l}\cdot\text{det}_q^{-m}\right).
    \end{split}
\end{equation}
Applying (\ref{eq:4.45}) repeatedly, we get
\begin{equation}\label{eq:4.46}
    \begin{split}
        &h\left(b^{m-r}c^{r}d^{m-l}e^{r-m+l}f^{m-r}g^{l}h^{m-l}\cdot\text{det}_q^{-m}\right)\\
        =&(-1)^{m-r}\frac{q^{(m-r)(l+2r-1-m)}}{{m\choose r}_{q^2}}\cdot h\left(c^{m}d^{m-l}e^{l}g^{l}h^{m-l}\cdot\text{det}_q^{-m}\right).
    \end{split}
\end{equation}
Substituting (\ref{eq:4.42}) into (\ref{eq:4.46}), we get
\begin{equation}\label{eq:4.47}
    \begin{split}
        &h\left(b^{m-r}c^{r}d^{m-l}e^{r-m+l}f^{m-r}g^{l}h^{m-l}\cdot\text{det}_q^{-m}\right)\\
        =&(-1)^{r+l}\frac{q^{(3m+1)(l+r)-2l^2-2r^2-rl-m^2-2m}}{{m\choose r}_{q^2}{m \choose l}_{q^2}}\cdot h\left(c^{m}e^{m}g^{m}\cdot\text{det}_q^{-m}\right)
    \end{split}
\end{equation}

%% file: Chapters/Chapter_4/Haar_states_special/Haar_3.tex
\subsubsection{The value of $h(m;s,r,l,0)$}
The monomial corresponds to $h(m;s,r,l,0)$ is
\begin{equation}\label{eq:4.48}
    a^sb^{m-r-s}c^{r}d^{m-s-l}e^{l+r+s-m}f^{m-r}g^{l}h^{m-l}\cdot\text{det}_q^{-m}.
\end{equation}
First, we assume that $l+r+s=m$, i.e., there is no generator $e$ in (\ref{eq:4.48}). Consider the left action of $e_1\in\mathcal{U}_q(gl_3)$ on the monomial $a^{s-1}b^{m-r-s+1}c^{r}d^{m-s-l}f^{m-r}g^{l}h^{m-l}\cdot\text{det}_q^{-m}$. We have
\begin{equation}\label{eq:4.49}
    \begin{split}
        &e_1\cdot a^{s-1}b^{m-r-s+1}c^{r}d^{m-s-l}f^{m-r}g^{l}h^{m-l}\cdot\text{det}_q^{-m}\\
        =&q^{\frac{5s-5-3m+3r-l}{2}}\frac{q^2-q^{2(m-r-s+2)}}{1-q^2}\\
        &\cdot a^{s}b^{m-r-s}c^{r}d^{m-s-l}f^{m-r}g^{l}h^{m-l}\cdot\text{det}_q^{-m}\\
        &+q^{\frac{s-3+r-3m+3l}{2}}\frac{q^2-q^{2m-2l+2}}{1-q^2}\\
        &\cdot a^{s-1}b^{m-r-s+1}c^{r}d^{m-s-l}f^{m-r}g^{l+1}h^{m-l-1}\cdot\text{det}_q^{-m}\\
    \end{split}
\end{equation}
Evaluating the Haar state on both sides of (\ref{eq:4.49}) and applying invariant property (\ref{eq:2.38}), we get
\begin{equation}\label{eq:4.50}
    \begin{split}
        &h\left(a^{s}b^{m-r-s}c^{r}d^{m-s-l}f^{m-r}g^{l}h^{m-l}\cdot\text{det}_q^{-m}\right)\\
        =&-q^{-2m+r+4l+1}\frac{1-q^{2m-2l}}{1-q^{2(l+1)}}\\
        &h\left(a^{s-1}b^{m-r-s+1}c^{r}d^{m-s-l}f^{m-r}g^{l+1}h^{m-l-1}\cdot\text{det}_q^{-m}\right).
    \end{split}
\end{equation}
Notice that in (\ref{eq:4.50}) the index of generator $c$ is fixed. Applying (\ref{eq:4.50}) repeatedly, we find that
\begin{equation}\label{eq:4.51}
    \begin{split}
        &h\left(a^{s}b^{m-r-s}c^{r}d^{m-s-l}f^{m-r}g^{l}h^{m-l}\cdot\text{det}_q^{-m}\right)\\
        =&(-1)^{m-r-l}q^{(m-r-l)(2l-1-r)}\frac{{m\choose r}_{q^2}}{{m\choose l}_{q^2}}h\left(b^{m-r}c^{r}d^{r}f^{m-r}g^{m-r}h^{r}\cdot\text{det}_q^{-m}\right)\\
    \end{split}
\end{equation}
Substituting (\ref{eq:4.47}) into (\ref{eq:4.51}), we find that
\begin{equation}\label{eq:4.52}
    \begin{split}
        &h\left(a^{s}b^{m-r-s}c^{r}d^{m-s-l}f^{m-r}g^{l}h^{m-l}\cdot\text{det}_q^{-m}\right)\\
        =&\frac{(-1)^{r+l}q^{(2m+1)(l+r)-2m-2l^2-2r^2-rl}}{{m\choose r}_{q^2}{m\choose l}_{q^2}}h\left(c^me^mg^m\cdot\text{det}_q^{-m}\right)
    \end{split}
\end{equation}
Then, we consider the case $l+r+s-m\ge 1$. Applying the left action of $e_2\in\mathcal{U}_q(gl_3)$ on the monomial $a^sb^{m-r-s}c^{r}d^{m-s-l}e^{l+r+s-m-1}f^{m-r+1}g^{l}h^{m-l}\cdot\text{det}_q^{-m}$, we get
\begin{equation}\label{eq:4.49+}
    \begin{split}
        &e_2\cdot a^sb^{m-r-s}c^{r}d^{m-s-l}e^{l+r+s-m-1}f^{m-r+1}g^{l}h^{m-l}\cdot\text{det}_q^{-m}\\
        =&q^{-3r-s+m+1/2}\frac{q^2-q^{2r+2}}{1-q^2}\\
        &\cdot a^sb^{m-r-s+1}c^{r-1}d^{m-s-l}e^{l+r+s-m-1}f^{m-r+1}g^{l}h^{m-l}\cdot\text{det}_q^{-m}\\
        &+q^{l-5/2-2m+r}\frac{q^2-q^{2m-2r+4}}{1-q^2}\\
        &\cdot a^sb^{m-r-s}c^{r}d^{m-s-l}e^{l+r+s-m}f^{m-r}g^{l}h^{m-l}\cdot\text{det}_q^{-m}\\
    \end{split}
\end{equation}
Evaluating the Haar state on both sides of (\ref{eq:4.49+}) and applying invariant property (\ref{eq:2.38}), we get
\begin{equation}\label{eq:4.50+}
    \begin{split}
        &h\left(a^sb^{m-r-s}c^{r}d^{m-s-l}e^{l+r+s-m}f^{m-r}g^{l}h^{m-l}\cdot\text{det}_q^{-m}\right)\\
        =&-q^{-4r-s+3-l+3m}\frac{1-q^{2r}}{1-q^{2m-2r+2}}\\
        &\cdot h\left(a^sb^{m-r-s+1}c^{r-1}d^{m-s-l}e^{l+r+s-m-1}f^{m-r+1}g^{l}h^{m-l}\cdot\text{det}_q^{-m}\right)
    \end{split}
\end{equation}
By (\ref{eq:4.3.b}), we must have $r\ge l+r+s-m$. Applying (\ref{eq:4.50+}) for $l+r+s-m$ times, we get
\begin{equation}\label{eq:4.51+}
    \begin{split}
        &h\left(a^sb^{m-r-s}c^{r}d^{m-s-l}e^{l+r+s-m}f^{m-r}g^{l}h^{m-l}\cdot\text{det}_q^{-m}\right)\\
        =&(-1)^{l+r+s-m}q^{(m-2r+l+s+1)(l+r+s-m)}\frac{\prod_{i=m-s-l+1}^{r}(1-q^{2i})}{\prod_{j=m-r+1}^{l+s}(1-q^{2j})}\\
        &h\left(a^sb^{l}c^{m-s-l}d^{m-s-l}f^{l+s}g^{l}h^{m-l}\cdot\text{det}_q^{-m}\right)
    \end{split}
\end{equation}
Substituting (\ref{eq:4.48}) into (\ref{eq:4.51+}), we get
\begin{equation}\label{eq:4.57}
    \begin{split}
        &h(a^sb^{m-r-s}c^{r}d^{m-s-l}e^{l+r+s-m}f^{m-r}g^{l}h^{m-l}\cdot\text{det}_q^{-m})\\
        =&\frac{(-1)^{r+l}q^{(2m+1)(l+r)-2l^2-2r^2-rl-2m+(m-s)(l+r+s-m)}}{{m\choose r}_{q^2}{m\choose l}_{q^2}}\\
        &\cdot h\left(c^me^mg^m\cdot\text{det}_q^{-m}\right).
    \end{split}
\end{equation}

%% file: Chapters/Chapter_4/Haar_states_special/Haar_4.tex
\subsubsection{The value of $h(m;s,r,l,t)$ with $s+r+l+t=m$}
The monomial corresponds to $h(m;s,r,l,t)$ with $s+r+l+t=m$ is
\begin{equation}
    a^sb^{m-s-r}c^rd^{m-s-l}f^{m-r-t}g^lh^{m-l-t}k^t\cdot\text{det}_q^{-m}.
\end{equation}
Consider the left action of $f_2\in\mathcal{U}_q(gl_3)$ on the monomial\\ $a^sb^{m-s-r}c^rd^{m-s-l}f^{m-r-t}g^lh^{m-l-t+1}k^{t-1}\cdot\text{det}_q^{-m}$. We get
\begin{equation}\label{eq:4.59}
    \begin{split}
        &f_2\cdot a^sb^{m-s-r}c^rd^{m-s-l}f^{m-r-t}g^lh^{m-l-t+1}k^{t-1}\cdot\text{det}_q^{-m}\\
        =&q^{\frac{m-s+1-r+l+t}{2}}\frac{q^{-2}-q^{-2(m-s-r+1)}}{1-q^{-2}}\\
        &\cdot a^sb^{m-s-r-1}c^{r+1}d^{m-s-l}f^{m-r-t}g^lh^{m-l-t+1}k^{t-1}\cdot\text{det}_q^{-m}\\
        &+q^{\frac{m-s-r-l+3+t}{2}}\frac{q^{-2}-q^{-2(m-l-t+2)}}{1-q^{-2}}\\
        &\cdot a^sb^{m-s-r}c^rd^{m-s-l}f^{m-r-t}g^lh^{m-l-t}k^{t}\cdot\text{det}_q^{-m}\\
    \end{split}
\end{equation}
Evaluating the Haar state on both sides of (\ref{eq:4.59}) and applying invariant property (\ref{eq:2.38}), we get
\begin{equation}\label{eq:4.60}
    \begin{split}
        &h\left(a^sb^{m-s-r}c^rd^{m-s-l}f^{m-r-t}g^lh^{m-l-t}k^{t}\cdot\text{det}_q^{-m}\right)\\
        =&-q^{2m-3l-4t+1}\frac{1-q^{2(l+t)}}{1-q^{2(m-l-t+1)}}\\
        &h\left(a^sb^{m-s-r-1}c^{r+1}d^{m-s-l}f^{m-r-t}g^lh^{m-l-t+1}k^{t-1}\cdot\text{det}_q^{-m}\right).
    \end{split}
\end{equation}
Notice that the index of generator $a$ and generator $g$ is fixed. Applying (\ref{eq:4.60}) repeatedly, we get
\begin{equation}\label{eq:4.61}
    \begin{split}
        &h\left(a^sb^{m-s-r}c^rd^{m-s-l}f^{m-r-t}g^lh^{m-l-t}k^{t}\cdot\text{det}_q^{-m}\right)\\
        =&(-1)^tq^{t(2m-3l-2t-1)}\prod_{i=1}^t\frac{1-q^{2(l+i+1)}}{1-q^{2(m-l-i)}}\\
        &\cdot h\left(a^sb^{m-s-r-t}c^{r+t}d^{m-s-l}f^{m-r-t}g^lh^{m-l}\cdot\text{det}_q^{-m}\right).
    \end{split}
\end{equation}
Substituting (\ref{eq:4.52}) into (\ref{eq:4.61}), we get
\begin{equation}\label{eq:4.63}
    \begin{split}
        &h\left(a^sb^{m-s-r}c^rd^{m-s-l}f^{m-r-t}g^lh^{m-l-t}k^{t}\cdot\text{det}_q^{-m}\right)\\
        =&\frac{(-1)^{r+l}q^{4st+(2m+1)(l+r)-2m-2l^2-2r^2-rl}}{{m\choose r+t}_{q^2}{m\choose l+t}_{q^2}}h\left(c^me^mg^m\cdot\text{det}_q^{-m}\right).
    \end{split}
\end{equation}

%% file: Chapters/Chapter_4/Haar_states_special/Haar_5.tex
\subsubsection{The value of $h(m;s,0,l,t)$}
The monomial corresponds to $h(m;s,0,l,t)$ is
\begin{equation}
    a^{s}b^{m-s}d^{m-s-l}e^{s+l+t-m}f^{m-t}g^{l}h^{m-l-t}k^{t}\cdot\text{det}_q^{-m}.
\end{equation}
Consider the left action of $e_2\in\mathcal{U}_q(gl_3)$ on the monomial\\ $a^{s}b^{m-s}d^{m-s-l}e^{s+l+t-m-1}f^{m-t+1}g^{l}h^{m-l-t}k^{t}$. We have
\begin{equation}\label{eq:4.65}
    \begin{split}
        &e_2\cdot a^{s}b^{m-s}d^{m-s-l}e^{s+l+t-m-1}f^{m-t+1}g^{l}h^{m-l-t}k^{t}\cdot\text{det}_q^{-m}\\
        =&q^{l+3t-2m-5/2}\frac{q^2-q^{2(m-t+2)}}{1-q^2}a^{s}b^{m-s}d^{m-s-l}e^{s+l+t-m}f^{m-t}g^{l}h^{m-l-t}k^{t}\cdot\text{det}_q^{-m}\\
        +&q^{-t-3/2}\frac{q^2-q^{2t+2}}{1-q^2}a^{s}b^{m-s}d^{m-s-l}e^{s+l+t-m-1}f^{m-t+1}g^{l}h^{m-l-t+1}k^{t-1}\cdot\text{det}_q^{-m}\\
    \end{split}
\end{equation}
Evaluating the Haar state on both sides of (\ref{eq:4.65}) and applying invariant property (\ref{eq:2.38}), we get
\begin{equation}\label{eq:4.66}
    \begin{split}
        &h\left(a^{s}b^{m-s}d^{m-s-l}e^{s+l+t-m}f^{m-t}g^{l}h^{m-l-t}k^{t}\cdot\text{det}_q^{-m}\right)\\
        =&-q^{-l-4t+2m+1}\frac{1-q^{2t}}{1-q^{2(m-t+1)}}\\
        &\cdot h\left(a^{s}b^{m-s}d^{m-s-l}e^{s+l+t-m-1}f^{m-t+1}g^{l}h^{m-l-t+1}k^{t-1}\cdot\text{det}_q^{-m}\right).
    \end{split}
\end{equation}
Recall (\ref{eq:4.3.b}) that $s+l\le m$. Hence, $s+l+t-m\le t$. We apply (\ref{eq:4.66}) for $s+l+t-m$ times and get
\begin{equation}\label{eq:4.67}
    \begin{split}
        &h\left(a^{s}b^{m-s}d^{m-s-l}e^{s+l+t-m}f^{m-t}g^{l}h^{m-l-t}k^{t}\cdot\text{det}_q^{-m}\right)\\
        =&(-1)^{s+l+t-m}\prod_{i=1}^{s+l+t-m} q^{-l-4(t+1-i)+2m+1}\frac{1-q^{2(t+1-i)}}{1-q^{2(m-(t+1-i)+1)}}\\
        &\cdot h\left(a^{s}b^{m-s}d^{m-s-l}f^{s+l}g^{l}h^{s}k^{m-s-l}\cdot\text{det}_q^{-m}\right)\\
    \end{split}
\end{equation}
Substituting (\ref{eq:4.63}) with $r=0$ into (\ref{eq:4.67}), we get
\begin{equation}\label{eq:4.69}
    \begin{split}
        &h\left(a^{s}b^{m-s}d^{m-s-l}e^{s+l+t-m}f^{m-t}g^{l}h^{m-l-t}k^{t}\cdot\text{det}_q^{-m}\right)\\
        =&(-1)^{m-s-t}
        \frac{q^{4st-(2s+2t+1-l)(s+l+t-m)+(2m+1)l-2m-2l^2}}{{m\choose t}_{q^2}{m\choose s}_{q^2}}\\
        &\cdot h\left(c^me^mg^m\cdot\text{det}_q^{-m}\right).
    \end{split}
\end{equation}

%% file: Chapters/Chapter_4/Haar_state_general.tex
\subsection{The Haar state value of the general form  $h(m;s,r,l,t)$}
To simplify the computation in this subsection, we denote the index of generator $e$ as
\begin{equation}
    n=s+r+l+t-m.
\end{equation}
We want to derive a recursive relation of $h(m;s,r,l,t)$ depending on the value of $n$. Consider the left action of $e_2\in\mathcal{U}_q(gl_3)$ on the monomial\\ $a^sb^{m-s-r}c^rd^{m-s-l}e^{n-1}f^{m-r-t+1}g^lh^{m-l-t}k^t\cdot\text{det}_q^{-m}$. We have
\begin{equation}\label{eq:4.71}
    \begin{split}
        &e_2\cdot a^sb^{m-s-r}c^rd^{m-s-l}e^{n-1}f^{m-r-t+1}g^lh^{m-l-t}k^t\cdot\text{det}_q^{-m}\\
        =&q^{m-s-3r+1/2}\frac{q^2-q^{2r+2}}{1-q^2}\\
        &\cdot a^sb^{m-s-r+1}c^{r-1}d^{m-s-l}e^{n-1}f^{m-r-t+1}g^lh^{m-l-t}k^t\cdot\text{det}_q^{-m}\\
        &+q^{-2m-5/2+l+3t+r}\frac{q^2-q^{2(m-r-t+2)}}{1-q^2}\\
        &\cdot a^sb^{m-s-r}c^rd^{m-s-l}e^{n}f^{m-r-t}g^lh^{m-l-t}k^t\cdot\text{det}_q^{-m}\\
        &+q^{-3/2-t}\frac{q^2-q^{2t+2}}{1-q^2}\\
        &\cdot a^sb^{m-s-r}c^rd^{m-s-l}e^{n-1}f^{m-r-t+1}g^lh^{m-l-t+1}k^{t-1}\cdot\text{det}_q^{-m}.
    \end{split}
\end{equation}
Evaluating the Haar state on both sides of (\ref{eq:4.71}) and applying invariant property (\ref{eq:2.38}), we get
\begin{equation}\label{eq:4.72}
    \begin{split}
        &h\left(a^sb^{m-s-r}c^rd^{m-s-l}e^{n}f^{m-r-t}g^lh^{m-l-t}k^t\cdot\text{det}_q^{-m}\right)\\
        =&-q^{3m-s-l-4r-3t+3}\frac{1-q^{2r}}{1-q^{2(m-r-t+1)}}\\
        &\cdot h\left(a^sb^{m-s-r+1}c^{r-1}d^{m-s-l}e^{n-1}f^{m-r-t+1}g^lh^{m-l-t}k^t\cdot\text{det}_q^{-m}\right)\\
        &-q^{2m-l-4t-r+1}\frac{1-q^{2t}}{1-q^{2(m-r-t+1)}}\\
        &\cdot h\left(a^sb^{m-s-r}c^rd^{m-s-l}e^{n-1}f^{m-r-t+1}g^lh^{m-l-t+1}k^{t-1}\cdot\text{det}_q^{-m}\right).
    \end{split}
\end{equation}
(\ref{eq:4.72}) is a recursive relation on the value of $n$ with initial condition (\ref{eq:4.63}) and boundary condition (\ref{eq:4.57}) and (\ref{eq:4.69}). Together with the solution of the Source matrix of order $m$, we state the main theorem.
\begin{theorem}[Main Theorem]
     Let $m,s,r,l,t$ be a set of integers satisfying condition (4.3.a-c) and set $n=s+r+l+t-m$. We have
     \begin{equation}\tag{4.23}
         h\left(c^me^mg^m\cdot\text{det}_q^{-m}\right)=\frac{(-q)^{3m}(q^2-1)^2(q^4-1)}{(q^{2m+2}-1)^2(q^{2m+4}-1)},
     \end{equation}
     and
    \begin{equation}\label{eq:4.73}
        \begin{split}
            &h(a^sb^{m-s-r}c^rd^{m-s-l}e^{n}f^{m-r-t}g^lh^{m-l-t}k^t\cdot\text{det}_q^{-m})\\
            =&\left(\sum_{k=0}^{n}(-1)^{k}q^{(n-k)(n-3k-1)+2k(s+t)}\frac{{r\choose k}_{q^2}{l\choose k}_{q^2}{s\choose n-k}_{q^2}{t\choose n-k}_{q^2}}{{n\choose k}_{q^2}}\right)\\
            &\frac{(-1)^{r+l+n}q^{(2m+1)(l+r)+(n-2)m-2l^2-2r^2-rl+4st-3ns-3nt}}{{m\choose n,m-l-s,m-r-t}_{q^2}{m\choose n,m-r-s,m-l-t}_{q^2}}\\
            &\cdot h\left(c^me^mg^m\cdot\text{det}_q^{-m}\right).
        \end{split}
    \end{equation}
\end{theorem}
\begin{remark}
    The product ${r\choose k}_{q^2}{l\choose k}_{q^2}{s\choose n-k}_{q^2}{t\choose n-k}_{q^2}$ is non-zero when 
    $$\text{max}(n-s,n-t,0)\le k\le \text{min}(r,l,n).$$
\end{remark}
\begin{proof}[Proof of (\ref{eq:4.73})]
    It is easy to check that our initial condition (\ref{eq:4.63}) and boundary condition (\ref{eq:4.57}) and (\ref{eq:4.69}) satisfies (\ref{eq:4.73}). 

    In the following computation, we will omit the part $h\left(c^me^mg^m\cdot\text{det}_q^{-m}\right)$ for simplicity. Recall the following expression from the theory of basic hypergeometric functions
    \begin{equation}\label{eq:4.71+}
        (a;q)_n=\prod_{i=0}^{n-1}(1-aq^i) \text{ and }(a;q)_0=1.
    \end{equation}
    With this notation we can rewrite (\ref{eq:4.73}) as
    \begin{equation}\label{eq:4.77}
        \begin{split}
            &h\left(a^sb^{m-s-r}c^rd^{m-s-l}e^{n}f^{m-r-t}g^lh^{m-l-t}k^t\cdot\text{det}_q^{-m}\right)\\
            =&\left(\sum_{k=0}^{n}(-1)^{k}q^{(n-k)(n-3k-1)+2k(s+t)}\right.\\
            &\left.\cdot{n\choose k}_{q^{2}}(q^{2r};q^{-2})_k(q^{2l};q^{-2})_k(q^{2s};q^{-2})_{n-k}(q^{2t};q^{-2})_{n-k}\right)\\
            &\frac{(q^2;q^2)_{r+t-n}(q^2;q^2)_{l+t-n}(q^2;q^2)_{r+s-n}(q^2;q^2)_{l+s-n}}{(q^2;q^2)_m^2}\\
            &(-1)^{r+l+n}q^{(2m+1)(l+r)+(n-2)m-2l^2-2r^2-rl+4st-3ns-3nt}.
        \end{split}
    \end{equation}
    Now, assume that (\ref{eq:4.77}) is true for all $n\le p-1$. We show that (\ref{eq:4.77}) also hold for $n=p$ using the recursive relation (\ref{eq:4.72}). Substituting the expression (\ref{eq:4.77}) into the $h(m;s,r-1,l,t)$ and $h(m;s,r,l,t-1)$ in (\ref{eq:4.72}) and after simplification, we get
    \begin{equation}\label{eq:4.78}
        \begin{split}
            &h\left(a^sb^{m-s-r}c^rd^{m-s-l}e^{p}f^{m-r-t}g^lh^{m-l-t}k^t\cdot\text{det}_q^{-m}\right)\\
            =&\frac{(q^2;q^2)_{r+t-p}(q^2;q^2)_{l+t-p}(q^2;q^2)_{r+s-p}(q^2;q^2)_{l+s-p}}{(q^2;q^2)_m^2}\\
            &(-1)^{r+l+p}q^{(2m+1)(l+r)+(p-2)m-2l^2-2r^2-rl+4st-3ps-3pt}\\
            &\left(\sum_{k=0}^{p-1}(-1)^{k}q^{(p-1-k)(p-3k-2)+2k(s+t)}\right.\\
            &\cdot{p-1\choose k}_{q^{2}}(q^{2l};q^{-2})_k(q^{2r};q^{-2})_k(q^{2s};q^{-2})_{p-1-k}(q^{2t};q^{-2})_{p-1-k}\\
            &\left(q^{2p-2k-2}(1-q^{2(r+s-p+1)})(1-q^{2(t-p-k-1)})\right.\\
            &\left.\left.-q^{2s}(1-q^{2(l+t-p+1)})(1-q^{2r-2k})\right)\right)\\
        \end{split}
    \end{equation}
    Direct computation shows that
    \begin{equation}\label{eq:4.79}
        \begin{split}
            &q^{2p-2k-2}(1-q^{2(r+s-p+1)})(1-q^{2(t-p-k-1)})\\
            &-q^{2s}(1-q^{2(l+t-p+1)})(1-q^{2r-2k})\\
            =&q^{2p-2k-2}(1-q^{2(s-p-k-1)})(1-q^{2(t-p-k-1)})\\
            &-q^{2t+2s-2(p-1-k)}(1-q^{2(l-k)})(1-q^{2r-2k})
        \end{split}
    \end{equation}
    Substituting (\ref{eq:4.79}) into (\ref{eq:4.78}), we get
    \begin{equation}
        \begin{split}
            &h\left(a^sb^{m-s-r}c^rd^{m-s-l}e^{p}f^{m-r-t}g^lh^{m-l-t}k^t\cdot\text{det}_q^{-m}\right)\\
            =&\frac{(q^2;q^2)_{r+t-p}(q^2;q^2)_{l+t-p}(q^2;q^2)_{r+s-p}(q^2;q^2)_{l+s-p}}{(q^2;q^2)_m^2}\\
            &(-1)^{r+l+p}q^{(2m+1)(l+r)+(p-2)m-2l^2-2r^2-rl+4st-3ps-3pt}\\
            &\left(\sum_{k=0}^{p}(-1)^{k}q^{(p-k)(p-3k-1)+2k(s+t)}\left(q^{2k}{p-1\choose k}_{q^{2}}+{p-1\choose k-1}_{q^{2}}\right)\right.\\
            &\left.\cdot(q^{2l};q^{-2})_k(q^{2r};q^{-2})_k(q^{2s};q^{-2})_{p-k}(q^{2t};q^{-2})_{p-k}\right)\\
            =&\frac{(q^2;q^2)_{r+t-p}(q^2;q^2)_{l+t-p}(q^2;q^2)_{r+s-p}(q^2;q^2)_{l+s-p}}{(q^2;q^2)_m^2}\\
            &(-1)^{r+l+p}q^{(2m+1)(l+r)+(p-2)m-2l^2-2r^2-rl+4st-3ps-3pt}\\
            &\left(\sum_{k=0}^{p}(-1)^{k}q^{(p-k)(p-3k-1)+2k(s+t)}\right.\\
            &\left.\cdot{p\choose k}_{q^{2}}(q^{2l};q^{-2})_k(q^{2r};q^{-2})_k(q^{2s};q^{-2})_{p-k}(q^{2t};q^{-2})_{p-k}\right)\\
        \end{split}
    \end{equation}
    The expression is consistent with (\ref{eq:4.77}). Hence, the formula (\ref{eq:4.73}) is true for all possible $n$ values.
\end{proof}

%% file: Chapters/Chapter_4/Generalization.tex
\subsection{Generalization to $\mathcal{O}(U_q(n))$}
Consider a copy of $\mathcal{O}(U_q(3))$, denoted as $\mathcal{U}_i$, located on the reversed diagonal of the generator matrix of $\mathcal{O}(U_q(n))$. In other words, we label the following generators of $\mathcal{O}(U_q(n))$ using the same labels used in $\mathcal{O}(U_q(3))$.
\begin{equation}
    \begin{bmatrix}
        x_{i,n-1-i}&x_{i,n-i}&x_{i,n+1-i}\\
        x_{i+1,n-1-i}&x_{i+1,n-i}&x_{i+1,n+1-i}\\
        x_{i+2,n-1-i}&x_{i+2,n-i}&x_{i+2,n+1-i}\\
    \end{bmatrix}=\begin{bmatrix}
        a&b&c\\
        d&e&f\\
        g&h&k
    \end{bmatrix}.
\end{equation}
We consider the monomials of order $m$ restricted on $\mathcal{U}_i$ in the following form 
\begin{equation}\label{eq:4.89}
    p_i^1\cdot\left(a^sb^{m-s-r}c^rd^{m-s-l}e^{s+r+l+t-m}f^{m-r-t}g^lh^{m-l-t}k^t\right)\cdot p_i^2\cdot\text{det}_q^{-m},
\end{equation}
where
\begin{equation}
        p_i^1=\prod_{j=1}^{i-1}x_{i,n+1-j} \quad\text{ and }\quad p_i^2=\prod_{j=i+3}^{n}x_{i,n+1-j}
\end{equation}
commutes with all generators in $\mathcal{U}_i$. Here, we take the convention that $p_1^1=1$ and $p_{n-2}^2=1$. We will denote the Haar state values of (\ref{eq:4.89}) as
\begin{equation}
    h_n(i;m;s,r,l,t).
\end{equation}
Notice that $h_n(i;m;0,m,m,0)=h_n(\prod_{i=1}^nx^m_{i,n+1-i}\cdot\text{det}_q^{-m})$ for all $1\le i\le n-2$.

Recall that in Subsection 4.2 and 4.3 we apply the action of $\mathcal{U}_q(gl_3)$ on $\mathcal{O}(U_q(3))$ to derive the linear relations between the Haar state values of different monomials in the form of (\ref{eq:4.4}). On $\mathcal{O}(U_q(n))$, the same strategy still works in the same way. Hence, the ratio between $h_n(i;m;s,r,l,t)$ and $h_n(i;m;0,m,m,0)$ can be computed using exactly the same computations used in Subsection 4.2 and 4.3. Therefore, we have the following result.
\begin{corollary}
    On $\mathcal{O}(U_q(n))$, we have
    \begin{equation}
        \frac{h_n(i;m;s,r,l,t)}{h_n(i;m;0,m,m,0)}=\frac{h(m;s,r,l,t)}{h(m;0,m,m,0)}
    \end{equation}
    and the ration $h(m;s,r,l,t)/h(m;0,m,m,0)$ can be found by (\ref{eq:4.73}).
\end{corollary}

%% file: Chapters/Chapter_5.tex
\section{A method to pick orthonormal bases on $\mathcal{O}({U_q(3)})$}
According to Theorem 2.5, to find an orthonormal basis on $\mathcal{O}({U_q(3)})$, it suffices to find an orthonormal basis (with respect to the inner product $\langle\cdot,\cdot\rangle_L$ in (\ref{eq:2.39})) of the irreducible co-representation $V^R(\lambda)$ (See (\ref{eq:irrep})) for each dominant integral weight $\lambda=\lambda_1\epsilon_1+\lambda_2\epsilon_2+\lambda_3\epsilon_3\in F_3$. Theorem 2.1 provides a method to pick basis vectors on $V^R(\lambda)$. In this section, basis vectors on $V^R(\lambda)$ are chosen in the form of 
\begin{equation}\label{eq:5.1}
    \text{det}_q^{-\lambda_3}\xi_{\textbf{T}}
\end{equation}
indexed by $\textbf{T}\in SSYT_3(\lambda-\lambda_3(\epsilon_1+\epsilon_2+\epsilon_3))$. Here, we take the convention that $\text{det}_q^{-\lambda_3}=D_q^{-\lambda_3}$ if $\lambda_3\le 0$. In (\ref{eq:5.1}), each \textbf{T} has at most tow rows. By (\ref{eq:2.29.a}) and (\ref{eq:2.11.a}), we have the following correspondence between columns of \textbf{T} and quantum minors in corresponding $\xi_\text{\textbf{T}}$.

\begin{center}

    \begin{tabular}{|c|c|c|c|c|c|c|}
    \hline
    \multirow{3}{*}{ in \textbf{T} }
    & \multirow{3}{*}{ \begin{ytableau} 1\cr 2\cr\end{ytableau} }
    & \multirow{3}{*}{ \begin{ytableau}1\cr 3\cr\end{ytableau} }
    & \multirow{3}{*}{ \begin{ytableau}2\cr 3\cr\end{ytableau} }
    & \multirow{3}{*}{ \begin{ytableau}1\cr\end{ytableau} }
    & \multirow{3}{*}{ \begin{ytableau}2\cr\end{ytableau} }
    & \multirow{3}{*}{ \begin{ytableau}3\cr\end{ytableau} } \\ 
    &&&&&& \\ 
    &&&&&& \\ \hline
    \multirow{2}{*}{ in $\xi_\text{\textbf{T}}$ }
    & \multirow{2}{*}{ $k^*\cdot D_q$ }
    & \multirow{2}{*}{ $-q\cdot h^*\cdot D_q$ }
    & \multirow{2}{*}{ $q^2\cdot g^*\cdot D_q$ }
    & \multirow{2}{*}{ $a$ }
    & \multirow{2}{*}{ $b$ }
    & \multirow{2}{*}{ $c$ }    \\ 
    &&&&&& \\ \hline
    
    \end{tabular}

\end{center}


\noindent There are two types of semi-standard Young tableaux in $SSYT_3((\lambda_1,\lambda_2,0))$. We list the correspondence between \textbf{T} and $\xi_\text{\textbf{T}}$ as follow.
\begin{enumerate}
    \item[1)] $\xi_\mathbf{T}=(k^*)^{d_1}(-qh^*)^{d_2}a^{c_1}b^{c_2}c^{c_3}D_q^{d_1+d_2}$ corresponding to
    \begin{align*}
    &\overbrace{\hspace{4em}}^{d_1}\hspace{0.5em}\overbrace{\hspace{4em}}^{d_2}\hspace{0.5em}\overbrace{\hspace{4em}}^{c_1}\hspace{0.5em}\overbrace{\hspace{4em}}^{c_2}\hspace{0.5em}\overbrace{\hspace{4em}}^{c_3}\\[-1.5ex]
    \mathbf{T}=&\begin{ytableau}\tag{5.1.a}
        1&\none[\cdots]&1&1&\none[\cdots]&1&1&\none[\cdots]&1&2&\none[\cdots]&2&3&\none[\cdots]&3 \cr
        2&\none[\cdots]&2&3&\none[\cdots]&3  \cr
    \end{ytableau}
    \end{align*}
    \item[2)] $\xi_\mathbf{T}=(k^*)^{d_1}(-qh^*)^{d_2}(q^2g^*)^{d_3}b^{c_2}c^{c_3}D_q^{d_1+d_2+d_3}$ corresponding to
    \begin{align*}
    &\overbrace{\hspace{4em}}^{d_1}\hspace{0.5em}\overbrace{\hspace{4em}}^{d_2}\hspace{0.5em}\overbrace{\hspace{4em}}^{d_3}\hspace{0.5em}\overbrace{\hspace{4em}}^{c_2}\hspace{0.5em}\overbrace{\hspace{4em}}^{c_3}\\[-1.5ex]
    \mathbf{T}=&\begin{ytableau}\tag{5.1.b}
        1&\none[\cdots]&1&1&\none[\cdots]&1&2&\none[\cdots]&2&2&\none[\cdots]&2&3&\none[\cdots]&3 \cr
        2&\none[\cdots]&2&3&\none[\cdots]&3&3&\none[\cdots]&3  \cr
    \end{ytableau}
    \end{align*}
\end{enumerate}
To find an orthonormal basis of $V^R(\lambda)$, we need to apply the Gram-Schmidt process. To perform the Gram-Schmidt process on $V^R(\lambda)$, it suffices to know the Gram matrices consisting of the inner products between $\text{det}_q^{-\lambda_3}\xi_{\textbf{T}}$'s in $V^R(\lambda)$. By (\ref{eq:2.44}), it suffices to compute the inner products between $\text{det}_q^{-\lambda_3}\xi_{\textbf{T}}\in V^R_\mu(\lambda)$'s. It is easy to see that $V^R_\mu(\lambda)$ is one dimensional if $\lambda_1=\lambda_2$ or $\lambda_2=\lambda_3$. In the following, we assume that $\lambda_1>\lambda_2>\lambda_3$. Notice that the $\text{det}_q^{-\lambda_3}$ part of $\text{det}_q^{-\lambda_3}\xi_{\textbf{T}}$ does not affect the inner products between basis vectors of $V_\mu^R(\lambda)$ since by (\ref{eq:2.11.b}) we have
\begin{equation}\label{eq:5.2}
    h\left(\left(\text{det}_q^{-\lambda_3}\xi_{\textbf{T}}\right)^*\left(\text{det}_q^{-\lambda_3}\xi_{\textbf{T'}}\right)\right)=h\left(\xi_{\textbf{T}}^*\text{det}_q^{\lambda_3}\text{det}_q^{-\lambda_3}\xi_{\textbf{T'}}\right)=h\left(\xi_{\textbf{T}}^*\xi_{\textbf{T'}}\right).
\end{equation}
Therefore, we may further assume that $\lambda=(\lambda_1,\lambda_2,0)$. There is a natural correspondence between the basis vectors of $V^R((\lambda_1,\lambda_2,\lambda_3))$ and the basis vectors of $V^R(\lambda_1-\lambda_3,\lambda_2-\lambda_3,0)$ such that the inner products between basis vectors of $V^R((\lambda_1,\lambda_2,\lambda_3))$ equals to that of their counterpart basis vectors of $V^R(\lambda_1-\lambda_3,\lambda_2-\lambda_3,0)$ in the sense of (\ref{eq:5.2}).

Given $\textbf{T}\in SSYT_3(\lambda_1,\lambda_2,0)$ of content $\mu\in F_3$, there are two possible operations, denoted as $O_1$ and $O_2$, to obtain other tableaux of the same content
\begin{enumerate}
    \item[1.] $O_1$: switching a box with label $3$ in the first row with a box with label $2$ in the second row, and
    \item[2.] $O_2$: switching a box with label $2$ in the first row with a box with label $3$ in the second row.
\end{enumerate}
For corresponding basis vectors, if $\xi_\textbf{T}=(k^*)^{d_1}(-qh^*)^{d_2}a^{c_1}b^{c_2}c^{c_3}D_q^{d_1+d_2}$, then \\$\xi_{O_1(\textbf{T})}=(k^*)^{d_1-1}(-qh^*)^{d_2+1}a^{c_1}b^{c_2+1}c^{c_3-1}D_q^{d_1+d_2}$ and\\ $\xi_{O_2(\textbf{T})}=(k^*)^{d_1+1}(-qh^*)^{d_2-1}a^{c_1}b^{c_2-1}c^{c_3+1}D_q^{d_1+d_2}$.

Notice that in $V_\mu^R(\lambda)$ basis vectors are either corresponding to \textbf{T} in the form of (5.1.a) or (5.1.b) but not both. Hence, we will consider the two types of basis vectors separately. The method to compute the inner products between basis vectors in (5.1.a) will be discussed in detail. The result is expressed as a finite double summation of basic hypergeometric terms. For the inner products between basis vectors in (5.1.b), the result is computed using exactly the same method and the we will give the result without detail.

Given $V^R_\mu(\lambda)$, we assume that the basis vectors are in the form of (\ref{eq:5.1}.a). We denote by $v_0$ the vector basis corresponding to the $\textbf{T}_0$ of content $\mu$ with the most number of boxes with label $2$ in the first row. In other words, $v_0=\xi_{\textbf{T}_0}$ contents no $k^*$ or no $c$, depending on $\lambda$ and $\mu$. Then, basis vectors of $V^R_\mu(\lambda)$ can be written as
\begin{equation}\label{eq:5.2+}
    v_i=O_2^i(v_0)=\xi_{O_2^i(\textbf{T}_0)}.
\end{equation}
If $d_0=\#$ of $h^*$ in $v_0$ and $c_0=\#$ of $b$ in $v_0$, then $0\le i\le \text{min}(c_0,d_0)$. If $v_i$ takes the form 
\begin{equation}\label{eq:5.4++}
    v_i=(k^*)^{d_1}(-qh^*)^{d_2}a^{c_1}b^{c_2}c^{c_3}D_q^{d_1+d_2},
\end{equation}
then we have
\begin{equation}\label{eq:5.5+}
    v_{i+k}=(k^*)^{d_1+k}(-qh^*)^{d_2-k}a^{c_1}b^{c_2-k}c^{c_3+k}D_q^{d_1+d_2}.
\end{equation}
We remark that decomposing $\xi_\mathbf{T}\xi_\mathbf{T'}^*$ into a linear combination of the psudo-basis (\ref{eq:4.4}) is easier than that of $\xi_\mathbf{T}^*\xi_\mathbf{T'}$. Hence, we will compute $\langle\xi_\mathbf{T},\xi_\mathbf{T'}\rangle_R=h(\xi_\mathbf{T}\xi_\mathbf{T'}^*)$ directly. $\langle\xi_\mathbf{T},\xi_\mathbf{T'}\rangle_L=h(\xi_\mathbf{T}^*\xi_\mathbf{T'})$ can be found using the relation
\begin{equation}\label{eq:5.4+}
    h(\xi_\mathbf{T}^*\xi_\mathbf{T'})=h(\rho(\xi_\mathbf{T'})\xi_\mathbf{T}^*),
\end{equation}
where
\begin{align*}
    &\rho(\xi_\mathbf{T'})=q^{4d_1+2d_2+4c_1+2c_2}\xi_\mathbf{T'}\text{ if }\xi_\mathbf{T'}=(k^*)^{d_1}(-qh^*)^{d_2}a^{c_1}b^{c_2}c^{c_3}D_q^{d_1+d_2},\\
    &\rho(\xi_\mathbf{T'})=q^{4d_1+2d_2+2c_2}\xi_\mathbf{T'}\text{ if }\xi_\mathbf{T'}=(k^*)^{d_1}(-qh^*)^{d_2}(q^2g^*)^{d_3}b^{c_2}c^{c_3}D_q^{d_1+d_2+d_3}.
\end{align*}

\input{Chapters/Chapter_5/Direct_computation}

\input{Chapters/Chapter_5/Algorithm_with_symmetry}

\input{Chapters/Chapter_5/Inner_products_R}

\input{Chapters/Chapter_5/Inner_product_L}

%% file: Chapters/Chapter_5/Direct_computation.tex
\subsection{Direct computation}
(\ref{eq:4.73}) allows us to write inner products of $v_i$'s in (\ref{eq:5.2+}) as a finite summation of basic hypergeometric terms (Appendix \ref{append_1}). However, the expression obtained from direct computation is cumbersome. In this subsection we first give an example to illustrate the steps to compute the inner product directly and then compare our approach with Noumi, Yamada, and Mimachi's method. Before actual computations, we prove a lemma that will be used in the decomposition of $v_iv_j^*$. 
\begin{lemma}
    Let $x_{i,k}, x_{j,l}, x_{j,k}, x_{i,l}\in\mathcal{O}(GL_q(n))$ with $i<j$ and $k<l$.\\
    For $n\in\mathbb{N}$, we have
    \begin{equation}\tag{5.5}\label{eq:5.4.a}\stepcounter{equation}
            \begin{split}
                &(x_{i,k}x_{j,l}-q\cdot x_{j,k}x_{i,l})^n\\
                =&\sum_{p=0}^n(-q)^{(1-2p)(n-p)}{n\choose p}_{q^2}x_{i,k}^p(x_{j,k}x_{i,l})^{n-p}x_{j,l}^p.
            \end{split}
        \end{equation}
\end{lemma}
\begin{proof}
    Without loss of generality, we may take $x_{i,k}=a$, $x_{j,l}=e$, $x_{i,l}=b$, and $x_{j,k}=d$ on $\mathcal{O}(GL_q(3))$ so that the computation steps are easy to read.\\
    \\
    Notice that $ae$ commutes with $bd$. Hence, in the expansion of $(ae-q\cdot bd)^n$ terms are in the form of $a^p(bd)^{n-p}e^p$. Assume that 
    \begin{equation}
        \begin{split}
            (ae-q\cdot bd)^n=\sum_{p=0}^nd(n,i)\cdot a^p(bd)^{n-p}e^p.
        \end{split}
    \end{equation}
    It is easy to see that $d(1,0)=1$ and $d(1,1)=-q$. We will derive a recursive relation between $d(n+1,p)$ and $d(n,p)$. We have
    \begin{equation}
    \begin{split}
        &(ae-q\cdot bd)^{n+1}=\left(\sum_{p=0}^{n}d(n,p)a^p(bd)^{n-p}e^p\right)(ae-q\cdot bd)\\
        =&d(n,n)a^{n+1}e^{n+1}-qd(n,0)\cdot (bd)^{n+1}\\
        &+\sum_{p=1}^{n}(q^{-2(n+1-p)}d(n,p-1)-qd(n,p))a^{p}(bd)^{n+1-p}e^p.
    \end{split}
    \end{equation}
    Hence, we find that
    \begin{equation}
    d(n+1,p)=q^{-2(n+1-p)}d(n,p-1)-qd(n,p)
    \end{equation}
    with $d(n,n)=1$ and $d(n,0)=(-q)^n$. Comparing with the recursive relation between $q$-combinatorial number, we get
    \begin{equation}
        d(n,p)=(-q)^{(1-2p)(n-p)}{n\choose p}_{q^2}
    \end{equation}
    and (\ref{eq:5.4.a}) is proved.\\
\end{proof}
\subsubsection{The square length of $\xi_\mathbf{T}=(k^*)^{d_1}a^{c_1}D_q^{d_1}$ with respect to $\langle\cdot,\cdot\rangle_R$}
As an example, we directly compute square length of $\xi_\mathbf{T}=(k^*)^{d_1}a^{c_1}D_q^{d_1}$ with respect to $\langle\cdot,\cdot\rangle_R$. Notice that $k^*a=ak^*$. We have
\begin{equation}
    \begin{split}
        h(\xi_\textbf{T}\xi_\textbf{T}^*)=h((k^*)^{d_1}a^{c_1}(a^*)^{c_1}k^{d_1})=h(a^{c_1}(k^*)^{d_1}(a^*)^{c_1}k^{d_1}).
    \end{split}
\end{equation}
Expanding $(k^*)^{d_1}$ and $(a^*)^{c_1}$ using (\ref{eq:5.4.a}), we get 
\begin{equation}\label{eq:5.19}
    \begin{split}
        &h(a^{c_1}(k^*)^{d_1}(a^*)^{c_1}k^{d_1})\\
        =&\sum_{i=0}^{d_1}(-q)^{(1-2i)(d_1-i)}{d_1\choose i}_{q^2}\sum_{j=0}^{c_1}(-q)^{(1-2j)(c_1-j)}{c_1\choose j}_{q^2}\\
        &h(a^{c_1+i}b^{d_1-i}d^{d_1-i}e^{i+j}f^{c_1-j}h^{c_1-j}k^{d_1+j}\cdot\text{det}_q^{-d_1-c_1})\\
        =&\sum_{i=0}^{d_1}(-q)^{(1-2i)(d_1-i)}{d_1\choose i}_{q^2}\sum_{j=0}^{c_1}(-q)^{(1-2j)(c_1-j)}{c_1\choose j}_{q^2}\\
        &(-1)^{i+j}\frac{q^{(2d_1+2d_2-1)(d_1+c_1+i+j)-(d_1+c_1)-2(c_1+i)^2-2(d_1+j)^2}}{{d_1+c_1\choose d_1+j}_{q^2}{d_1+c_1\choose c_1+i}_{q^2}}\\
        &\cdot h\left(c^{d_1+c_1}e^{d_1+c_1}g^{d_1+c_1}\cdot\text{det}_q^{-d_1-c_1}\right)\\
        =&(-1)^{d_1+c_1}q^{-3d_1-3c_1}h\left(c^{d_1+c_1}e^{d_1+c_1}g^{d_1+c_1}\cdot\text{det}_q^{-d_1-c_1}\right)\\
        &\left(\sum_{i=0}^{d_1}\frac{q^{2(c_1+1)i}{d_1\choose i}_{q^2}}{{d_1+c_1\choose i}_{q^2}}\right)\left(\sum_{j=0}^{c_1}\frac{q^{2(d_1+1)j}{c_1\choose j}_{q^2}}{{d_1+c_1\choose j}_{q^2}}\right).
    \end{split}
\end{equation}
We denote
\begin{equation}\label{eq:5.20}
    S(d_1)=\sum_{i=0}^{d_1}\frac{q^{2(c_1+1)i}{d_1\choose i}_{q^2}}{{d_1+c_1\choose i}_{q^2}}.
\end{equation}
The $q$-Zeilberger's Algorithm~\cite{paule1997mathematica,riese1997contributions} is the standard method to find the closed form of a finite summation of basic hypergeometric terms. Applying the $q$-Zeilberger's algorithm to (\ref{eq:5.20}), we find that 
\begin{equation}
    S(d_1)=q^{2d_1}+\frac{(1-q^{2d_1})S(d_1-1)}{1-q^{2d_1+2c_1}}
\end{equation}
with $S(0)=1$. Hence, we get
\begin{equation}\label{eq:5.22}
    S(d_1)=\frac{1-q^{2d_1+2c_1+2}}{1-q^{2c_1+2}}.
\end{equation}
Substituting (\ref{eq:5.22}) into (\ref{eq:5.19}) and applying (\ref{eq:4.73}), we get
\begin{equation}\label{eq:5.23}
    h((k^*)^{d_1}a^{c_1}(a^*)^{c_1}k^{d_1})=\frac{(q^2-1)^2(q^4-1)}{(q^{2c_1+2}-1)(q^{2d_1+2}-1)(q^{2(d_1+c_1)+4}-1)}.
\end{equation}

\subsubsection{Comparing with NYM's result}
Here we compare the result of Noumi, Yamada, and Mimachi with our result on the square length of the vector $\xi_\textbf{T}=(h^*)^{d_2}(k^*)^{d_1}$ with respect to $\langle\cdot,\cdot\rangle_R$. Noumi, Yamada, and Mimachi showed that
\begin{proposition}[NYM, Proposition 4.5]
    \begin{equation}\label{eq:5.24}
        h((g^*)^{d_3}(h^*)^{d_2}(k^*)^{d_1}k^{d_1}h^{d_2}g^{d_3})=\frac{(q^{2};q^{2})_{d_1}(q^{2};q^{2})_{d_2}(q^{2};q^{2})_{d_3}(q^{2};q^{2})_{2}}{(q^{2};q^{2})_{d_1+d_2+d_3+2}}
    \end{equation}
    where $(a;q^2)_n=\prod_{i=0}^{n-1}(1-aq^{2i})$.
\end{proposition}
\noindent On the other hand, noticing that $k^*a=ak^*$ we have
\begin{equation}\label{eq:5.25}
    \begin{split}
        &h((h^*)^{d_2}(k^*)^{d_1}k^{d_1}h^{d_2})=q^{2d_12d_2}h((k^*)^{d_1}(h^*)^{d_2}h^{d_2}k^{d_1})\\
        =&(-q)^{-d_2}\sum_{j=0}^{d_2}(-q)^{(1-2j)(d_2-j)}{d_2\choose j}_{q^2}h((k^*)^{d_1}a^j(cd)^{d_2-j}f^jh^{d_2}k^{d_1})\\
        =&(-q)^{-d_2}\sum_{j=0}^{d_2}(-q)^{(1-2j)(d_2-j)}{d_2\choose j}_{q^2}h(a^j(k^*)^{d_1}(cd)^{d_2-j}f^jh^{d_2}k^{d_1})\\
        =&(-q)^{-d_2}\sum_{i=0}^{d_1}\sum_{j=0}^{d_2}(-q)^{(1-2i)(d_1-i)+(1-2j)(d_2-j)}{d_1\choose i}_{q^2}{d_2\choose j}_{q^2}\\
        &q^{-i(d_2-j)}h(a^{i+j}b^{d_1-i}c^{d_2-j}d^{d_1-i+d_2-j}e^if^jh^{d_2}k^{d_1})\\
        =&q^{2d_1d_2}\frac{(q^2;q^2)_{d_1}(q^2;q^2)_{d_2}}{(q^2;q^2)_{d_1+d_2}}\frac{q^{2(d_1+d_2)}(q^2-1)^2(q^4-1)}{(q^{2(d_1+d_2)+2}-1)^2(q^{2(d_1+d_2)+4}-1)}\\
        &\sum_{i=0}^{d_1}\sum_{j=0}^{d_2}q^{2d_1j-2ij-2i-2j}\cdot \frac{{d_1\choose i}_{q^2}{d_2\choose j}_{q^2}}{{d_1+d_2\choose i+j}_{q^2}}
    \end{split}
\end{equation}
Comparing (\ref{eq:5.24}) and (\ref{eq:5.25}) we find that
\begin{proposition}
    \begin{equation}\label{eq:5.26}
        \begin{split}
        \sum_{i=0}^{d_1}\sum_{j=0}^{d_2}q^{2d_1j-2ij-2i-2j}\frac{{d_1\choose i}_{q^2}{d_2\choose j}_{q^2}}{{d_1+d_2\choose i+j}_{q^2}}=\frac{1-q^{2(d_1+d_2)+2}}{q^{2d_1+2d_2}(1-q^2)}
        \end{split}
    \end{equation}
\end{proposition}
To prove (\ref{eq:5.24}), Noumi, Yamada, and Mimachi derived a recursive relation on the index $d_1$, $d_2$, and $d_3$ of $h((g^*)^{d_3}(h^*)^{d_2}(k^*)^{d_1}k^{d_1}h^{d_2}g^{d_3})$ using the invariant property (\ref{eq:2.38}) of the Haar state under the action of $\mathcal{U}_q(gl_3)$ on $\mathcal{O}(GL_q(3))$. Although we can use the $q$MultiSum~\cite{riese2003qmultisum} program to obtain a recursive relation on the index $d_1$ and $d_2$ of $h((h^*)^{d_2}(k^*)^{d_1}k^{d_1}h^{d_2})$ from our direct computation (\ref{eq:5.25}) and then prove (\ref{eq:5.26}), Noumi, Yamada, and Mimachi's approach provides a quantum group explanation of the equality (\ref{eq:5.24}) and the recursive relation that (\ref{eq:5.24}) obeys. In the following computation, we will develop an algorithm applying the invariant property of the Haar state under the action of $\mathcal{U}_q(gl_3)$ on $\mathcal{O}(GL_q(3))$ to compute the inner products between basis vectors of $V^R_\mu(\lambda)$.

%% file: Chapters/Chapter_5/Algorithm_with_symmetry.tex
\subsection{An algorithm applying property (\ref{eq:2.38}) }
Given $V_\mu^R(\lambda)$, to find the inner products between all basis vectors, it suffices to compute $\langle v_i,v_j\rangle_L=h(v_i^*v_j)$ with $i\le j$. The case $i>j$ can be found using the relation $\langle v_i,v_j\rangle_L=\overline{\langle v_j,v_i\rangle_L}$. To keep consistent with previous computations, we will first compute $\langle v_j,v_i\rangle_R$ and then obtain $\langle v_i,v_j\rangle_L$ using relation (\ref{eq:5.4+}). The $q$ factors in $v_i$'s will be omitted during the computation and put back at the last step of our computation.
\subsubsection{The square length of $(k^*)^{d_1}(h^*)^{d_2}a^{c_1}D_q^{d_1+d_2}$ with respect to $\langle \cdot,\cdot\rangle_R$}
Consider the left action of $f_2\in\mathcal{U}_q(gl_3)$ on $(k^*)^{d_1+1}(h^*)^{d_2}a^{c_1}(a^*)^{c_1}h^{d_2+1}k^{d_1}$. We have
\begin{equation}\label{eq:5.29}
    \begin{split}
        &f_2\cdot (k^*)^{d_1+1}(h^*)^{d_2}a^{c_1}(a^*)^{c_1}h^{d_2+1}k^{d_1}\\
        =&-q^{-d_1+1/2}\frac{1-q^{2(d_1+1)}}{1-q^2}(k^*)^{d_1}(h^*)^{d_2+1}a^{c_1}(a^*)^{c_1}h^{d_2+1}k^{d_1}\\
        &+q^{d_1+1/2-2d_2}\frac{1-q^{2(d_2+1)}}{1-q^2}(k^*)^{d_1+1}(h^*)^{d_2}a^{c_1}(a^*)^{c_1}h^{d_2}k^{d_1+1}.
    \end{split}
\end{equation}
Evaluating the Haar state on both sides of (\ref{eq:5.29}) and applying invariant property (\ref{eq:2.38}), we get
\begin{equation}\label{eq:5.30}
    \begin{split}
        &h\left((k^*)^{d_1}(h^*)^{d_2+1}a^{c_1}(a^*)^{c_1}h^{d_2+1}k^{d_1}\right)\\
        =&q^{2d_1-2d_2}\frac{1-q^{2(d_2+1)}}{1-q^{2(d_1+1)}}h\left((k^*)^{d_1+1}(h^*)^{d_2}a^{c_1}(a^*)^{c_1}h^{d_2}k^{d_1+1}\right).
    \end{split}
\end{equation}
Applying (\ref{eq:5.30}) repeatedly, we get
\begin{equation}\label{eq:5.31}
    \begin{split}
         &h\left((k^*)^{d_1}(h^*)^{d_2}a^{c_1}(a^*)^{c_1}h^{d_2}k^{d_1}\right)\\
         =&\frac{q^{2d_1d_2}(1-q^2)^2(1-q^4)(q^2;q^2)_{d_1}(q^2;q^2)_{d_2}}{(1-q^{2c_1+2})(q^2;q^2)_{d_1+d_2+1}(1-q^{2(d_1+d_2+c_1)+4})}.
    \end{split}
\end{equation}
\subsubsection{The square length of $(k^*)^{d_1}(h^*)^{d_2}a^{c_1}b^{c_2}D_q^{d_1+d_2}$ with respect to $\langle \cdot,\cdot\rangle_R$}
We will design an induction scheme on the value of $c_2$ to compute the square length of $(k^*)^{d_1}(h^*)^{d_2}a^{c_1}b^{c_2}$. Consider the left action of $e_1\in \mathcal{U}_q(gl_3)$ on\\ $(k^*)^{d_1}(h^*)^{d_2}a^{c_1}b^{c_2+1}(b^*)^{c_2}(a^*)^{c_1+1}h^{d_2}k^{d_1}$. We have
\begin{equation}\label{eq:5.32}
    \begin{split}
        &e_1\cdot (k^*)^{d_1}(h^*)^{d_2}a^{c_1}b^{c_2+1}(b^*)^{c_2}(a^*)^{c_1+1}h^{d_2}k^{d_1}\\
        =&q^{d_2+c_1+1/2}\frac{1-q^{-2(c_2+1)}}{1-q^{-2}}(k^*)^{d_1}(h^*)^{d_2}a^{c_1+1}b^{c_2}(b^*)^{c_2}(a^*)^{c_1+1}h^{d_2}k^{d_1}\\
        &-q^{d_2+c_1-3/2}\frac{1-q^{-2(c_1+1)}}{1-q^{-2}}(k^*)^{d_1}(h^*)^{d_2}a^{c_1}b^{c_2+1}(b^*)^{c_2+1}(a^*)^{c_1}h^{d_2}k^{d_1}\\
        &+q^{d_2-3/2}\frac{1-q^{-2d_2}}{1-q^{-2}}(k^*)^{d_1}(h^*)^{d_2}a^{c_1}b^{c_2+1}(b^*)^{c_2}(a^*)^{c_1+1}gh^{d_2-1}k^{d_1}.
    \end{split}
\end{equation}
For the extra term $(k^*)^{d_1}(h^*)^{d_2}a^{c_1}b^{c_2+1}(b^*)^{c_2}(a^*)^{c_1+1}gh^{d_2-1}k^{d_1}$, we consider the left action of $f_1\in \mathcal{U}_q(gl_3)$ on it and get
\begin{equation}\label{eq:5.33}
    \begin{split}
        &f_1\cdot (k^*)^{d_1}(h^*)^{d_2}a^{c_1+1}b^{c_2}(b^*)^{c_2}(a^*)^{c_1+1}gh^{d_2-1}k^{d_1}\\
        =&-q^{3/2-d_2}\frac{1-q^{2d_2}}{1-q^2}(k^*)^{d_1}(h^*)^{d_2-1}g^*a^{c_1+1}b^{c_2}(b^*)^{c_2}(a^*)^{c_1+1}gh^{d_2-1}k^{d_1}\\
        &+q^{d_2-c_1-1/2}\frac{1-q^{2(c_1+1)}}{1-q^2}(k^*)^{d_1}(h^*)^{d_2}a^{c_1}b^{c_2+1}(b^*)^{c_2}(a^*)^{c_1+1}gh^{d_2-1}k^{d_1}\\
        &-q^{d_2+c_1-2c_2+5/2}\frac{1-q^{2c_2}}{1-q^2}(k^*)^{d_1}(h^*)^{d_2}a^{c_1+1}b^{c_2}(b^*)^{c_2-1}(a^*)^{c_1+2}gh^{d_2-1}k^{d_1}\\
        &+q^{d_2-1/2}(k^*)^{d_1}(h^*)^{d_2}a^{c_1+1}b^{c_2}(b^*)^{c_2}(a^*)^{c_1+1}h^{d_2}k^{d_1}.
    \end{split}
\end{equation}
Evaluating the Haar state on both sides of (\ref{eq:5.32}) and (\ref{eq:5.33}) and applying invariant property (\ref{eq:2.38}), we get
\begin{equation}\label{eq:5.34}
    \begin{split}
        &q^{c_1}\frac{1-q^{-2(c_1+1)}}{1-q^{-2}}h\left((k^*)^{d_1}(h^*)^{d_2}a^{c_1}b^{c_2+1}(b^*)^{c_2+1}(a^*)^{c_1}h^{d_2}k^{d_1}\right)\\
        =&q^{c_1+2}\frac{1-q^{-2(c_2+1)}}{1-q^{-2}}h\left((k^*)^{d_1}(h^*)^{d_2}a^{c_1+1}b^{c_2}(b^*)^{c_2}(a^*)^{c_1+1}h^{d_2}k^{d_1}\right)\\
        &+\frac{1-q^{-2d_2}}{1-q^{-2}}h\left((k^*)^{d_1}(h^*)^{d_2}a^{c_1}b^{c_2+1}(b^*)^{c_2}(a^*)^{c_1+1}gh^{d_2-1}k^{d_1}\right)
    \end{split}
\end{equation}
and
\begin{equation}\label{eq:5.35}
    \begin{split}
        &q^{d_2-c_1}\frac{1-q^{2(c_1+1)}}{1-q^2}h\left((k^*)^{d_1}(h^*)^{d_2}a^{c_1}b^{c_2+1}(b^*)^{c_2}(a^*)^{c_1+1}gh^{d_2-1}k^{d_1}\right)\\
        =&q^{d_2+c_1-2c_2+3}\frac{1-q^{2c_2}}{1-q^2}h\left((k^*)^{d_1}(h^*)^{d_2}a^{c_1+1}b^{c_2}(b^*)^{c_2-1}(a^*)^{c_1+2}gh^{d_2-1}k^{d_1}\right)\\
        &+q^{2-d_2}\frac{1-q^{2d_2}}{1-q^2}h\left((k^*)^{d_1}(h^*)^{d_2-1}g^*a^{c_1+1}b^{c_2}(b^*)^{c_2}(a^*)^{c_1+1}gh^{d_2-1}k^{d_1}\right)\\
        &-q^{d_2}h\left((k^*)^{d_1}(h^*)^{d_2}a^{c_1+1}b^{c_2}(b^*)^{c_2}(a^*)^{c_1+1}h^{d_2}k^{d_1}\right).
    \end{split}
\end{equation}
Recall that $1=kk^*+hh^*+gg^*$. Hence, we have
\begin{equation}\label{eq:5.36}
    \begin{split}
        &h\left((k^*)^{d_1}(h^*)^{d_2-1}g^*a^{c_1+1}b^{c_2}(b^*)^{c_2}(a^*)^{c_1+1}gh^{d_2-1}k^{d_1}\right)\\
        =&q^{2(d_1+d_2-1)}h\left(\sigma(g)g^*(k^*)^{d_1}(h^*)^{d_2-1}a^{c_1+1}b^{c_2}(b^*)^{c_2}(a^*)^{c_1+1}h^{d_2-1}k^{d_1}\right)\\
        =&q^{2(d_1+d_2-1)}h\left((k^*)^{d_1}(h^*)^{d_2-1}a^{c_1+1}b^{c_2}(b^*)^{c_2}(a^*)^{c_1+1}h^{d_2-1}k^{d_1}\right)\\
        &-q^{2d_2}h\left((k^*)^{d_1}(h^*)^{d_2}a^{c_1+1}b^{c_2}(b^*)^{c_2}(a^*)^{c_1+1}h^{d_2}k^{d_1}\right)\\
        &-q^{2(d_1+d_2+1)}h\left((k^*)^{d_1+1}(h^*)^{d_2-1}a^{c_1+1}b^{c_2}(b^*)^{c_2}(a^*)^{c_1+1}h^{d_2-1}k^{d_1+1}\right).
    \end{split}
\end{equation}
Substituting (\ref{eq:5.36}) into (\ref{eq:5.35}), we get
\begin{align}\label{eq:5.37}\stepcounter{equation}
        &q^{-c_1}\frac{1-q^{2(c_1+1)}}{1-q^2}h\left((k^*)^{d_1}(h^*)^{d_2}a^{c_1}b^{c_2+1}(b^*)^{c_2}(a^*)^{c_1+1}gh^{d_2-1}k^{d_1}\right)\notag\\
        =&q^{c_1-2c_2+3}\frac{1-q^{2c_2}}{1-q^2}h\left((k^*)^{d_1}(h^*)^{d_2}a^{c_1+1}b^{c_2}(b^*)^{c_2-1}(a^*)^{c_1+2}gh^{d_2-1}k^{d_1}\right)\notag\\
        &+q^{2d_1}\frac{1-q^{2d_2}}{1-q^2}h\left((k^*)^{d_1}(h^*)^{d_2-1}a^{c_1+1}b^{c_2}(b^*)^{c_2}(a^*)^{c_1+1}h^{d_2-1}k^{d_1}\right)\tag{5.29}\\
        &-\frac{1-q^{2d_2+2}}{1-q^2}h\left((k^*)^{d_1}(h^*)^{d_2}a^{c_1+1}b^{c_2}(b^*)^{c_2}(a^*)^{c_1+1}h^{d_2}k^{d_1}\right)\notag\\
        &-q^{2(d_1+2)}\frac{1-q^{2d_2}}{1-q^2}h\left((k^*)^{d_1+1}(h^*)^{d_2-1}a^{c_1+1}b^{c_2}(b^*)^{c_2}(a^*)^{c_1+1}h^{d_2-1}k^{d_1+1}\right).\notag
\end{align}
Combining (\ref{eq:5.34}) and (\ref{eq:5.37}), we find a recursive scheme on the value of $c_2=0,1,2,\cdots$ as follow.
\begin{enumerate}
    \item[1.] Compute the Haar state of $h\left((k^*)^{d_1}(h^*)^{d_2}a^{c_1}(a^*)^{c_1}h^{d_2}k^{d_1}\right)$ and\\ $h\left((k^*)^{d_1}(h^*)^{d_2}a^{c_1}b(a^*)^{c_1+1}gh^{d_2-1}k^{d_1}\right)$. This is the case $c_2=0$.
    \item[2.] Assume that the general form of $h\left((k^*)^{d_1}(h^*)^{d_2}a^{c_1}b^{n}(b^*)^{n}(a^*)^{c_1}h^{d_2}k^{d_1}\right)$ and $h\left((k^*)^{d_1}(h^*)^{d_2}a^{c_1}b^{n+1}(b^*)^{n}(a^*)^{c_1+1}gh^{d_2-1}k^{d_1}\right)$ is known. This is the case $c_2=n$.
    \item[3.] Compute the general form of $h\left((k^*)^{d_1}(h^*)^{d_2}a^{c_1}b^{n+1}(b^*)^{n+1}(a^*)^{c_1}h^{d_2}k^{d_1}\right)$ using Equation (\ref{eq:5.34}).
    \item[4.] Compute the general form of $h\left((k^*)^{d_1}(h^*)^{d_2}a^{c_1}b^{n+2}(b^*)^{n+1}(a^*)^{c_1+1}gh^{d_2-1}k^{d_1}\right)$ using Equation (\ref{eq:5.37}). This complete the case $c_2=n+1$.
\end{enumerate}
Based on the inductive scheme, we achieve the following result.
\begin{proposition}\label{prop:5.4}
We have
    \begin{equation}\label{eq:5.38}
            \begin{split}
                &h\left((k^*)^{d_1}(h^*)^{d_2}a^{c_1}b^{c_2}(b^*)^{c_2}(a^*)^{c_1}h^{d_2}k^{d_1}\right)\\
                =&\frac{q^{2d_1d_2+2c_1c_2+2c_2}(1-q^2)^2(1-q^4)(q^2;q^2)_{d_2}(q^2;q^2)_{c_2}}{(q^2;q^2)_{d_1+d_2+1}(q^2;q^2)_{c_1+c_2+1}(q^{2(d_1+d_2+c_1)+4};q^2)_{c_2+1}}\\
                &\left(\sum_{i=0}^{c_2}q^{2(d_1+d_2+1)i}(q^2;q^2)_{d_1+c_2-i}(q^2;q^2)_{c_1+i}{c_2\choose i}_{q^2}\right)
                \end{split}
            \end{equation} 
         and   
       \begin{equation}\label{eq:5.39}
            \begin{split}
                &h\left((k^*)^{d_1}(h^*)^{d_2}a^{c_1}b^{c_2+1}(b^*)^{c_2}(a^*)^{c_1+1}gh^{d_2-1}k^{d_1}\right)\\
                =&-\frac{q^{2d_1d_2+2d_1+2d_2+c_1+2+2c_1c_2+4c_2}(1-q^2)^2(1-q^4)(q^2;q^2)_{c_2+1}(q^2;q^2)_{d_2}}{(q^2;q^2)_{d_1+d_2+1}(q^2;q^2)_{c_1+c_2+2}(q^{2(d_1+d_2+c_1)+4};q^2)_{c_2+2}}\\
                &\left(\sum_{i=0}^{c_2}q^{2(d_1+d_2)i}(q^2;q^2)_{d_1+c_2-i}(q^2;q^2)_{c_1+1+i}{c_2\choose i}_{q^2}\right).
            \end{split}
        \end{equation}
\end{proposition}
\begin{proof}
    (\ref{eq:5.38}) and (\ref{eq:5.39}) are proved by following the induction scheme consisting of (\ref{eq:5.34}) and (\ref{eq:5.37}). We start with the case $c_2=0$. The value of \\$h\left((k^*)^{d_1}(h^*)^{d_2}a^{c_1}(a^*)^{c_1}h^{d_2}k^{d_1}\right)$ can be found in (\ref{eq:5.31}) and is consistent with (\ref{eq:5.38}) when $c_2=0$. $h\left((k^*)^{d_1}(h^*)^{d_2}a^{c_1}b(a^*)^{c_1+1}gh^{d_2-1}k^{d_1}\right)$ can be computed by the relation derived from $f_1\cdot (k^*)^{d_1}(h^*)^{d_2}a^{c_1+1}(a^*)^{c_1+1}gh^{d_2-1}k^{d_1}$. The derived relation is consistent with the relation obtained by putting $c_2=0$ in (\ref{eq:5.37}) and erasing the term $h\left((k^*)^{d_1}(h^*)^{d_2}a^{c_1+1}b^{c_2}(b^*)^{c_2-1}(a^*)^{c_1+2}gh^{d_2-1}k^{d_1}\right)$. Hence, we have
    \begin{equation}\label{eq:5.40}
        \begin{split}
        &q^{-c_1}\frac{1-q^{2(c_1+1)}}{1-q^2}h\left((k^*)^{d_1}(h^*)^{d_2}a^{c_1}b(a^*)^{c_1+1}gh^{d_2-1}k^{d_1}\right)\\
        =&q^{2d_1}\frac{1-q^{2d_2}}{1-q^2}h\left((k^*)^{d_1}(h^*)^{d_2-1}a^{c_1+1}(a^*)^{c_1+1}h^{d_2-1}k^{d_1}\right)\\
        &-\frac{1-q^{2d_2+2}}{1-q^2}h\left((k^*)^{d_1}(h^*)^{d_2}a^{c_1+1}(a^*)^{c_1+1}h^{d_2}k^{d_1}\right)\\
        &-q^{2(d_1+2)}\frac{1-q^{2d_2}}{1-q^2}h\left((k^*)^{d_1+1}(h^*)^{d_2-1}a^{c_1+1}(a^*)^{c_1+1}h^{d_2-1}k^{d_1+1}\right).
        \end{split}
    \end{equation}
    Substituting (\ref{eq:5.31}) into (\ref{eq:5.40}), we get
    \begin{equation}\label{eq:5.41}
        \begin{split}
        &h\left((k^*)^{d_1}(h^*)^{d_2}a^{c_1}b(a^*)^{c_1+1}gh^{d_2-1}k^{d_1}\right)\\
        =&-\frac{q^{2d_1d_2+2d_2+2d_1+c_1+2}(1-q^2)^3(1-q^4)(q^2;q^2)_{d_2}(q^2;q^2)_{d_1}}{(1-q^{2c_1+4})(q^2;q^2)_{d_1+d_2+1}(q^{2(d_1+d_2+c_1)+4};q^2)_{2}}\\
        \end{split}
    \end{equation}
     (\ref{eq:5.40}) is consistent with (\ref{eq:5.39}) when $c_2=0$. Hence, the case $c_2=0$ is verified.

     The proof from the case $c_2=n$ to the case $c_2=n+1$ follows from direct computation. For the detail of this computation, see Appendix \ref{append_a}.
\end{proof}

\subsubsection{The inner product between $(k^*)^{d_1+k}(h^*)^{d_2-k}a^{c_1}b^{c_2-k}c^kD_q^{d_1+d_2}$ and\\ $(k^*)^{d_1}(h^*)^{d_2}a^{c_1}b^{c_2}D_q^{d_1+d_2}$ with respect to $\langle \cdot,\cdot\rangle_R$}
Consider the left action of $f_2\in\mathcal{U}_q(gl_3)$ on $(k^*)^{d_1+k}(h^*)^{d_2-k}a^{c_1}b^{c_2-k+1}c^{k-1}(b^*)^{c_2}(a^*)^{c_1}h^{d_2}k^{d_1}$. We have
\begin{align}\label{eq:5.31+}\stepcounter{equation}
        &f_2\cdot (k^*)^{d_1+k}(h^*)^{d_2-k}a^{c_1}b^{c_2-k+1}c^{k-1}(b^*)^{c_2}(a^*)^{c_1}h^{d_2}k^{d_1}\tag{5.34}\\
        =&-q^{d_1+k+3/2}\frac{q^{-2}-q^{-2(d_1+k+1)}}{1-q^{-2}}(k^*)^{d_1+k-1}(h^*)^{d_2-k+1}a^{c_1}b^{c_2-k+1}c^{k-1}(b^*)^{c_2}(a^*)^{c_1}h^{d_2}k^{d_1}\notag\\
        &+q^{d_1-d_2+k+c_2+3/2}\frac{q^{-2}-q^{-2(c_2-k+2)}}{1-q^{-2}}(k^*)^{d_1+k}(h^*)^{d_2-k}a^{c_1}b^{c_2-k}c^{k}(b^*)^{c_2}(a^*)^{c_1}h^{d_2}k^{d_1}\notag\\
        &+q^{d_1+5/2}\frac{q^{-2}-q^{-2(d_2+1)}}{1-q^{-2}}(k^*)^{d_1+k}(h^*)^{d_2-k}a^{c_1}b^{c_2-k+1}c^{k-1}(b^*)^{c_2}(a^*)^{c_1}h^{d_2-1}k^{d_1+1}.\notag
\end{align}
Evaluating the Haar state on both sides of (\ref{eq:5.31+}) and applying invariant property (\ref{eq:2.38}), we get
\begin{equation}\label{eq:5.32+}
    \begin{split}
        &h((k^*)^{d_1+k}(h^*)^{d_2-k}a^{c_1}b^{c_2-k}c^{k}(b^*)^{c_2}(a^*)^{c_1}h^{d_2}k^{d_1})\\
        =&q^{-2d_1+d_2+c_2-4k+2}\frac{1-q^{2(d_1+k)}}{1-q^{2(c_2-k+1)}}\\
        &h((k^*)^{d_1+k-1}(h^*)^{d_2-k+1}a^{c_1}b^{c_2-k+1}c^{k-1}(b^*)^{c_2}(a^*)^{c_1}h^{d_2}k^{d_1})\\
        &-q^{-d_2+c_2-3k+3}\frac{1-q^{2d_2}}{1-q^{2(c_2-k+1)}}\\
        &h((k^*)^{d_1+k}(h^*)^{d_2-k}a^{c_1}b^{c_2-k+1}c^{k-1}(b^*)^{c_2}(a^*)^{c_1}h^{d_2-1}k^{d_1+1}).
    \end{split}
\end{equation}
(\ref{eq:5.32+}) gives a recursive relation on the value of $0\le k\le \text{min}(d_2,c_2)$. The initial value of the recursive relation comes from (\ref{eq:5.38}). Using (\ref{eq:5.32+}) and (\ref{eq:5.38}), we get the following result.
\begin{proposition} We have
    \begin{equation}\label{eq:5.33+}
        \begin{split}
            &h((k^*)^{d_1+k}(h^*)^{d_2-k}a^{c_1}b^{c_2-k}c^{k}(b^*)^{c_2}(a^*)^{c_1}h^{d_2}k^{d_1})\\
            =&(-1)^k\frac{q^{2d_1d_2+2c_1c_2+(k+2)c_2+k(d_2-k+1)}(1-q^2)^2(1-q^4)(q^2;q^2)_{d_2}(q^2;q^2)_{c_2}}{(q^2;q^2)_{c_1+c_2+1}(q^2;q^2)_{d_1+d_2+1}(q^{2(d_1+d_2+c_1+i)+4};q^2)_{c_2+1}}\\
            &\left(\sum_{i=0}^{c_2-k}q^{2(d_1+d_2+1)i}(q^2;q^2)_{d_1+c_2-i}(q^2;q^2)_{c_1+i}{c_2-k\choose i}_{q^2}\right)\\
        \end{split}
    \end{equation}
\end{proposition}
\begin{proof}
    It is easy to see that (\ref{eq:5.33+}) is true for $d_2,c_2\ge 0$ in the case $k=0$. Assume that (\ref{eq:5.33+}) is true for all $d_2,c_2\ge n-1$ in the case $k=n-1$. We compute the case $k=n$ using recursive relation (\ref{eq:5.32+}). For all $d_2,c_2\ge n$, we have
    \begin{align}\label{eq:5.36+}\stepcounter{equation}
            &h((k^*)^{d_1+n}(h^*)^{d_2-n}a^{c_1}b^{c_2-n}c^{n}(b^*)^{c_2}(a^*)^{c_1}h^{d_2}k^{d_1})\notag\\
            =&(-1)^{n-1}\frac{q^{2d_1d_2+2c_1c_2+(n+2)c_2+n(d_2-n+1)-2d_1-2n}(1-q^2)^2(1-q^4)(q^2;q^2)_{d_2}(q^2;q^2)_{c_2}}{(q^2;q^2)_{c_1+c_2+1}(q^2;q^2)_{d_1+d_2+1}(q^{2(d_1+d_2+c_1+i)+4};q^2)_{c_2+1}}\notag\\
            &\left(\sum_{i=0}^{c_2-n+1}q^{2(d_1+d_2+1)i}(q^2;q^2)_{c_1+i}\frac{(q^2;q^2)_{c_2-n}}{(q^2;q^2)_{c_2-n+1-i}(q^2;q^2)_{i}}\right.\tag{5.37}\\
            &\left.\left((1-q^{2(d_1+n)})(q^2;q^2)_{d_1+c_2-i}-(q^2;q^2)_{d_1+1+c_2-i}\right)\right)\notag\\
            =&(-1)^n\frac{q^{2d_1d_2+2c_1c_2+(n+2)c_2+n(d_2-n+1)}(1-q^2)^2(1-q^4)(q^2;q^2)_{d_2}(q^2;q^2)_{c_2}}{(q^2;q^2)_{c_1+c_2+1}(q^2;q^2)_{d_1+d_2+1}(q^{2(d_1+d_2+c_1+i)+4};q^2)_{c_2+1}}\notag\\
            &\left(\sum_{i=0}^{c_2-n}q^{2(d_1+d_2+1)i}(q^2;q^2)_{d_1+c_2-i}(q^2;q^2)_{c_1+i}{c_2-n\choose i}_{q^2}\right).\notag
    \end{align}
    This is consistent with (\ref{eq:5.33+}). Hence, the case $k=n$ is verified for all $d_2,c_2\ge n$.
\end{proof}

\subsubsection{The inner product between $(k^*)^{d_1+k}(h^*)^{d_2-k}a^{c_1}b^{c_2-k}c^{c_3+k}D_q^{d_1+d_2}$ and\\ $(k^*)^{d_1}(h^*)^{d_2}a^{c_1}b^{c_2}c^{c_3}D_q^{d_1+d_2}$ with respect to $\langle \cdot,\cdot\rangle_R$}
Recall that $1=aa^*+bb^*+cc^*$. Hence, we have
\begin{equation}\label{eq:5.41+}
    \begin{split}
        &h((k^*)^{d_1+k}(h^*)^{d_2-k}a^{c_1}b^{c_2-k}c^{c_3+1+k}(c^*)^{c_3+1}(b^*)^{c_2}(a^*)^{c_1}h^{d_2}k^{d_1})\\
        =&h((k^*)^{d_1+k}(h^*)^{d_2-k}a^{c_1}b^{c_2-k}c^{c_3+k}(c^*)^{c_3}(b^*)^{c_2}(a^*)^{c_1}h^{d_2}k^{d_1})\\
        &-q^{-2(c_2+c_3)}\\
        &h((k^*)^{d_1+k}(h^*)^{d_2-k}a^{c_1+1}b^{c_2-k}c^{c_3+k}(c^*)^{c_3}(b^*)^{c_2}(a^*)^{c_1+1}h^{d_2}k^{d_1})\\
        &-q^{-(2c_3+k)}\\
        &h((k^*)^{d_1+k}(h^*)^{d_2-k}a^{c_1}b^{c_2+1-k}c^{c_3+k}(c^*)^{c_3}(b^*)^{c_2+1}(a^*)^{c_1}h^{d_2}k^{d_1}).
    \end{split}
\end{equation}
(\ref{eq:5.41+}) is a recursive relation on the value of $c_3$ with initial condition (\ref{eq:5.33+}). By induction, we get
\begin{proposition}\label{prop:5.6}
    We have
    \begin{align}\label{eq:5.35+}\stepcounter{equation}
            &h((k^*)^{d_1+k}(h^*)^{d_2-k}a^{c_1}b^{c_2-k}c^{c_3+k}(c^*)^{c_3}(b^*)^{c_2}(a^*)^{c_1}h^{d_2}k^{d_1})\notag\\
            =&(-1)^kq^{2d_1d_2+2c_1c_2+2c_1c_3+2c_2c_3+2c_2+4c_3+k(d_2+c_2-k+1)}\tag{5.39}\\
            &\frac{(1-q^2)^2(1-q^4)(q^2;q^2)_{d_2}(q^2;q^2)_{c_2}}{(q^2;q^2)_{d_1+d_2+1}(q^2;q^2)_{c_1+c_2+c_3+1}}\left(\sum_{j=0}^{c_3}\frac{(-1)^{j}q^{j^2-j}{d_1\choose j}_{q^{2}}{c_3\choose j}_{q^2}(q^2;q^2)_{j}}{(q^{2d_1+2d_2+2c_1+2c_3-2j+4};q^2)_{c_2+1}}\right.\notag\\
            &\left.\cdot\sum_{i=0}^{c_2-k}q^{(2d_1+2d_2+2c_3-2j+2)i}(q^2;q^2)_{c_1+i}(q^2;q^2)_{d_1+c_2+c_3-j-i}{c_2-k\choose i}_{q^2}\right).\notag
    \end{align}
\end{proposition}
\begin{proof}
    The proof of general form (\ref{eq:5.35+}) comes from direct computation following the recursive relation (\ref{eq:5.41+}). For detail of this computation, see Appendix(\ref{append_b}).
\end{proof}
\begin{corollary}
    \begin{align}\label{eq:5.40+}\stepcounter{equation}
            &h((k^*)^{d_1}a^{c_1}b^{c_2}c^{c_3}(c^*)^{c_3}(b^*)^{c_2}(a^*)^{c_1}k^{d_1})\tag{5.40}\\
            =&\frac{q^{2c_1c_2+2c_1c_3+2c_2c_3+2c_2+4c_3+2d_1c_3}(1-q^2)^2(1-q^4)(q^2;q^2)_{c_1}(q^2;q^2)_{c_2}(q^2;q^2)_{c_3}}{(q^2;q^2)_{c_1+c_2+1}(1-q^{2d_1+2})(q^{2d_1+2c_1+2c_2+4};q^2)_{c_3+1}}.\notag
    \end{align}
\end{corollary}
\begin{proof}
    We can prove (\ref{eq:5.40+}) by applying the $q$-Zeilberger's algorithm. First, we sum over index $i$ and then sum over index $j$.
\end{proof}
We remark that in the commutative case ($q=1$), the object in (\ref{eq:5.39}), (\ref{eq:5.40+}) and (\ref{eq:5.35+}) can be computed using the Weingarten calculus~\cite{collins2003moments} on $U(3)$. Readers may compare (\ref{eq:5.39}), (\ref{eq:5.40+}) and (\ref{eq:5.35+}) with their commutative counter part on $U(3)$.

%% file: Chapters/Chapter_5/Inner_products_R.tex
\subsection{The inner products of basis vectors on $V^R_\mu(\lambda)$}
\subsubsection{Inner products between basis vectors of form (\ref{eq:5.1}.a)}
Combining (\ref{eq:5.4+}) and (\ref{eq:5.35+}), we have the following result on the inner products between basis vectors of form (\ref{eq:5.1}.a).
\begin{theorem}
    Assume that we enumerate basis vectors of $V^R_\mu(\lambda)$ as in (\ref{eq:5.2+}) and these basis vectors take the form of (\ref{eq:5.1}.a). If $v_i,v_{i+k}\in V^R_\mu(\lambda)$, $0\le k\le\text{min}(d_2,c_2)$, takes the form (\ref{eq:5.4++}) and (\ref{eq:5.5+}), respectively, then the inner product between $v_i$ and $v_{i+k}$ with respect to the norm $\langle\cdot,\cdot\rangle_L$ is
    \begin{align}\stepcounter{equation}
            &\langle v_{i},v_{i+k}\rangle_L=h(v^*_iv_{i+k})\notag\\
            =&q^{2d_1d_2+4d_1+4d_2+2c_1c_2+2c_1c_3+2c_2c_3+4c_1+4c_2+4c_3+k(d_2+c_2-k)}\tag{5.41}\\
            &\frac{(1-q^2)^2(1-q^4)(q^2;q^2)_{d_2}(q^2;q^2)_{c_2}}{(q^2;q^2)_{d_1+d_2+1}(q^2;q^2)_{c_1+c_2+c_3+1}}\left(\sum_{j=0}^{c_3}\frac{(-1)^{j}q^{j^2-j}{d_1\choose j}_{q^{2}}{c_3\choose j}_{q^2}(q^2;q^2)_{j}}{(q^{2d_1+2d_2+2c_1+2c_3-2j+4};q^2)_{c_2+1}}\right.\notag\\
            &\left.\cdot\sum_{i=0}^{c_2-k}q^{(2d_1+2d_2+2c_3-2j+2)i}(q^2;q^2)_{c_1+i}(q^2;q^2)_{d_1+c_2+c_3-j-i}{c_2-k\choose i}_{q^2}\right).\notag
    \end{align}
\end{theorem}
\subsubsection{Inner products between basis vectors of form (\ref{eq:5.1}.b)}
If basis vectors of $V^R_\mu(\lambda)$ take the form of (\ref{eq:5.1}.b), we may enumerate these vectors in the same way as we did in (\ref{eq:5.2+}). If $v_i$ is written as
\begin{equation}\label{eq:5.44+}
    v_i=(k^*)^{d_1}(-qh^*)^{d_2}(q^2g^*)^{d_3}b^{c_2}c^{c_3}D_q^{d_1+d_2+d_3},
\end{equation}
then
\begin{equation}\label{eq:5.45+}
    v_{i+k}=(k^*)^{d_1+k}(-qh^*)^{d_2-k}(q^2g^*)^{d_3}b^{c_2-k}c^{c_3+k}D_q^{d_1+d_2+d_3}.
\end{equation}
The inner products between $v_i$ and $v_{i+k}$ can be computed using the same strategy as in Subsection 5.2. We list the recursive relations used in the computation and intermediate results as follow but omit the computation steps. \\
\\
1. Direct computation gives
\begin{equation}\label{eq:5.43}
    \begin{split}
        h((g^*)^{d_3}c^{c_3}(c^*)^{c_3}g^{d_3})=\frac{q^{4c_3}(q^2-1)^2(q^4-1)}{(1-q^{2c_3+2})(1-q^{2d_3+2})(q^{2(d_3+c_3)+4}-1)}.
    \end{split}
\end{equation}
2. Using the linear relation derived from the left action of $f_2\in\mathcal{U}_q(gl(3))$ on $(g^*)^{d_3}b^{c_2+1}c^{c_3-1}(c^*)^{c_3}(b^*)^{c_2}g^{d_3}$ and the invariant property of the Haar state, we get the recursive relation
\begin{equation}\label{eq:5.44}
    \begin{split}
        &h((g^*)^{d_3}b^{c_2}c^{c_3}(c^*)^{c_3}(b^*)^{c_2}g^{d_3})\\
        =&q^{2c_3-2c_2}\frac{1-q^{2c_2}}{1-q^{2c_3+2}}h((g^*)^{d_3}b^{c_2-1}c^{c_3+1}(c^*)^{c_3+1}(b^*)^{c_2-1}g^{d_3}).
    \end{split}
\end{equation}
With (\ref{eq:5.43}) as the initial condition, recursive relation (\ref{eq:5.44}) gives
\begin{equation}\label{eq:5.45}
    \begin{split}
        &h((g^*)^{d_3}b^{c_2}c^{c_3}(c^*)^{c_3}(b^*)^{c_2}g^{d_3})\\
        =&\frac{q^{2c_2c_3+2c_2+4c_3}(q^2-1)^2(q^4-1)(q^2;q^2)_{c_2}(q^2;q^2)_{c_3}}{(q^2;q^2)_{c_2+c_3+1}(1-q^{2(d_3+1)})(1-q^{2(d_3+c_2+c_3)+4})}.
    \end{split}
\end{equation}
3. Using the relation derived from the left action of $f_1\in\mathcal{U}_q(gl(3))$ on\\ $(g^*)^{d_3+1}ab^{c_2-1}c^{c_3}(c^*)^{c_3}(b^*)^{c_2}g^{d_3+1}$ and the invariant property of the Haar state, we get
\begin{equation}\label{eq:5.46}
    \begin{split}
        &q^{-d_3+3}\frac{1-q^{2d_3}}{1-q^2}h\left((g^*)^{d_3}ab^{c_2-1}c^{c_3}(c^*)^{c_3}(b^*)^{c_2}g^{d_3-1}h\right)\\
        =&q^{2c_3+2}\frac{1-q^{2c_2}}{1-q^2}h\left((g^*)^{d_3}b^{c_2-1}c^{c_3}(c^*)^{c_3}(b^*)^{c_2-1}g^{d_3}\right)\\
        &-\frac{1-q^{2c_2+2}}{1-q^2}h\left((g^*)^{d_3}b^{c_2}c^{c_3}(c^*)^{c_3}(b^*)^{c_2}g^{d_3}\right)\\
        &-q^{2c_3+2}\frac{1-q^{2c_2}}{1-q^2}h\left((g^*)^{d_3}b^{c_2-1}c^{c_3+1}(c^*)^{c_3+1}(b^*)^{c_2-1}g^{d_3}\right).
    \end{split}
\end{equation}
Substituting (\ref{eq:5.45}) into (\ref{eq:5.46}), we get
\begin{equation}\label{eq:5.47}
    \begin{split}
         &h\left((g^*)^{d_3}ab^{c_2-1}c^{c_3}(c^*)^{c_3}(b^*)^{c_2}g^{d_3-1}h\right)\\
        =&-\frac{q^{2c_2c_3+4c_2+6c_3+d_3-1}(q^2-1)^2(q^4-1)}{(1-q^{2(c_2+c_3+1)})(1-q^{2(d_3+1)})}\frac{(q^2;q^2)_{c_1}(q^2;q^2)_{c_2}}{(q^2;q^2)_{c_1+c_2}}\\
        &\left(\frac{(1-q^2)}{(1-q^{2(d_3+c_2+c_3)+2})(1-q^{2(d_3+c_2+c_3)+4})}\right).
    \end{split}
\end{equation}
Then, using the linear relations derived from the left action of $e_1, f_1\in\mathcal{U}_q(gl(3))$ on $(h^*)^{d_2}(g^*)^{d_3+1}b^{c_2}c^{c_3}(c^*)^{c_3}(b^*)^{c_2}g^{d_3}h^{d_2+1}$ and\\ $(h^*)^{d_2+1}(g^*)^{d_3}ab^{c_2-1}c^{c_3}(c^*)^{c_3}(b^*)^{c_2}g^{d_3}h^{d_2+1}$, respectively, the relation $1=aa^*+bb^*+cc^*$, and the invariant property of the Haar state, we get
\begin{equation}\label{eq:5.48}
    \begin{split}
        &h\left((h^*)^{d_2+1}(g^*)^{d_3}b^{c_2}c^{c_3}(c^*)^{c_3}(b^*)^{c_2}(g)^{d_3}h^{d_2+1}\right)\\
        =&q^{d_3-2c_2+2}\frac{1-q^{2c_2}}{1-q^{2d_3+2}}h\left((h^*)^{d_2}(g^*)^{d_3+1}ab^{c_2-1}c^{c_3}(c^*)^{c_3}(b^*)^{c_2}(g)^{d_3}h^{d_2+1}\right)\\
        &+q^{2d_3-2d_2}\frac{1-q^{2d_2+2}}{1-q^{2d_3+2}}h\left((h^*)^{d_2}(g^*)^{d_3+1}b^{c_2}c^{c_3}(c^*)^{c_3}(b^*)^{c_2}g^{d_3+1}h^{d_2}\right),
    \end{split}
\end{equation}
and
\begin{align}\label{eq:5.49}\stepcounter{equation}
        &h\left((h^*)^{d_2+1}(g^*)^{d_3}ab^{c_2-1}c^{c_3}(c^*)^{c_3}(b^*)^{c_2}g^{d_3-1}h^{d_2+2}\right)\notag\\
        =&q^{2d_3-2d_2-}\frac{1-q^{2d_2+2}}{1-q^{2d_3}}h\left((h^*)^{d_2}(g^*)^{d_3+1}ab^{c_2-1}c^{c_3}(c^*)^{c_3}(b^*)^{c_2}g^{d_3}h^{d_2+1}\right)\notag\\
        &+q^{2c_3+d_3-1}\frac{1-q^{2c_2}}{1-q^{2d_3}}h\left((h^*)^{d_2+1}(g^*)^{d_3}b^{c_2-1}c^{c_3}(c^*)^{c_3}(b^*)^{c_2-1}g^{d_3}h^{d_2+1}\right)\tag{5.50}\\
        &-q^{d_3-3}\frac{1-q^{2c_2+2}}{1-q^{2d_3}}h\left((h^*)^{d_2+1}(g^*)^{d_3}b^{c_2}c^{c_3}(c^*)^{c_3}(b^*)^{c_2}g^{d_3}h^{d_2+1}\right)\notag\\
        &-q^{2c_3+d_3-1}\frac{1-q^{2c_2}}{1-q^{2d_3}}h\left((h^*)^{d_2+1}(g^*)^{d_3}b^{c_2-1}c^{c_3+1}(c^*)^{c_3+1}(b^*)^{c_2-1}g^{d_3}h^{d_2+1}\right).\notag
\end{align}
With (\ref{eq:5.45}) and (\ref{eq:5.47}) as initial condition, recursive relation (\ref{eq:5.48}) and (\ref{eq:5.49}) gives
\begin{equation}\label{eq:5.50}
    \begin{split}
        &h\left((h^*)^{d_2}(g^*)^{d_3}b^{c_2}c^{c_3}(c^*)^{c_3}(b^*)^{c_2}(g)^{d_3}h^{d_2}\right)\\
        =&\frac{q^{2c_2c_3+2c_2+4c_3+2d_3d_2}(q^2-1)^2(q^4-1)(q^2;q^2)_{c_2}(q^2;q^2)_{d_2}}{(q^2;q^2)_{c_2+c_3+1}(q^2;q^2)_{d_2+d_3+1}(q^{2(d_3+c_2+c_3+i)+4};q^2)_{d_2+1}}\\
        &\left(\sum_{i=0}^{d_2}q^{(2c_2+2c_3+2)i}(q^2;q^2)_{c_3+d_2-i}(q^2;q^2)_{d_3+i}{d_2\choose i}_{q^2}\right),
    \end{split}
\end{equation}
and
\begin{equation}\label{eq:5.51}
    \begin{split}
        &h\left((h^*)^{d_2}(g^*)^{d_3}ab^{c_2-1}c^{c_3}(c^*)^{c_3}(b^*)^{c_2}g^{d_3-1}h^{d_2+1}\right)\\
        =&-\frac{q^{2c_2c_3+4c_2+6c_3+2d_2d_3+d_3-1}(q^2-1)^2(q^4-1)(q^2;q^2)_{c_2}(q^2;q^2)_{d_2+1}}{(q^2;q^2)_{c_2+c_3+1}(q^2;q^2)_{d_2+d_3+1}(q^{2(d_3+c_2+c_3+i)+2};q^2)_{d_2+2}}\\
        &\left(\sum_{i=0}^{d_2}q^{(2c_2+2c_3)i}(q^2;q^2)_{c_3+d_2-i}(q^2;q^2)_{d_3+i}{d_2\choose i}_{q^2}\right).
    \end{split}
\end{equation}
4. Using the linear relation derived from the left action of $e_2\in\mathcal{U}_q(gl(3))$ on $(k^*)^{k-1}(h^*)^{d_2+1-k}(g^*)^{d_3}b^{c_2-k}c^{c_3+k}(c^*)^{c_3}(b^*)^{c_2}g^{d_3}h^{d_2}$ and the invariant property of the Haar state, we get
\begin{equation}\label{eq:5.52}
    \begin{split}
        &h\left((k^*)^{k}(h^*)^{d_2-k}(g^*)^{d_3}b^{c_2-k}c^{c_3+k}(c^*)^{c_3}(b^*)^{c_2}g^{d_3}h^{d_2}\right)\\
        =&q^{c_2+d_2-4k+2-2c_3}\frac{1-q^{2(c_3+k)}}{1-q^{2(d_2-k+1)}}\\
        &\cdot h\left((k^*)^{k-1}(h^*)^{d_2+1-k}(g^*)^{d_3}b^{c_2-k+1}c^{c_3+k-1}(c^*)^{c_3}(b^*)^{c_2}g^{d_3}h^{d_2}\right)\\
        &-q^{d_2-3k+1-c_2}\frac{1-q^{2c_2}}{1-q^{2(d_2-k+1)}}\\
        &\cdot h\left((k^*)^{k-1}(h^*)^{d_2+1-k}(g^*)^{d_3}b^{c_2-k}c^{c_3+k}(c^*)^{c_3+1}(b^*)^{c_2-1}g^{d_3}h^{d_2}\right).
    \end{split}
\end{equation}
With (\ref{eq:5.50}) as initial condition, recursive relation (\ref{eq:5.52}) gives
\begin{align}\label{eq:5.53}\stepcounter{equation}
        &h\left((k^*)^{k}(h^*)^{d_2-k}(g^*)^{d_3}b^{c_2-k}c^{c_3+k}(c^*)^{c_3}(b^*)^{c_2}g^{d_3}h^{d_2}\right)\tag{5.54}\\
        =&(-1)^k\frac{q^{2c_2c_3+2c_2+4c_3+2d_3d_2+k(d_2+c_2-k+1)}(q^2-1)^2(q^4-1)(q^2;q^2)_{c_2}(q^2;q^2)_{d_2}}{(q^2;q^2)_{c_2+c_3+1}(q^2;q^2)_{d_2+d_3+1}(q^{2(d_3+c_2+c_3+i)+4};q^2)_{d_2+1}}\notag\\
        &\left(\sum_{i=0}^{d_2-k}q^{(2c_2+2c_3+2)i}(q^2;q^2)_{c_3+d_2-i}(q^2;q^2)_{d_3+i}{d_2-k\choose i}_{q^2}\right).\notag
\end{align}
5. Using the relation $1=kk^*+hh^*+gg^*$ and the modular automorphism $\rho$, we have
\begin{align}\label{eq:5.54}\stepcounter{equation}
        &h\left((k^*)^{d_1+k}(h^*)^{d_2-k}(g^*)^{d_3}b^{c_2-k}c^{c_3+k}(c^*)^{c_3}(b^*)^{c_2}g^{d_3}h^{d_2}k^{d_1}\right)\notag\\
        =&q^{-4}h\left((k^*)^{d_1+k-1}(h^*)^{d_2-k}(g^*)^{d_3}b^{c_2-k}c^{c_3+k}(c^*)^{c_3}(b^*)^{c_2}g^{d_3}h^{d_2}k^{d_1-1}\right)\tag{5.55}\\
        &-q^{-2d_1-k}h\left((k^*)^{d_1+k-1}(h^*)^{d_2+1-k}(g^*)^{d_3}b^{c_2-k}c^{c_3+k}(c^*)^{c_3}(b^*)^{c_2}g^{d_3}h^{d_2+1}k^{d_1-1}\right)\notag\\
        &-q^{-2-2d_1-2d_2}h\left((k^*)^{d_1+k-1}(h^*)^{d_2-k}(g^*)^{d_3+1}b^{c_2-k}c^{c_3+k}(c^*)^{c_3}(b^*)^{c_2}g^{d_3+1}h^{d_2}k^{d_1-1}\right).\notag
\end{align}
With (\ref{eq:5.53}) as initial condition, recursive relation (\ref{eq:5.54}) gives
\begin{align}\label{eq:5.55}\stepcounter{equation}
        &h\left((k^*)^{d_1+k}(h^*)^{d_2-k}(g^*)^{d_3}b^{c_2-k}c^{c_3+k}(c^*)^{c_3}(b^*)^{c_2}g^{d_3}h^{d_2}k^{d_1}\right)\notag\\
        =&(-1)^kq^{2d_2d_3+2d_1d_2+2d_1d_3+2c_2c_3+2c_2+4c_3+k(d_2+c_2-k+1)}\tag{5.56}\\
        &\frac{(1-q^2)^2(1-q^4)(q^2;q^2)_{d_2}(q^2;q^2)_{c_2}}{(q^2;q^2)_{c_2+c_3+1}(q^2;q^2)_{d_1+d_2+d_3+1}}\left(\sum_{j=0}^{d_1}\frac{(-1)^jq^{j^2-j}{d_1\choose j}_{q^2}{c_3\choose j}_{q^2}(q^2;q^2)_j}{(q^{2d_3+2c_2+2c_3+2d_1-2j+4};q^2)_{d_2+1}}\right.\notag\\
        &\left.\sum_{i=0}^{d_2-k}q^{(2c_2+2c_3+2d_1-2j+2)i}(q^2;q^2)_{c_3+d_2+d_1-j-i}(q^2;q^2)_{d_3+i}{d_2-k\choose i}_{q^2}\right).\notag
\end{align}

Based on (\ref{eq:5.55}), we get the following result.
\begin{theorem}
    Assume that we enumerate basis vectors of $V^R_\mu(\lambda)$ as in (\ref{eq:5.2+}) and these basis vectors take the form of (\ref{eq:5.1}.b). If $v_i,v_{i+k}\in V^R_\mu(\lambda)$, $0\le k\le\text{min}(d_2,c_2)$, takes the form (\ref{eq:5.44+}) and (\ref{eq:5.45+}), respectively, then the inner product between $v_i$ and $v_{i+k}$ with respect to the norm $\langle\cdot,\cdot\rangle_L$ is
    \begin{align}\stepcounter{equation}
            &\langle v_i,v_{i+k}\rangle_L=h(v^*_iv_{i+k})\notag\\
            =&q^{2d_2d_3+2d_1d_2+2d_1d_3+4d_3+4d_1+4d_2+2c_2c_3+4c_2+4c_3+k(d_2+c_2-k)}\tag{5.57}\\
            &\frac{(1-q^2)^2(1-q^4)(q^2;q^2)_{d_2}(q^2;q^2)_{c_2}}{(q^2;q^2)_{c_2+c_3+1}(q^2;q^2)_{d_1+d_2+d_3+1}}\left(\sum_{j=0}^{d_1}\frac{(-1)^jq^{j^2-j}{d_1\choose j}_{q^2}{c_3\choose j}_{q^2}(q^2;q^2)_j}{(q^{2d_3+2c_2+2c_3+2d_1-2j+4};q^2)_{d_2+1}}\right.\notag\\
            &\left.\sum_{i=0}^{d_2-k}q^{(2c_2+2c_3+2d_1-2j+2)i}(q^2;q^2)_{c_3+d_2+d_1-j-i}(q^2;q^2)_{d_3+i}{d_2-k\choose i}_{q^2}\right).\notag
    \end{align}
\end{theorem}

%% file: Chapters/Chapter_5/Inner_product_L.tex
\subsection{The inner products of basis vectors on $V^L_\mu(\lambda)$}
According to Theorem \ref{thm:2.1}, we pick basis vectors of $V^L_\mu(\lambda)$ as
\begin{equation}
    \xi^{\textbf{T}}\cdot\text{det}_q^{-\lambda_3}.
\end{equation}
Similar to $V^R_\mu(\lambda)$, it is suffice to consider the case $\lambda=(\lambda_1,\lambda_2,0)$. Notice that $\gamma(\xi^I_J)=\xi^J_I$ where $\gamma$ is the diagonal flip homomorphism in (\ref{eq:3.50}). Hence, we have $\xi^{\textbf{T}}=\gamma(\xi_{\textbf{T}})$. On $\mathcal{O}(U_q(3))$, there are two types of basis vectors of $V^L_\mu(\lambda)$ depending on $\lambda$ and $\mu$ corresponding to (\ref{eq:5.1}.a) and (\ref{eq:5.1}.b), respectively.
\begin{enumerate}
    \item[1)] $\xi^{\textbf{T}}=(k^*)^{d_1}(-q^{-1}f^*)^{d_2}a^{c_1}d^{c_2}g^{c_3}D_q^{d_1+d_2}$, and
    \item[2)] $\xi^{\textbf{T}}=(k^*)^{d_1}(-q^{-1}f^*)^{d_2}(q^{-2}c^*)^{d_3}d^{c_2}g^{c_3}D_q^{d_1+d_2+d_3}$.
\end{enumerate}
We enumerate basis vectors of $V^L_\mu(\lambda)$ using $v_i$ in (\ref{eq:5.2+}) as
\begin{equation}\label{eq:5.63}
    v^i=\gamma(v_i).
\end{equation}
According to Theorem \ref{thm:3.10}, we have
\begin{equation}\label{eq:5.64}
    \langle v^{i+k},v^i\rangle_R=h(v^{i+k}(v^{i})^*)=h(\gamma(v^{i+k})\gamma((v^{i})^*))=h(v_{i+k}\gamma((v^{i})^*)).
\end{equation}
If $v_i$ takes the form (\ref{eq:5.1}.a), then
\begin{equation}
    \gamma((v^{i})^*)=q^{-4c_3-2c_2-2d_2}v_i^*.
\end{equation}
If $v_i$ takes the form (\ref{eq:5.1}.b), then
\begin{equation}
    \gamma((v^{i})^*)=q^{-4c_3-2c_2-4d_3-2d_2}v_i^*.
\end{equation}
Combining (\ref{eq:5.64}) with (\ref{eq:5.35+}) and (\ref{eq:5.55}), we get the following result.
\begin{theorem}
    Assume that we enumerate basis vectors of $V^R_\mu(\lambda)$ as in (\ref{eq:5.63}). If $v_i,v_{i+k}$ takes the form (\ref{eq:5.4++}) and (\ref{eq:5.5+}), respectively, then the inner product with respect to $\langle\cdot,\cdot\rangle_R$ between corresponding $v^i$ and $v^{i+k}$ is
    \begin{equation}
        \begin{split}
            &\langle v^{i+k},v^i\rangle_R=q^{-4c_3-2c_2-2d_2}h(v_{i+k}v^*_i)\\
            =&q^{2d_1d_2+2c_1c_2+2c_1c_3+2c_2c_3+k(d_2+c_2-k)}\frac{(1-q^2)^2(1-q^4)(q^2;q^2)_{d_2}(q^2;q^2)_{c_2}}{(q^2;q^2)_{d_1+d_2+1}(q^2;q^2)_{c_1+c_2+c_3+1}}\\
            &\left(\sum_{j=0}^{c_3}\frac{(-1)^{j}q^{j^2-j}{d_1\choose j}_{q^{2}}{c_3\choose j}_{q^2}(q^2;q^2)_{j}}{(q^{2d_1+2d_2+2c_1+2c_3-2j+4};q^2)_{c_2+1}}\right.\\
            &\left.\cdot\sum_{i=0}^{c_2-k}q^{(2d_1+2d_2+2c_3-2j+2)i}(q^2;q^2)_{c_1+i}(q^2;q^2)_{d_1+c_2+c_3-j-i}{c_2-k\choose i}_{q^2}\right).
        \end{split}
    \end{equation}
    If $v_i,v_{i+k}$ takes the form (\ref{eq:5.44+}) and (\ref{eq:5.45+}), respectively, then the inner product with respect to $\langle\cdot,\cdot\rangle_R$ between corresponding $v^i$ and $v^{i+k}$ is
    \begin{align}\stepcounter{equation}
            &\langle v^{i+k},v^i\rangle_R=q^{-4c_3-2c_2-4d_3-2d_2}h(v_{i+k}v^*_i)\notag\\
            =&q^{2d_2d_3+2d_1d_2+2d_1d_3+2c_2c_3+k(d_2+c_2-k)}\frac{(1-q^2)^2(1-q^4)(q^2;q^2)_{d_2}(q^2;q^2)_{c_2}}{(q^2;q^2)_{c_2+c_3+1}(q^2;q^2)_{d_1+d_2+d_3+1}}\notag\\
            &\left(\sum_{j=0}^{d_1}\frac{(-1)^jq^{j^2-j}{d_1\choose j}_{q^2}{c_3\choose j}_{q^2}(q^2;q^2)_j}{(q^{2d_3+2c_2+2c_3+2d_1-2j+4};q^2)_{d_2+1}}\right.\tag{5.64}\\
            &\left.\left(\sum_{i=0}^{d_2-k}q^{(2c_2+2c_3+2d_1-2j+2)i}(q^2;q^2)_{c_3+d_2+d_1-j-i}(q^2;q^2)_{d_3+i}{d_2-k\choose i}_{q^2}\right)\right).\notag
    \end{align}
\end{theorem}


%% file: Chapters/Appendix_A.tex
\section{The formula for direct computation}\label{append_1}
In this appendix, we give the formula to compute the the Haar state value of $(k^*)^{d_1+k}(h^*)^{d_2-k}a^{c_1}b^{c_2-k}c^{c_3+k}(c^*)^{c_3}(b^*)^{c_2}(a^*)^{c_1}h^{d_2}k^{d_1}$. Readers may compare this formula with (\ref{eq:5.35+}). Before actual computation, we prove the following lemma which will be used in the computation.
\begin{lemma}\label{lemma:append}
    Let $x_{i,k}, x_{j,l}, x_{j,k}, x_{i,l}\in\mathcal{O}(GL_q(n))$ with $i<j$ and $k<l$.\\
    For $s,t\in\mathbb{N}$, we have 
    \begin{equation}\label{eq:app1.1}
        \begin{split}
            x_{j,l}^sx_{i,k}^t=\sum_{n}q^{3n^2-2(i+j)n}{s\choose n}_{q^2}{t\choose n}_{q^2}(q^2;q^2)_nx_{i,k}^{t-n}(x_{j,k}x_{i,l})^{n}x_{j,l}^{s-n}.
        \end{split}
    \end{equation}
\end{lemma}
\begin{proof}
     Without loss of generality, we may take $x_{i,k}=a$, $x_{j,l}=e$, $x_{i,l}=b$, and $x_{j,k}=d$ on $\mathcal{O}(GL_q(3))$ so that the computation steps are easy to read.\\
     To prove (\ref{eq:app1.1}), we consider the case $s\ge t$ first. We have
    \begin{equation}\label{eq:app1.2}
        \begin{split}
            e^sa^t&=(ae^i)a^{j-1}-(q-q^{-1})\left(\sum_{k=0}^{s-1}e^{k}bde^{s-1-k}\right)a^{t-1}\\
            &=(ae-q(1-q^{-2s})\cdot bd)e^{s-1}a^{t-1}\\
            &=\left(\prod_{k=0}^{t-1}(ae-q(1-q^{-2(s-k)})bd)\right)e^{s-t}.
        \end{split}
    \end{equation}
    Notice that
    \begin{equation}\label{eq:app1.3}
        \begin{split}
            \prod_{k=0}^{t-1}(ae-q(1-q^{-2(s-k)})bd)=\prod_{k=0}^{t-1}\left([ae-q\cdot bd]+q^{2k}(q^{1-2s}bd)\right).
        \end{split}
    \end{equation}
    Hence in the expansion of (\ref{eq:app1.3}) the coefficient of $(ae-q\cdot bd)^{t-p}(q^{1-2s}bd)^p$ is the same as the coefficient of $z^p$ in the expansion of $(-z;q^2)_t$. Recall that 
    \begin{equation}
        (z;q)_t=\sum_{j=0}^t{t\choose j}q^{j(j-1)/2}(-z)^j.
    \end{equation}
    Therefore, we have
    \begin{equation}\label{eq:app1.5}
        \begin{split}
            &\prod_{k=0}^{t-1}(ae-q(1-q^{-2(s-k)})bd)\\
            =&\sum_{p=0}^t{t\choose p}q^{p(p-1)}(ae-q\cdot bd)^{t-p}(q^{1-2s}bd)^p\\
            =&\sum_{p=0}^t{t\choose p}_{q^2}q^{p^2-2sp}\left(\sum_{r=0}^{t-p}(-q)^{(1-2r)(t-p-r)}{t-p\choose r}_{q^2}a^r(bd)^{t-p-r}e^r\right)bd^p.
        \end{split}
    \end{equation}
    By further simplification of (\ref{eq:app1.5}), we get
    \begin{equation}
        \begin{split}
            &\prod_{k=0}^{t-1}(ae-q(1-q^{-2(s-k)})bd)\\
            =&\sum_{p=0}^t{t\choose p}_{q^2}q^{p^2-2sp}\left(\sum_{r=0}^{t-p}(-q)^{(1-2r)(t-r)-p}{t-p\choose r}_{q^2}a^r(bd)^{t-r}e^r\right)\\
            =&\sum_{r=0}^t(-q)^{(1-2r)(t-r)}{t\choose r}_{q^2}\left(\sum_{p=0}^{t-r}q^{p^2-p}{t-r\choose p}_{q^2}(-q^{-2s})^p\right)a^r(bd)^{t-r}e^r\\
            =&\sum_{r=0}^t(-q)^{(1-2r)(t-r)}{t\choose r}_{q^2}(q^{-2s};q^2)_{t-r}a^r(bd)^{t-r}e^r
        \end{split}
    \end{equation}
    Hence, we get
    \begin{equation}\label{eq:app1.7}
        \begin{split}
            e^sa^t=&\sum_{r=0}^t(-q)^{(1-2r)(t-r)}{t\choose r}_{q^2}(q^{-2s};q^2)_{t-r}a^r(bd)^{t-r}e^{r+s-t}\\
            =&\sum_{n=0}^{t}q^{3n^2-2(s+t)n}{s\choose n}_{q^2}{t\choose n}_{q^2}(q^2;q^2)_na^{t-n}(bd)^{n}e^{s-n}.
        \end{split}
    \end{equation}
     For the case $s<t$, we have
    \begin{equation}
        \begin{split}
            e^sa^t&=e^{s-1}a^te-(q-q^{-1})\sum_{k=0}^{t-1}e^{s-1}a^{t-1-k}bda^{k}\\
            &=e^{s-1}a^{t-1}(ae-q(1-q^{-2t})bd)\\
            &=a^{t-s}\prod_{k=0}^{s-1}(ae-q(1-q^{-2t+2k})bd).
        \end{split}
    \end{equation}
    Through similar computation steps, we find that
    \begin{equation}\label{eq:app1.9}
        e^sa^t=\sum_{n=0}^{s}q^{3n^2-2(s+t)n}{s\choose n}_{q^2}{t\choose n}_{q^2}(q^2;q^2)_na^{t-n}(bd)^{n}e^{s-n}.
    \end{equation}
    Since (\ref{eq:app1.7}) and (\ref{eq:app1.9}) are consistent, we get (\ref{eq:app1.1}).
\end{proof}
Recall that $ak^*=k^*a$, $ah^*=h^*a$, $c^*d=dc^*$. We have
\begin{equation}
    \begin{split}
        &h((k^*)^{d_1+k}(h^*)^{d_2-k}a^{c_1}b^{c_2-k}c^{c_3+k}(c^*)^{c_3}(b^*)^{c_2}(a^*)^{c_1}h^{d_2}k^{d_1})\\
        =&h(a^{c_1}(k^*)^{d_1+k}(h^*)^{d_2-k}b^{c_2-k}c^{c_3+k}(c^*)^{c_3}(b^*)^{c_2}(a^*)^{c_1}h^{d_2}k^{d_1})\\
        =&(-q)^{-d_2+c_2}\sum_{j=0}^{d_2-k}(-q)^{(1-2j)(d_2-k-j)}{d_2-k\choose j}_{q^{2}}\sum_{l=0}^{c_2}(-q)^{(1-2l)(c_2-l)}{c_2\choose l}_{q^2}\\
        &h(a^{c_1}(k^*)^{d_1+k}\cdot a^j(cd)^{d_2-k-j}f^j\cdot b^{c_2-k}c^{c_3+k}\\
        &\cdot (c^*)^{c_3}\cdot d^l(fg)^{c_2-l}k^l\cdot (a^*)^{c_1} h^{d_2}k^{d_1})\\
        =&(-q)^{-d_2+c_2}\sum_{j=0}^{d_2-k}(-q)^{(1-2j)(d_2-k-j)}{d_2-k\choose j}_{q^{2}}\sum_{l=0}^{c_2}(-q)^{(1-2l)(c_2-l)}{c_2\choose l}_{q^2}\\
        &h(a^{c_1+j}(k^*)^{d_1+k}\cdot (cd)^{d_2-k-j}f^j\cdot b^{c_2-k}c^{c_3+k}\\
        &\cdot d^l(c^*)^{c_3}\cdot (fg)^{c_2-l}\cdot (a^*)^{c_1} k^lh^{d_2}k^{d_1})\\
        =&(-q)^{-d_2+c_2+2c_3}\sum_{i=0}^{d_1+k}(-q)^{(1-2i)(d_1+k-i)}{d_1+k\choose i}_{q^{2}}\\
        &\sum_{j=0}^{d_2-k}(-q)^{(1-2j)(d_2-k-j)}{d_2-k\choose j}_{q^{2}}\sum_{r=0}^{c_3}(-q)^{(1-2r)(c_3-r)}{c_3\choose r}_{q^2}\\
        &\sum_{l=0}^{c_2}(-q)^{(1-2l)(c_2-l)}{c_2\choose l}_{q^2}\sum_{m=0}^{c_1}(-q)^{(1-2m)(c_1-m)}{c_1\choose m}_{q^2}\\
        &\sum_{s}q^{3s^2-2(j+c_2-k)s}{j\choose s}_{q^2}{c_2-k\choose s}_{q^2}(q^2;q^2)_s\\
        &q^{-i(d_2-k-j)-(c_2-k-s)(i+d_2-k-j)-(c_3+k)(j-s)-(r+l)(i+j)-(c_3-r)(j-s)}\\
        &q^{-r(c_2-l+m)-m(j-s+c_2-l)-d_2(m+l)}\\
        &h(a^{i+j+c_1}b^{d_1+k-i+c_2-k-s}c^{d_2-j+s+c_3}d^{d_1-i+d_2-j+r+l}\\
        &\cdot e^{i+s+c_3-r+m}f^{j-s+c_2-l+c_1-m}g^{c_3-r+c_2-l}h^{c_1-m+r+d_2}k^{m+l+d_1})\\
    \end{split}
\end{equation}

%% file: Chapters/Appendix_B.tex
\section{Proof of Proposition \ref{prop:5.4}}\label{append_a}
In this appendix, we give the computation steps to prove (\ref{eq:5.38}) and (\ref{eq:5.39}).

Assume that (\ref{eq:5.38}) and (\ref{eq:5.39}) hold in the case $c_2=n$. First, we compute $h\left((k^*)^{d_1}(h^*)^{d_2}a^{c_1}b^{n+1}(b^*)^{n+1}(a^*)^{c_1}h^{d_2}k^{d_1}\right)$ using (\ref{eq:5.34}). We have
\begin{align}\label{eq:a.1}\stepcounter{equation}
        &q^{c_1}\frac{1-q^{-2(c_1+1)}}{1-q^{-2}}h\left((k^*)^{d_1}(h^*)^{d_2}a^{c_1}b^{n+1}(b^*)^{n+1}(a^*)^{c_1}h^{d_2}k^{d_1}\right)\notag\\
        =&q^{c_1+2}\frac{1-q^{-2(n+1)}}{1-q^{-2}}h\left((k^*)^{d_1}(h^*)^{d_2}a^{c_1+1}b^{n}(b^*)^{n}(a^*)^{c_1+1}h^{d_2}k^{d_1}\right)\notag\\
        &+\frac{1-q^{-2d_2}}{1-q^{-2}}h\left((k^*)^{d_1}(h^*)^{d_2}a^{c_1}b^{n+1}(b^*)^{n}(a^*)^{c_1+1}gh^{d_2-1}k^{d_1}\right)\tag{B.1}\\
        =&\frac{1}{1-q^{2}}\frac{q^{2d_1d_2+2c_1n+2n+c_1+2}(1-q^2)^2(1-q^4)(q^2;q^2)_{d_2}(q^2;q^2)_{n+1}}{(q^2;q^2)_{d_1+d_2+1}(q^2;q^2)_{c_1+n+2}(q^{2(d_1+d_2+c_1)+4};q^2)_{n+2}}\notag\\
        &\left((1-q^{2(d_1+d_2+c_1)+4})\left(\sum_{i=0}^{n}q^{2(d_1+d_2+1)i}(q^2;q^2)_{d_1+n-i}(q^2;q^2)_{c_1+1+i}{n\choose i}_{q^2}\right)\right.\notag\\
        &\left.-q^{2d_1+2+2n}(1-q^{2d_2})\left(\sum_{i=0}^{n}q^{2(d_1+d_2)i}(q^2;q^2)_{d_1+n-i}(q^2;q^2)_{c_1+1+i}{n\choose i}_{q^2}\right)\right).\notag
\end{align}
Denote the difference of in the last two lines of (\ref{eq:a.1}) as $D_1$. We have
\begin{equation}
    \begin{split}
        &D_1=\\
        &(1-q^{2(d_1+d_2+c_1)+4})\left(\sum_{i=0}^{n}q^{2(d_1+d_2+1)i}(q^2;q^2)_{d_1+n-i}(q^2;q^2)_{c_1+1+i}{n\choose i}_{q^2}\right)\\
        &-q^{2d_1+2+2n}(1-q^{2d_2})\left(\sum_{i=0}^{n}q^{2(d_1+d_2)i}(q^2;q^2)_{d_1+n-i}(q^2;q^2)_{c_1+1+i}{n\choose i}_{q^2}\right)\\
        =&\left(\sum_{i=0}^{n}q^{2(d_1+d_2+1)i}(q^2;q^2)_{d_1+n-i}(q^2;q^2)_{c_1+1+i}{n\choose i}_{q^2}\right)\\
        &-q^{2c_1+2}\left(\sum_{i=0}^{n}q^{2(d_1+d_2+1)(i+1)}(q^2;q^2)_{d_1+n-i}(q^2;q^2)_{c_1+1+i}{n\choose i}_{q^2}\right)\\
        &-\left(\sum_{i=0}^{n}q^{2(d_1+d_2+1)i}q^{2d_1+2+2n-2i}(q^2;q^2)_{d_1+n-i}(q^2;q^2)_{c_1+1+i}{n\choose i}_{q^2}\right)\\
        &+\left(\sum_{i=0}^{n}q^{2(d_1+d_2+1)(i+1)}q^{2n-2i}(q^2;q^2)_{d_1+n-i}(q^2;q^2)_{c_1+1+i}{n\choose i}_{q^2}\right)\\
    \end{split}
\end{equation}
\begin{equation}\label{eq:a.3}
    \begin{split}
        D_1=&\left(\sum_{i=0}^{n}q^{2(d_1+d_2+1)i}(q^2;q^2)_{d_1+n+1-i}(q^2;q^2)_{c_1+1+i}{n\choose i}_{q^2}\right)\\
        &-q^{2c_1+2}\left(\sum_{i=1}^{n+1}q^{2(d_1+d_2+1)i}(q^2;q^2)_{d_1+n+1-i}(q^2;q^2)_{c_1+i}{n\choose i-1}_{q^2}\right)\\
        &+\left(\sum_{i=1}^{n+1}q^{2(d_1+d_2+1)i}q^{2n+2-2i}(q^2;q^2)_{d_1+n+1-i}(q^2;q^2)_{c_1+i}{n\choose i-1}_{q^2}\right)\\
         =&\left(\sum_{i=0}^{n}q^{2(d_1+d_2+1)i}(q^2;q^2)_{d_1+n+1-i}(q^2;q^2)_{c_1+i}{n\choose i}_{q^2}\right)\\
         &-q^{2c_1+2}\left(\sum_{i=0}^{n}q^{2(d_1+d_2+1)i}q^{2i}(q^2;q^2)_{d_1+n+1-i}(q^2;q^2)_{c_1+i}{n\choose i}_{q^2}\right)\\
        &-q^{2c_1+2}\left(\sum_{i=1}^{n+1}q^{2(d_1+d_2+1)i}(q^2;q^2)_{d_1+n+1-i}(q^2;q^2)_{c_1+i}{n\choose i-1}_{q^2}\right)\\
        &+\left(\sum_{i=1}^{n+1}q^{2(d_1+d_2+1)i}q^{2n+2-2i}(q^2;q^2)_{d_1+n+1-i}(q^2;q^2)_{c_1+i}{n\choose i-1}_{q^2}\right)\\
        =&(1-q^{2c_1+2})\left(\sum_{i=0}^{n+1}q^{2(d_1+d_2+1)i}(q^2;q^2)_{d_1+n+1-i}(q^2;q^2)_{c_1+i}{n+1\choose i}_{q^2}\right)
    \end{split}
\end{equation}
Substituting (\ref{eq:a.3}) into (\ref{eq:a.1}) and after simplification, we get
\begin{equation}\label{eq:a.4}
    \begin{split}
        &h\left((k^*)^{d_1}(h^*)^{d_2}a^{c_1}b^{n+1}(b^*)^{n+1}(a^*)^{c_1}h^{d_2}k^{d_1}\right)\\
        =&\frac{q^{2d_1d_2+2c_1(n+1)+2n+2}(1-q^2)^2(1-q^4)(q^2;q^2)_{d_2}(q^2;q^2)_{n+1}}{(q^2;q^2)_{d_1+d_2+1}(q^2;q^2)_{c_1+n+2}(q^{2(d_1+d_2+c_1)+4};q^2)_{n+2}}\\
        &\left(\sum_{i=0}^{n+1}q^{2(d_1+d_2+1)i}(q^2;q^2)_{d_1+n+1-i}(q^2;q^2)_{c_1+i}{n+1\choose i}_{q^2}\right).
    \end{split}
\end{equation}
(\ref{eq:a.4}) is consistent with (\ref{eq:5.38}) when we put $c_2=n+1$. Then, we compute $h\left((k^*)^{d_1}(h^*)^{d_2}a^{c_1}b^{n+2}(b^*)^{n+1}(a^*)^{c_1+1}gh^{d_2-1}k^{d_1}\right)$ using (\ref{eq:5.37}). We have
\begin{align}\stepcounter{equation}
    &q^{-c_1}\frac{1-q^{2(c_1+1)}}{1-q^2}h\left((k^*)^{d_1}(h^*)^{d_2}a^{c_1}b^{n+2}(b^*)^{n+1}(a^*)^{c_1+1}gh^{d_2-1}k^{d_1}\right)\tag{B.5}\\
    =&q^{c_1-2n+1}\frac{1-q^{2n+2}}{1-q^2}h\left((k^*)^{d_1}(h^*)^{d_2}a^{c_1+1}b^{n+1}(b^*)^{n}(a^*)^{c_1+2}gh^{d_2-1}k^{d_1}\right)\tag{$p_1$}\\
    &+q^{2d_1}\frac{1-q^{2d_2}}{1-q^2}h\left((k^*)^{d_1}(h^*)^{d_2-1}a^{c_1+1}b^{n+1}(b^*)^{n+1}(a^*)^{c_1+1}h^{d_2-1}k^{d_1}\right)\tag{$p_2$}\\
    &-\frac{1-q^{2d_2+2}}{1-q^2}h\left((k^*)^{d_1}(h^*)^{d_2}a^{c_1+1}b^{n+1}(b^*)^{n+1}(a^*)^{c_1+1}h^{d_2}k^{d_1}\right)\tag{$p_3$}\\
    &-q^{2(d_1+2)}\frac{1-q^{2d_2}}{1-q^2}h\left((k^*)^{d_1+1}(h^*)^{d_2-1}a^{c_1+1}b^{n+1}(b^*)^{n+1}(a^*)^{c_1+1}h^{d_2-1}k^{d_1+1}\right)\tag{$p_4$}
\end{align}
We will compute $p_1+p_3+p_4$ first. We have
\begin{align}\label{eq:a.6}\stepcounter{equation}
        &p_1+p_3+p_4\tag{B.6}\\
        =&\frac{-1}{1-q^2}\frac{q^{2d_1d_2+2(c_1+1)(n+1)+2(n+1)}(1-q^2)^2(1-q^4)(q^2;q^2)_{d_2}(q^2;q^2)_{n+1}}{(q^2;q^2)_{c_1+n+3}(q^2;q^2)_{d_1+d_2+1}(q^{2(d_1+d_2+c_1+1)+4};q^2)_{n+2}}\notag\\
        &\left(q^{2d_1+2d_2}(1-q^{2n+2})\left(\sum_{i=0}^{n}q^{2(d_1+d_2)i}(q^2;q^2)_{d_1+n-i}(q^2;q^2)_{c_1+2+i}{n\choose i}_{q^2}\right)\right.\notag\\
        &+(1-q^{2d_2+2})\left(\sum_{i=0}^{n+1}q^{2(d_1+d_2+1)i}(q^2;q^2)_{d_1+n+1-i}(q^2;q^2)_{c_1+1+i}{n+1\choose i}_{q^2}\right)\notag\\
        &\left.+q^{2d_2+2}\left(\sum_{i=0}^{n+1}q^{2(d_1+d_2+1)i}(q^2;q^2)_{d_1+n+2-i}(q^2;q^2)_{c_1+1+i}{n+1\choose i}_{q^2}\right)\right).\notag
\end{align}
Denote the sum in the last three lines of (\ref{eq:a.6}) as $S_1$. We have
\begin{align}\label{eq:a.7}\stepcounter{equation}
        &S_1\tag{B.7}\\
        =&q^{2d_1+2d_2}(1-q^{2n+2})\left(\sum_{i=0}^{n}q^{2(d_1+d_2)i}(q^2;q^2)_{d_1+n-i}(q^2;q^2)_{c_1+2+i}{n\choose i}_{q^2}\right)\notag\\
        &+(1-q^{2d_2+2})\left(\sum_{i=0}^{n+1}q^{2(d_1+d_2+1)i}(q^2;q^2)_{d_1+n+1-i}(q^2;q^2)_{c_1+1+i}{n+1\choose i}_{q^2}\right)\notag\\
        &+q^{2d_2+2}\left(\sum_{i=0}^{n+1}q^{2(d_1+d_2+1)i}(q^2;q^2)_{d_1+n+2-i}(q^2;q^2)_{c_1+1+i}{n+1\choose i}_{q^2}\right)\notag\\
        =&(1-q^{2n+2})\left(\sum_{i=0}^{n}q^{2(d_1+d_2)(i+1)}(q^2;q^2)_{d_1+n-i}(q^2;q^2)_{c_1+2+i}{n\choose i}_{q^2}\right)\notag\\
        &+\left(\sum_{i=0}^{n+1}q^{2(d_1+d_2+1)i}(q^2;q^2)_{d_1+n+1-i}(q^2;q^2)_{c_1+1+i}{n+1\choose i}_{q^2}\right)\notag\\
        &-q^{2d_1+2d_2+2n+6}\left(\sum_{i=0}^{n+1}q^{2(d_1+d_2)i}(q^2;q^2)_{d_1+n+1-i}(q^2;q^2)_{c_1+1+i}{n+1\choose i}_{q^2}\right)\notag\\
        =&(1-q^{2n+2})\left(\sum_{i=1}^{n+1}q^{2(d_1+d_2)i}(q^2;q^2)_{d_1+n+1-i}(q^2;q^2)_{c_1+1+i}{n\choose i-1}_{q^2}\right)\notag\\
        &+\left(\sum_{i=0}^{n+1}q^{2(d_1+d_2)i}(q^2;q^2)_{d_1+n+1-i}(q^2;q^2)_{c_1+1+i}\cdot q^{2i}{n+1\choose i}_{q^2}\right)\notag\\
        &-q^{2d_1+2d_2+2n+6}\left(\sum_{i=0}^{n+1}q^{2(d_1+d_2)i}(q^2;q^2)_{d_1+n+1-i}(q^2;q^2)_{c_1+1+i}{n+1\choose i}_{q^2}\right)\notag\\
        =&(1-q^{2d_1+2d_2+2n+6})\left(\sum_{i=0}^{n+1}q^{2(d_1+d_2)i}(q^2;q^2)_{d_1+n+1-i}(q^2;q^2)_{c_1+1+i}{n+1\choose i}_{q^2}\right).\notag
\end{align}
Substituting (\ref{eq:a.6}) and (\ref{eq:a.7}) back to (B.5), we get
\begin{align}\stepcounter{equation}
        &q^{-c_1}\frac{1-q^{2(c_1+1)}}{1-q^2}h\left((k^*)^{d_1}(h^*)^{d_2}a^{c_1}b^{n+2}(b^*)^{n+1}(a^*)^{c_1+1}gh^{d_2-1}k^{d_1}\right)\tag{B.8}\\
        =&\frac{1}{1-q^2}\frac{q^{2d_1d_2+2(c_1+1)(n+1)+2n+2}(1-q^2)^2(1-q^4)(q^2;q^2)_{d_2}(q^2;q^2)_{n+1}}{(q^2;q^2)_{d_1+d_2}(q^2;q^2)_{c_1+n+3}(q^{2(d_1+d_2+c_1)+4};q^2)_{n+2}}\notag\\
        &\left(\sum_{i=0}^{n+1}q^{2(d_1+d_2)i}(q^2;q^2)_{d_1+n+1-i}(q^2;q^2)_{c_1+1+i}{n+1\choose i}_{q^2}\right)\notag\\
        &\frac{-1}{1-q^2}\frac{q^{2d_1d_2+2(c_1+1)(n+1)+2(n+1)}(1-q^2)^2(1-q^4)(q^2;q^2)_{d_2}(q^2;q^2)_{n+1}}{(q^2;q^2)_{c_1+n+3}(q^2;q^2)_{d_1+d_2+1}(q^{2(d_1+d_2+c_1+1)+4};q^2)_{n+2}}\notag\\
        &(1-q^{2d_1+2d_2+2n+6})\left(\sum_{i=0}^{n+1}q^{2(d_1+d_2)i}(q^2;q^2)_{d_1+n+1-i}(q^2;q^2)_{c_1+1+i}{n+1\choose i}_{q^2}\right)\notag\\
        =&-\frac{(1-q^{2c_1+2})}{1-q^2}\frac{q^{2d_1d_2+2d_1+2d_2+2+2c_1(n+1)+4(n+1)}(1-q^2)^2(1-q^4)}{(q^2;q^2)_{c_1+n+3}(q^2;q^2)_{d_1+d_2+1}(q^{2(d_1+d_2+c_1)+4};q^2)_{n+3}}\notag\\
        &(q^2;q^2)_{d_2}(q^2;q^2)_{n+2}\left(\sum_{i=0}^{n+1}q^{2(d_1+d_2)i}(q^2;q^2)_{d_1+n+1-i}(q^2;q^2)_{c_1+1+i}{n+1\choose i}_{q^2}\right).\notag
\end{align}
Hence, we get
\begin{equation}
    \begin{split}
        &h\left((k^*)^{d_1}(h^*)^{d_2}a^{c_1}b^{n+2}(b^*)^{n+1}(a^*)^{c_1+1}gh^{d_2-1}k^{d_1}\right)\\
        =&\frac{q^{2d_1d_2+2d_1+2d_2+c_1+2+2c_1(n+1)+4(n+1)}(1-q^2)^2(1-q^4)(q^2;q^2)_{d_2}(q^2;q^2)_{n+2}}{(q^2;q^2)_{c_1+n+3}(q^2;q^2)_{d_1+d_2+1}(q^{2(d_1+d_2+c_1)+4};q^2)_{n+3}}\\
        &\left(\sum_{i=0}^{n+1}q^{2(d_1+d_2)i}(q^2;q^2)_{d_1+n+1-i}(q^2;q^2)_{c_1+1+i}{n+1\choose i}_{q^2}\right)
    \end{split}
\end{equation}
which is consistent with (\ref{eq:5.39}) when we put $c_2=n+1$. Therefore, (\ref{eq:5.38}) and (\ref{eq:5.39}) still hold in the case $c_2=n+1$ and we proved Proposition \ref{prop:5.4} by induction.

%% file: Chapters/Appendix_C.tex
\section{Proof of Proposition \ref{prop:5.6}}\label{append_b}
In this appendix we give the computation steps to prove (\ref{eq:5.35+}).

Assume that (\ref{eq:5.35+}) holds in the case $c_3=n$. We compute the case $c_3=n+1$ following recursive relation (\ref{eq:5.41+}). We have
\begin{align}\label{eq:b.1}\stepcounter{equation}
        &h((k^*)^{d_1+k}(h^*)^{d_2-k}a^{c_1}b^{c_2-k}c^{n+1+k}(c^*)^{n+1}(b^*)^{c_2}(a^*)^{c_1}h^{d_2}k^{d_1})\tag{C.1}\\
        =&h((k^*)^{d_1+k}(h^*)^{d_2-k}a^{c_1}b^{c_2-k}c^{n+k}(c^*)^{n}(b^*)^{c_2}(a^*)^{c_1}h^{d_2}k^{d_1})\notag\\
        &-q^{-2(c_2+n)}h((k^*)^{d_1+k}(h^*)^{d_2-k}a^{c_1+1}b^{c_2-k}c^{n+k}(c^*)^{n}(b^*)^{c_2}(a^*)^{c_1+1}h^{d_2}k^{d_1})\notag\\
        &-q^{-(2n+k)}h((k^*)^{d_1+k}(h^*)^{d_2-k}a^{c_1}b^{c_2+1-k}c^{n+k}(c^*)^{n}(b^*)^{c_2+1}(a^*)^{c_1}h^{d_2}k^{d_1})\notag\\
        =&(-1)^kq^{2d_1d_2+2c_1c_2+2c_1n+2c_2n+2c_2+4n+k(d_2+c_2-k+1)}\notag\\
        &\frac{(1-q^2)^2(1-q^4)(q^2;q^2)_{d_2}(q^2;q^2)_{c_2}}{(q^2;q^2)_{d_1+d_2+1}(q^2;q^2)_{c_1+c_2+n+1}}\left(\sum_{j=0}^{n}\frac{(-1)^{j}q^{j^2-j}{d_1\choose j}_{q^{2}}{n\choose j}_{q^2}(q^2;q^2)_{j}}{(q^{2d_1+2d_2+2c_1+2n-2j+4};q^2)_{c_2+1}}\right.\notag\\
        &\left.\cdot\underbrace{\sum_{i=0}^{c_2-k}q^{(2d_1+2d_2+2n-2j+2)i}(q^2;q^2)_{c_1+i}(q^2;q^2)_{d_1+c_2+n-j-i}{c_2-k\choose i}_{q^2}}_{p_1(j)}\right)\notag\\
        &-(-1)^kq^{2d_1d_2+2c_1c_2+2c_1n+2c_2n+2c_2+4n+k(d_2+c_2-k+1)}\notag\\
        &\frac{(1-q^2)^2(1-q^4)(q^2;q^2)_{d_2}(q^2;q^2)_{c_2}}{(q^2;q^2)_{d_1+d_2+1}(q^2;q^2)_{c_1+1+c_2+n+1}}\left(\sum_{j=0}^{n}\frac{(-1)^{j}q^{j^2-j}{d_1\choose j}_{q^{2}}{n\choose j}_{q^2}(q^2;q^2)_{j}}{(q^{2d_1+2d_2+2c_1+2+2n-2j+4};q^2)_{c_2+1}}\right.\notag\\
        &\left.\cdot\underbrace{\sum_{i=0}^{c_2-k}q^{(2d_1+2d_2+2n-2j+2)i}(q^2;q^2)_{c_1+1+i}(q^2;q^2)_{d_1+c_2+n-j-i}{c_2-k\choose i}_{q^2}}_{p_2(j)}\right)\notag\\
        &-(-1)^kq^{2d_1d_2+2c_1(c_2+1)+2c_1n+2c_2n+2c_2+2+4n+k(d_2+c_2-k+1)}\notag\\
        &\frac{(1-q^2)^2(1-q^4)(q^2;q^2)_{d_2}(q^2;q^2)_{c_2+1}}{(q^2;q^2)_{d_1+d_2+1}(q^2;q^2)_{c_1+c_2+1+n+1}}\left(\sum_{j=0}^{n}\frac{(-1)^{j}q^{j^2-j}{d_1\choose j}_{q^{2}}{n\choose j}_{q^2}(q^2;q^2)_{j}}{(q^{2d_1+2d_2+2c_1+2n-2j+4};q^2)_{c_2+2}}\right.\notag\\
        &\left.\cdot\underbrace{\sum_{i=0}^{c_2+1-k}q^{(2d_1+2d_2+2n-2j+2)i}(q^2;q^2)_{c_1+i}(q^2;q^2)_{d_1+c_2+1+n-j-i}{c_2+1-k\choose i}_{q^2}}_{p_3(j)}\right).\notag
\end{align}
For a fixed index $j$, we first compute $p_2(j)+p_3(j)$. The common factor $p_2(j)$ and $p_3(j)$ is 
\begin{align}\stepcounter{equation}
        cf_1(j)=&(-1)^kq^{2d_1d_2+2c_1c_2+2c_1n+2c_2n+2c_2+4n+k(d_2+c_2-k+1)}\tag{C.2}\\
        &\frac{(1-q^2)^2(1-q^4)(q^2;q^2)_{d_2}(q^2;q^2)_{c_2}}{(q^2;q^2)_{d_1+d_2+1}(q^2;q^2)_{c_1+c_2+n+2}}\frac{(-1)^{j}q^{j^2-j}{d_1\choose j}_{q^{2}}{n\choose j}_{q^2}(q^2;q^2)_{j}}{(q^{2d_1+2d_2+2c_1+2n-2j+4};q^2)_{c_2+2}}.\notag
\end{align}
Taking out the common factor $cf_1(j)$ from $p_2(j)+p_3(j)$, we have
\begin{align}\label{eq:b.3}\stepcounter{equation}
        &(p_2(j)+p_3(j))/cf_1(j)\tag{C.3}\\
        =&(1-q^{2d_1+2d_2+2c_1+2n-2j+4})\notag\\
        &\sum_{i=0}^{c_2-k}q^{(2d_1+2d_2+2n-2j+2)i}(q^2;q^2)_{c_1+1+i}(q^2;q^2)_{d_1+c_2+n-j-i}{c_2-k\choose i}_{q^2}\notag\\
        &+q^{2c_1+2}(1-q^{2c_2+2})\notag\\
        &\sum_{i=0}^{c_2+1-k}q^{(2d_1+2d_2+2n-2j+2)i}(q^2;q^2)_{c_1+i}(q^2;q^2)_{d_1+c_2+1+n-j-i}{c_2+1-k\choose i}_{q^2}\notag\\
        =&(1-q^{2d_1+2d_2+2c_1+2n-2j+4})\notag\\
        &\sum_{i=0}^{c_2-k}q^{(2d_1+2d_2+2n-2j+2)i}(q^2;q^2)_{c_1+1+i}(q^2;q^2)_{d_1+c_2+n-j-i}{c_2-k\choose i}_{q^2}\notag\\
        &+q^{2c_1+2}(1-q^{2c_2+2})\notag\\
        &\sum_{i=1}^{c_2+1-k}q^{(2d_1+2d_2+2n-2j+2)i}(q^2;q^2)_{c_1+i}(q^2;q^2)_{d_1+c_2+1+n-j-i}{c_2-k\choose i-1}_{q^2}\notag\\
        &+q^{2c_1+2}(1-q^{2c_2+2})\notag\\
        &\sum_{i=0}^{c_2-k}q^{(2d_1+2d_2+2n-2j+2)i}(q^2;q^2)_{c_1+i}(q^2;q^2)_{d_1+c_2+1+n-j-i}\cdot q^{2i}{c_2-k\choose i}_{q^2}\notag\\
        =&(1-q^{2d_1+2d_2+2c_1+2n-2j+4+2c_2+2})\notag\\
        &\sum_{i=0}^{c_2-k}q^{(2d_1+2d_2+2n-2j+2)i}(q^2;q^2)_{c_1+1+i}(q^2;q^2)_{d_1+c_2+n-j-i}{c_2-k\choose i}_{q^2}\notag\\
        &+q^{2c_1+2}(1-q^{2c_2+2})\notag\\
        &\sum_{i=0}^{c_2-k}q^{(2d_1+2d_2+2n+2-2j+2)i}(q^2;q^2)_{c_1+i}(q^2;q^2)_{d_1+c_2+1+n-j-i}\cdot {c_2-k\choose i}_{q^2}.\notag
\end{align}
Substituting (\ref{eq:b.3}) into (\ref{eq:b.1}) we get
\begin{align}\label{eq:b.4}\stepcounter{equation}
        &h((k^*)^{d_1+k}(h^*)^{d_2-k}a^{c_1}b^{c_2-k}c^{n+1+k}(c^*)^{n+1}(b^*)^{c_2}(a^*)^{c_1}h^{d_2}k^{d_1})\tag{C.4}\\
        =&(-1)^kq^{2d_1d_2+2c_1c_2+2c_1n+2c_2n+2c_2+4n+k(d_2+c_2-k+1)}\notag\\
        &\frac{(1-q^2)^2(1-q^4)(q^2;q^2)_{d_2}(q^2;q^2)_{c_2}}{(q^2;q^2)_{d_1+d_2+1}(q^2;q^2)_{c_1+c_2+n+1}}\left(\sum_{j=0}^{n}\frac{(-1)^{j}q^{j^2-j}{d_1\choose j}_{q^{2}}{n\choose j}_{q^2}(q^2;q^2)_{j}}{(q^{2d_1+2d_2+2c_1+2n-2j+4};q^2)_{c_2+1}}\right.\notag\\
        &\left.\cdot\underbrace{\sum_{i=0}^{c_2-k}q^{(2d_1+2d_2+2n-2j+2)i}(q^2;q^2)_{c_1+i}(q^2;q^2)_{d_1+c_2+n-j-i}{c_2-k\choose i}_{q^2}}_{p_1(j)}\right)\notag\\
        &-(-1)^kq^{2d_1d_2+2c_1c_2+2c_1n+2c_2n+2c_2+4n+k(d_2+c_2-k+1)}\notag\\
        &\frac{(1-q^2)^2(1-q^4)(q^2;q^2)_{d_2}(q^2;q^2)_{c_2}}{(q^2;q^2)_{d_1+d_2+1}(q^2;q^2)_{c_1+c_2+n+2}}\left(\sum_{j=0}^{n}\frac{(-1)^{j}q^{j^2-j}{d_1\choose j}_{q^{2}}{n\choose j}_{q^2}(q^2;q^2)_{j}}{(q^{2d_1+2d_2+2c_1+2n-2j+4};q^2)_{c_2+1}}\right.\notag\\
        &\left.\cdot\underbrace{\sum_{i=0}^{c_2-k}q^{(2d_1+2d_2+2n-2j+2)i}(q^2;q^2)_{c_1+1+i}(q^2;q^2)_{d_1+c_2+n-j-i}{c_2-k\choose i}_{q^2}}_{p_2'(j)}\right)\notag\\
        &-(-1)^kq^{2d_1d_2+2c_1c_2+2c_1n+2c_2n+2c_2+4n+k(d_2+c_2-k+1)+2c_1+2}\notag\\
        &\frac{(1-q^2)^2(1-q^4)(q^2;q^2)_{d_2}(q^2;q^2)_{c_2+1}}{(q^2;q^2)_{d_1+d_2+1}(q^2;q^2)_{c_1+c_2+n+2}}\left(\sum_{j=0}^{n}\frac{(-1)^{j}q^{j^2-j}{d_1\choose j}_{q^{2}}{n\choose j}_{q^2}(q^2;q^2)_{j}}{(q^{2d_1+2d_2+2c_1+2n-2j+4};q^2)_{c_2+2}}\right.\notag\\
        &\left.\cdot\underbrace{\sum_{i=0}^{c_2-k}q^{(2d_1+2d_2+2n-2j+4)i}(q^2;q^2)_{c_1+i}(q^2;q^2)_{d_1+c_2+1+n-j-i}{c_2-k\choose i}_{q^2}}_{p_3'(j)}\right).\notag
\end{align}
Then, we compute $p_1(j)-p_2'(j)$ for a fixed index $j$. The common factor of $p_1(j)$ and $p_2'(j)$ is 
\begin{align}\stepcounter{equation}
        cf_2(j)=&(-1)^kq^{2d_1d_2+2c_1c_2+2c_1n+2c_2n+2c_2+4n+k(d_2+c_2-k+1)}\tag{C.5}\\
        &\frac{(1-q^2)^2(1-q^4)(q^2;q^2)_{d_2}(q^2;q^2)_{c_2}}{(q^2;q^2)_{d_1+d_2+1}(q^2;q^2)_{c_1+c_2+n+2}}\frac{(-1)^{j}q^{j^2-j}{d_1\choose j}_{q^{2}}{n\choose j}_{q^2}(q^2;q^2)_{j}}{(q^{2d_1+2d_2+2c_1+2n-2j+4};q^2)_{c_2+1}}.\notag
\end{align}
Taking out the common factor $cf_2(j)$ from $p_1(j)-p_2'(j)$, we have
\begin{align}\label{eq:b.6}\stepcounter{equation}
        &(p_1(j)-p_2'(j))/cf_2(j)\tag{C.6}\\
        =&(1-q^{2c_1+2c_2+2n+4})\notag\\
        &\sum_{i=0}^{c_2-k}q^{(2d_1+2d_2+2n-2j+2)i}(q^2;q^2)_{c_1+i}(q^2;q^2)_{d_1+c_2+n-j-i}{c_2-k\choose i}_{q^2}\notag\\
        &-\sum_{i=0}^{c_2-k}q^{(2d_1+2d_2+2n-2j+2)i}(q^2;q^2)_{c_1+1+i}(q^2;q^2)_{d_1+c_2+n-j-i}{c_2-k\choose i}_{q^2}\notag\\
        =&-q^{2c_1+2c_2+2n+4}\notag\\
        &\sum_{i=0}^{c_2-k}q^{(2d_1+2d_2+2n-2j+2)i}(q^2;q^2)_{c_1+i}(q^2;q^2)_{d_1+c_2+n-j-i}{c_2-k\choose i}_{q^2}\notag\\
        &+q^{2c_1+2}\sum_{i=0}^{c_2-k}q^{(2d_1+2d_2+2n-2j+4)i}(q^2;q^2)_{c_1+i}(q^2;q^2)_{d_1+c_2+n-j-i}{c_2-k\choose i}_{q^2}\notag\\
        =&q^{2c_1+2c_2+2n+4}(q^{2d_1-2j}-1)\notag\\
        &\sum_{i=0}^{c_2-k}q^{(2d_1+2d_2+2n-2j+2)i}(q^2;q^2)_{c_1+i}(q^2;q^2)_{d_1+c_2+n-j-i}{c_2-k\choose i}_{q^2}\notag\\
        &+q^{2c_1+2}\underbrace{\sum_{i=0}^{c_2-k}q^{(2d_1+2d_2+2n-2j+4)i}(q^2;q^2)_{c_1+i}(q^2;q^2)_{d_1+c_2+n+1-j-i}{c_2-k\choose i}_{q^2}}_{p_3'(j)}.\notag
\end{align}
Substituting (\ref{eq:b.6}) into (\ref{eq:b.4}) and after simplification, we get
\begin{align}\label{eq:b.7}\stepcounter{equation}
        &h((k^*)^{d_1+k}(h^*)^{d_2-k}a^{c_1}b^{c_2-k}c^{n+1+k}(c^*)^{n+1}(b^*)^{c_2}(a^*)^{c_1}h^{d_2}k^{d_1})\tag{C.7}\\
        =&(-1)^kq^{2d_1d_2+2c_1c_2+2c_1(n+1)+2c_2(n+1)+2c_2+4(n+1)+k(d_2+c_2-k+1)}\notag\\
        &\frac{(1-q^2)^2(1-q^4)(q^2;q^2)_{d_2}(q^2;q^2)_{c_2}}{(q^2;q^2)_{d_1+d_2+1}(q^2;q^2)_{c_1+c_2+n+2}}\left(\sum_{j=0}^{n}\frac{(-1)^{j}q^{j^2-j}{d_1\choose j}_{q^{2}}{n\choose j}_{q^2}(q^2;q^2)_{j}q^{2n}(q^{2d_1-2j}-1)}{(q^{2d_1+2d_2+2c_1+2n-2j+4};q^2)_{c_2+1}}\right.\notag\\
        &\left.\cdot\sum_{i=0}^{c_2-k}q^{(2d_1+2d_2+2n-2j+2)i}(q^2;q^2)_{c_1+i}(q^2;q^2)_{d_1+c_2+n-j-i}{c_2-k\choose i}_{q^2}\right)\notag\\
        &+(-1)^kq^{2d_1d_2+2c_1c_2+2c_1(n+1)+2c_2(n+1)+2c_2+4(n+1)+k(d_2+c_2-k+1)}\notag\\
        &\frac{(1-q^2)^2(1-q^4)(q^2;q^2)_{d_2}(q^2;q^2)_{c_2}}{(q^2;q^2)_{d_1+d_2+1}(q^2;q^2)_{c_1+c_2+n+2}}\left(\sum_{j=0}^{n}\frac{(-1)^{j}q^{j^2-j}{d_1\choose j}_{q^{2}}{n\choose j}_{q^2}(q^2;q^2)_{j}}{(q^{2d_1+2d_2+2c_1+2n-2j+6};q^2)_{c_2+1}}\right.\notag\\
        &\left.\cdot\sum_{i=0}^{c_2-k}q^{(2d_1+2d_2+2n-2j+4)i}(q^2;q^2)_{c_1+i}(q^2;q^2)_{d_1+c_2+n+1-j-i}{c_2-k\choose i}_{q^2}\right).\notag
\end{align}
Notice that in the first term
\begin{equation}\label{eq:b.8}
    \begin{split}
        &(-1)^{j}q^{j^2-j}{d_1\choose j}_{q^{2}}(q^2;q^2)_{j}(q^{2d_1-2j}-1)\\
        =&(-1)^{j+1}q^{(j+1)^2-(j+1)}q^{-2j}{d_1\choose j+1}_{q^{2}}(q^2;q^2)_{j+1}.
    \end{split}
\end{equation}
Substituting (\ref{eq:b.8}) into (\ref{eq:b.7}), we get
\begin{align}\stepcounter{equation}
        &h((k^*)^{d_1+k}(h^*)^{d_2-k}a^{c_1}b^{c_2-k}c^{n+1+k}(c^*)^{n+1}(b^*)^{c_2}(a^*)^{c_1}h^{d_2}k^{d_1})\tag{C.9}\\
        =&(-1)^kq^{2d_1d_2+2c_1c_2+2c_1(n+1)+2c_2(n+1)+2c_2+4(n+1)+k(d_2+c_2-k+1)}\notag\\
        &\frac{(1-q^2)^2(1-q^4)(q^2;q^2)_{d_2}(q^2;q^2)_{c_2}}{(q^2;q^2)_{d_1+d_2+1}(q^2;q^2)_{c_1+c_2+n+2}}\notag\\
        &\left(\sum_{j=0}^{n}\frac{(-1)^{j}q^{j^2-j}{d_1\choose j}_{q^{2}}\left((q^{2n-2j+2}{n\choose j-1}_{q^2}+{n\choose j}_{q^2}\right)(q^2;q^2)_{j}}{(q^{2d_1+2d_2+2c_1+2n-2j+6};q^2)_{c_2+1}}\right.\notag\\
        &\left.\cdot\sum_{i=0}^{c_2-k}q^{(2d_1+2d_2+2n-2j+4)i}(q^2;q^2)_{c_1+i}(q^2;q^2)_{d_1+c_2+n+1-j-i}{c_2-k\choose i}_{q^2}\right)\notag\\
        =&(-1)^kq^{2d_1d_2+2c_1c_2+2c_1(n+1)+2c_2(n+1)+2c_2+4(n+1)+k(d_2+c_2-k+1)}\notag\\
        &\frac{(1-q^2)^2(1-q^4)(q^2;q^2)_{d_2}(q^2;q^2)_{c_2}}{(q^2;q^2)_{d_1+d_2+1}(q^2;q^2)_{c_1+c_2+n+2}}\left(\sum_{j=0}^{n}\frac{(-1)^{j}q^{j^2-j}{d_1\choose j}_{q^{2}}{n+1\choose j}_{q^2}(q^2;q^2)_{j}}{(q^{2d_1+2d_2+2c_1+2n-2j+6};q^2)_{c_2+1}}\right.\notag\\
        &\left.\cdot\sum_{i=0}^{c_2-k}q^{(2d_1+2d_2+2n-2j+4)i}(q^2;q^2)_{c_1+i}(q^2;q^2)_{d_1+c_2+n+1-j-i}{c_2-k\choose i}_{q^2}\right).\notag
\end{align}
This is consistent with (\ref{eq:5.35+}) when we put $c_3=n+1$. Hence, (\ref{eq:5.35+}) holds in the case $c_3=n+1$ and the Proposition is proved by induction.

%% file: main.bbl
\begin{thebibliography}{10}

\bibitem{andrews1998theory}
George~E Andrews.
\newblock {\em The theory of partitions}, volume~2.
\newblock Cambridge university press, 1998.

\bibitem{banica2007integration}
Teodor Banica and Beno{\^i}t Collins.
\newblock Integration over compact quantum groups.
\newblock {\em Publications of the Research Institute for Mathematical Sciences}, 43(2):277--302, 2007.

\bibitem{banica2007permutationintegration}
Teodor Banica and Beno{\^i}t Collins.
\newblock Integration over quantum permutation groups.
\newblock {\em Journal of Functional Analysis}, 242(2):641--657, 2007.

\bibitem{banica2008integration}
Teodor Banica and Beno{\^i}t Collins.
\newblock {Integration over the Pauli quantum group}.
\newblock {\em Journal of Geometry and Physics}, 58(8):942--961, 2008.

\bibitem{bergeron2019suq}
Geoffroy Bergeron, Erik Koelink, and Luc Vinet.
\newblock {$SU_q(3)$ corepresentations and bivariate q-Krawtchouk polynomials}.
\newblock {\em Journal of Mathematical Physics}, 60(5), 2019.

\bibitem{brezin1980external}
Edouard Br{\'e}zin and David~J Gross.
\newblock {The external field problem in the large N limit of QCD}.
\newblock {\em Physics Letters B}, 97(1):120--124, 1980.

\bibitem{brouwer1996diagrammatic}
PW~Brouwer and CWJ Beenakker.
\newblock Diagrammatic method of integration over the unitary group, with applications to quantum transport in mesoscopic systems.
\newblock {\em arXiv preprint cond-mat/9604059}, 1996.

\bibitem{collins2003moments}
Beno{\^i}t Collins.
\newblock Moments and cumulants of polynomial random variables on unitarygroups, the itzykson-zuber integral, and free probability.
\newblock {\em International Mathematics Research Notices}, 2003(17):953--982, 2003.

\bibitem{collins2022moment}
Beno{\^i}t Collins.
\newblock Moment methods on compact groups: Weingarten calculus and its applications.
\newblock {\em ICM Vol. IV. Sections 5}, 8:3142--3164, 2022.

\bibitem{collins2023tensor}
Beno{\^i}t Collins, Razvan Gurau, and Luca Lionni.
\newblock {The tensor Harish-Chandra--Itzykson--Zuber integral I: Weingarten calculus and a generalization of monotone Hurwitz numbers}.
\newblock {\em Journal of the European Mathematical Society}, 26(5):1851--1897, 2023.

\bibitem{collins2022weingarten}
Beno{\^i}t Collins, Sho Matsumoto, and Jonathan Novak.
\newblock The weingarten calculus.
\newblock {\em Notices of the American Mathematical Society}, 69(05):1, 2022.

\bibitem{collins2006integration}
Beno{\^i}t Collins and Piotr {\'S}niady.
\newblock Integration with respect to the haar measure on unitary, orthogonal and symplectic group.
\newblock {\em Communications in Mathematical Physics}, 264(3):773--795, 2006.

\bibitem{creutz1978invariant}
Michael Creutz.
\newblock {On invariant integration over SU(N)}.
\newblock {\em Journal of Mathematical Physics}, 19(10):2043--2046, 1978.

\bibitem{diaconis1994eigenvalues}
Persi Diaconis and Mehrdad Shahshahani.
\newblock On the eigenvalues of random matrices.
\newblock {\em Journal of Applied Probability}, 31(A):49--62, 1994.

\bibitem{dijkhuizen1994cqg}
Mathijs~S Dijkhuizen and Tom~H Koornwinder.
\newblock {CQG algebras: a direct algebraic approach to compact quantum groups}.
\newblock {\em Letters in Mathematical Physics}, 32(4):315--330, 1994.

\bibitem{drinfeld1986quantum}
Vladimir~Gershonovich Drinfeld.
\newblock Quantum groups.
\newblock {\em Proc. Int. Congr. Math.}, 1:798--820, 1986.

\bibitem{faddeev1988quantization}
Ludwig~D Faddeev, N~Yu Reshetikhin, and Leon~A Takhtajan.
\newblock Quantization of lie groups and lie algebras.
\newblock In {\em Algebraic analysis}, pages 129--139. Elsevier, 1988.

\bibitem{jimbo1985aq}
Michio Jimbo.
\newblock {A $q$-difference analogue of $U_q(g)$ and the Yang-Baxter equation}.
\newblock {\em Letters in Mathematical Physics}, 10(1):63--69, 1985.

\bibitem{klimyk2012quantum}
Anatoli Klimyk and Konrad Schm{\"u}dgen.
\newblock {\em Quantum groups and their representations}.
\newblock Springer Science \& Business Media, 2012.

\bibitem{koelink1989clebsch}
Hendrik~Tjerk Koelink and Tom~H Koornwinder.
\newblock {The Clebsch-Gordan coefficients for the quantum group S$_\mu$U(2) and $q$-Hahn polynomials}.
\newblock In {\em Indagationes Mathematicae (Proceedings)}, volume~92, pages 443--456. Elsevier, 1989.

\bibitem{kumar2015connections}
Shrawan Kumar and Joseph~M Landsberg.
\newblock {Connections between conjectures of Alon--Tarsi, Hadamard--Howe, and integrals over the special unitary group}.
\newblock {\em Discrete Mathematics}, 338(7):1232--1238, 2015.

\bibitem{magee2019matrix}
Michael Magee and Doron Puder.
\newblock Matrix group integrals, surfaces, and mapping class groups i: U(n).
\newblock {\em Inventiones mathematicae}, 218(2):341--411, 2019.

\bibitem{masuda1991representations}
Tetsuya Masuda, Katsuhisa Mimachi, Yoshiomi Nakagami, Masatoshi Noumi, and Kimio Ueno.
\newblock {Representations of the quantum group $SU_q(2)$ and the little $q$-Jacobi polynomials}.
\newblock {\em Journal of Functional Analysis}, 99(2):357--386, 1991.

\bibitem{noumi1993finite}
Masatoshi Noumi, Hirofumi Yamada, and Katsuhisa Mimachi.
\newblock {Finite Dimensional representations of the quantum group $GL_q(n;C)$ and the zonal spherical functions on $U_q(n-1) \setminus U_q(n)$}.
\newblock {\em Japanese journal of mathematics. New series}, 19(1):31--80, 1993.

\bibitem{odama2001enumeration}
Yumi Odama and Gregg Musiker.
\newblock Enumeration of (0, 1) and integer doubly stochastic matrices.
\newblock {\em Science Direct Working Paper}, (S1574-0358):04, 2001.

\bibitem{paule1997mathematica}
Peter Paule and Axel Riese.
\newblock {A Mathematica $q$-analogue of Zeilberger’s algorithm based on an algebraically motivated approach to $q$-hypergeometric telescoping}.
\newblock {\em in Special Functions, q-Series and Related Topics, Fields Inst. Commun.}, 14:179--210, 1997.

\bibitem{reshetikhin2001quantum}
Nicolai Reshetikhin and Milen Yakimov.
\newblock Quantum invariant measures.
\newblock {\em Communications in Mathematical Physics}, 224(2):399--426, 2001.

\bibitem{riese1997contributions}
Axel Riese.
\newblock {\em {Contributions to symbolic $q$-hypergeometric summation}}.
\newblock PhD Dissertation, 1997.

\bibitem{riese2003qmultisum}
Axel Riese.
\newblock {qMultiSum—a package for proving q-hypergeometric multiple summation identities}.
\newblock {\em Journal of Symbolic Computation}, 35(3):349--376, 2003.

\bibitem{sl1987twisted}
Woronowicz SL.
\newblock {Twisted SU(2) group. An example of a non-commutative differential calculus}.
\newblock {\em Publications of the Research Institute for Mathematical Sciences}, 23(1):117--181, 1987.

\bibitem{stein1970enumeration}
ML~Stein and PR~Stein.
\newblock Enumeration of stochastic matrices with integer elements.
\newblock Technical report, Los Alamos National Lab.(LANL), Los Alamos, NM (United States), 1970.

\bibitem{sweedler1969hopf}
M.E. Sweedler.
\newblock {\em Hopf Algebras}.
\newblock Mathematics lecture note series. W. A. Benjamin, 1969.

\bibitem{vaksman1990algebra}
Leonid~L'vovich Vaksman and Yakov~Sergeevich Soibel'man.
\newblock {Algebra of functions on the quantum group SU(n+1) and odd-dimensional quantum spheres}.
\newblock {\em Algebra i analiz}, 2(5):101--120, 1990.

\bibitem{wang1995free}
Shuzhou Wang.
\newblock Free products of compact quantum groups.
\newblock {\em Communications in Mathematical Physics}, 167(3):671--692, 1995.

\bibitem{weingarten1978asymptotic}
Don Weingarten.
\newblock Asymptotic behavior of group integrals in the limit of infinite rank.
\newblock {\em J. Math. Phys.(NY);(United States)}, 19(5), 1978.

\bibitem{woronowicz1987compact}
Stanis{\l}aw~L Woronowicz.
\newblock Compact matrix pseudogroups.
\newblock {\em Communications in Mathematical Physics}, 111(4):613--665, 1987.

\bibitem{woronowicz1988tannaka}
Stanis{\l}aw~L Woronowicz.
\newblock {Tannaka-Krein duality for compact matrix pseudogroups. Twisted SU(N) groups}.
\newblock {\em Inventiones mathematicae}, 93(1):35--76, 1988.

\end{thebibliography}
